\documentclass[12pt]{amsart}
\usepackage{amssymb}
\usepackage[all]{xy}
\newcommand\Diff{\operatorname{Diff}}

\newcommand\Ti{{\mathcal T}^\infty}

\def \crpr {\rtimes}
\def  \nuint {\raise10pt\hbox{$\nu$}\kern-6pt\int}
\def \OM {\overline{\widehat{\Omega}}_\ast(T, \Bi_\Gamma)}
\newcommand\Hom{\operatorname{Hom}}
\newcommand\wn{\widehat{N}}
\newcommand\wf{\widehat{F}}
\newcommand\STR{\operatorname{STR}}
\newcommand\Str{\operatorname{Str}}
\newcommand\Tr{\operatorname{Tr}}
\newcommand\ind{\operatorname{ind}}
\newcommand\Ind{\operatorname{Ind}}
\newcommand\ch{\operatorname{Ch}}

\def \AA {\mathbb A}

\newcommand\wm{\widehat{M}}
\newcommand\we{\widehat{E}}

\newcommand\R{\mathcal R}

\newcommand\B{\mathcal B}
\newcommand\Bi{\B^\infty}
\newcommand\D{\mathcal D}
\newcommand\Di{D\kern-6pt/}
\newcommand\cDi{{\mathcal D}\kern-6pt/}
\newcommand\CC{\mathbb C}
\def \n {\noindent}
\def \s {\smallskip}

\def \m {\medskip}
\def \cal {\mathcal}
\def \BB {\mathbb B}

\newcommand\NN{\mathbb N}

\newcommand\RR{\mathbb R}
\newcommand\ZZ{\mathbb Z}
\newcommand\KK{\mathbb K}
\newcommand\HH{\mathbb H}

\newcommand\pa{\partial}

\newcommand\Id{\operatorname{Id}}

\newtheorem{theorem}{Theorem}

\newtheorem{proposition}{Proposition}
\newtheorem{corollary}{Corollary}
\newtheorem{lemma}{Lemma}

\newtheorem{definition}{Definition}

\theoremstyle{remark}
\newtheorem{remark}{Remark}
\begin{document}

\title[ ]
{Etale Groupoids, eta invariants and index theory}
\author[E. Leichtnam]{Eric Leichtnam}
\address{Institut de Jussieu et CNRS,
Etage 7E,
175 rue du Chevaleret,
75013, Paris,  France}
        \email{leicht@math.jussieu.fr}
\author[P. Piazza]{Paolo Piazza}
        \address{Dipartimento di Matematica G. Castelnuovo,
Universit\`a di Roma ``La Sapienza'', P.le Aldo Moro 2, 00185 Rome, Italy}
       \email{piazza@mat.uniroma1.it}


\begin{abstract}
Let $\Gamma$ be a discrete finitely generated group.
Let
$\widehat{M}\rightarrow T$
be  a $\Gamma$-equivariant fibration,
with fibers diffeomorphic to
a fixed even dimensional manifold with boundary $Z$. We assume that
$\Gamma\rightarrow \wm\rightarrow \wm/\Gamma$ is
a Galois covering of a compact manifold with boundary.
Let $(D^+ (\theta))_{\theta\in T}$ be a $\Gamma$-equivariant family
of Dirac-type operators. Under the assumption that
 the boundary family is
$L^2$-invertible , we define an index class
in $K_0 (C^0 (T)\rtimes_r \Gamma)$. If, in addition,
$\Gamma$ is of polynomial growth,  we define higher indeces
by pairing the index class with suitable cyclic cocycles.
Our main result is then  a   formula
for these higher indeces: the structure of the formula is
as in the seminal work of Atiyah, Patodi and Singer, with an
interior geometric contribution and a boundary contribution
in the form of a higher eta invariant associated to the boundary family.
Under similar assumptions we extend our theorem to any $G$-proper manifold,
with $G$ an \'etale groupoid. We employ this generalization in order
to establish a higher Atiyah-Patodi-Singer index formula on certain
foliations with boundary.
Fundamental to our work is  a suitable generalization of Melrose
$b$-pseudodifferential calculus as well as
 the superconnection proof of the index theorem on $G$-proper manifolds
recently given by  Gorokhovsky and Lott in \cite{Go-Lo}.

\end{abstract}

\maketitle \tableofcontents

\section{{\bf Introduction and main results}}
Connes' index theorem for $G$-proper manifolds \cite{Co}, with $G$
an \'etale groupoid, unifies under a single statement most of the
existing index theorems. For this introduction we shall focus on a
particular case of such a theorem: $\Gamma$ is a discrete group
acting on an even dimensional  manifold  $\widehat{M}$ and on a
compact manifold $T$; $\widehat{M}\rightarrow T$ is a
$\Gamma$-equivariant fibration and  the action of $\Gamma$ on
$\widehat{M}$ is free, properly discontinuous and cocompact, thus
 $\widehat{M}/\Gamma:=M$ is a smooth compact manifold and  $\Gamma\rightarrow
 \widehat{M}\rightarrow M$
is a $\Gamma$-Galois covering. (This is an example of $G$-proper
manifold with $G$ equal to the groupoid $T\rtimes \Gamma$.)

\smallskip
 If
$T=$ point and $\Gamma=\{1\}$ we have a compact manifold and
Connes' index theorem reduces to the Atiyah-Singer index theorem.

\smallskip
 If $\Gamma=\{1\}$ we simply have a fibration and the theorem
 reduces to the Atiyah-Singer family index theorem.

\smallskip
 Finally, if $T=$ point then we have a Galois
covering and Connes' index theorem reduces to the Connes-Moscovici
higher index theorem.

\smallskip
 If $\dim T >0$ and $\Gamma\not= \{1\}$, then
Connes' index theorem can be seen as a higher foliation index
theorem on the foliated manifold, $(M,\mathcal{F})$, obtained by
foliating $M$ by the images of the fibers of
$\widehat{M}\rightarrow T$. Recall that
Connes' index theorem can be seen
as a cohomological version of a K-theoretical statement, namely
the Connes-Skandalis longitudinal
index theorem for foliations \cite{CSk}.

\smallskip
In the case  $\Gamma=\{1\}$, $T=$ point (compact closed
manifolds), in the case $\Gamma=\{1\}$ (fibrations) and in the
case $T=$ point (Galois coverings), the  index theorems mentioned
above have been sharpened into {\it local index theorems.} For a
single closed compact manifold this goes back to the work of
McKean-Singer, Patodi, Gilkey and many others. In the case of a
fibration the result is due to Bismut, see \cite{Bismut}
\cite{BGV}, whereas in the  case of Galois coverings the local
version of the higher index theorem it is due to Lott, see
\cite{Lott I}. Notice that for fibrations and Galois coverings the
new proofs make use of the heat kernel associated to a
superconnection. One of the most interesting outcome of these
improved proofs is the possibility of extending the index theorem
to manifolds with boundary. In the case of a single manifold and
of a Dirac-type operator on it, such an index theorem is due to
Atiyah-Patodi-Singer. For families of Dirac operators on manifolds
with boundary, the theorem is due to Bismut and Cheeger and, more
generally, to Melrose and Piazza (\cite{BC},\cite{BC2},\cite{MP
I}, \cite{MP II}). Finally, for Galois coverings of a compact
manifold with boundary, the result is due to the authors of the
present article, see \cite{LPMEMOIRS}, \cite{LPGAFA}, following a
conjecture of Lott \cite{Lott II}.\\ {\it Geometric applications} of these
results 
have been given
\begin{itemize}
\item to the problem of
defining higher signatures on manifolds with boundary \cite{Lott
II} \cite{Lott III} and proving their homotopy invariance \cite{LLP},
\item  to
uniqueness problems in positive scalar curvature \cite{LPPSC},
\item  to
the problem of cut-and-paste invariance of Novikov higher
signatures on closed manifolds \cite{LLP} (see also \cite{LLK}
\cite{Hilsum} \cite{LPCUT}),
\item to the homotopy invariance of the Atiyah-Patodi-Singer
and Cheeger-Gromov rho-invariants for closed compact manifolds
having a  torsion-free fundamental group $\Gamma$
satisfying the bijectivity
of the Baum-Connes map for $C^*_{{\rm max}} \Gamma$
\cite{PS}, a result due originally to Keswani \cite{Kes}.
\end{itemize}

\smallskip
Recently  Gorokhovsky and Lott have given a
superconnection-heat-kernel  proof of Connes' index theorem, see
\cite{Go-Lo}, and raised the question of extending such a theorem
to manifolds with boundary. {\it The goal of this paper is to
establish such a result,  thus proving a higher
Atiyah-Patodi-Singer higher index theorem on a $G$-proper manifold
with boundary.} Geometric applications of our theorem will be
considered in a future publication.

\bigskip We shall now give a detailed description of the content of
this paper. In  Section \ref{fibration}  we  introduce our basic
geometric object: a $\Gamma$-equivariant fibration
$\widehat{M}\rightarrow T$ with $\widehat{M}$ a manifold with
boundary, $T$ a closed compact manifold and  $\Gamma\rightarrow
 \widehat{M}\rightarrow \widehat{M}/\Gamma$
 a $\Gamma$-Galois covering. This is done in Subsection \ref{data}
 In Subsection \ref{groupoid} we introduce the groupoid $T\rtimes \Gamma$
 and the associated $C^*$-algebra $C^*_r (T\rtimes \Gamma)$.
The latter is simply the reduced cross-product algebra $C^0 (T)
\rtimes_r \Gamma$ and can be described as a suitable completion of
the algebraic cross-product  $$C^\infty (T)\rtimes
\Gamma=\{\sum_{{\rm finite}} f_g (\theta) g \} \;\;\text{with}
\;\;g \cdot (f_h (\theta) h):= f_h (\theta \cdot g) gh\,.$$

 We next explain how
to  associate to our $\Gamma$-equivariant fibration
$\widehat{M}\rightarrow T$ some natural Sobolev $(C^0 (T)\rtimes_r
\Gamma)$-Hilbert modules $H^m_{b,C^0(T)\rtimes_r \Gamma}$; we also
introduce $\Gamma$-equivariant families  of Dirac-type operators
and explain how  such a $\Gamma$-equivariant family $\mathcal{D}$
defines  a morphism of Sobolev $(C^0 (T)\rtimes_r \Gamma)$-Hilbert
modules $H^m_{b,C^0(T)\rtimes_r \Gamma}\to
H^{m-1}_{b,C^0(T)\rtimes_r \Gamma}$; this is done in Subsections
\ref{g-families} and \ref{h-modules}.
A $\Gamma$-equivariant family of Dirac operators
defines in a natural way a longitudinal  elliptic operator on the
foliated manifold  $(M,\mathcal{F})$; for the sake of completeness
we explain in Subsection \ref{foliations}  the connection between the
noncommutative framework just explained and the one associated to
such a longitudinal operator on $(M,\mathcal{F})$. We end Section
\ref{data} introducing our basic assumption on the family of
operators, denoted $\mathcal{D}_0$, induced on the boundary by a
$\Gamma$-equivariant family of Dirac operators $\mathcal{D}$; this
is an invertibility assumption. We also draw some important
consequences out of this assumption for the morphism of $(C^0 (T)
\rtimes_r \Gamma)$-Hilbert modules defined by $\mathcal{D}_0$.

\smallskip
In Section \ref{bcalc} we develop  a $\Gamma$-equivariant
fiber-$b$-calculus on $\widehat{M}\to T$ and we employ it in order
to prove that, under our invertibility assumption, a
$\Gamma$-equivariant family of Dirac operators defines an index
class in $K_* (C^0 (T)\rtimes_r \Gamma)$.

\smallskip
Once the index class is defined we need a way to extract numerical
invariants out of it; these are  the {\it higher indeces}. Our
main concern  is then to prove a Atiyah-Patodi-Singer-type formula
for these higher indeces,  equating a higher index with the sum of
an interior geometric contribution and a boundary contribution.

In order to orient the   reader, we recall the structure of these
index theorems in the three special cases considered above. To fix
the notation we assume that we are in the spin context and that
our operators are Dirac operators acting on spinors.

\smallskip
 If
$T=$ point and $\Gamma=\{e\}$, then our $(T\rtimes \Gamma)$-proper
manifold is simply equal to an even dimensional compact manifold
with boundary $M$. In this case the index of the
Atiyah-Patodi-Singer boundary value problem associated to the
Dirac operator $D$ is expressed as follows: $$\ind_{{\rm APS}}
(D)= \int_M \widehat{A}(TM,\nabla^{{\rm LC}})
-\frac{\eta(D_0)+\dim {\rm Ker}(D_0)}{2}$$ with $\nabla^{{\rm LC}}
$ equal to the Levi-Civita
connection and  
$\eta (D_0)$ the eta invariant associated to the boundary operator
$D_0$.

\smallskip
If $\Gamma=\{e\}$ then the $(T\rtimes \Gamma)$-proper manifold  we
are dealing with is  a fibration $$Z\to \widehat{M}\to T$$ with
fibers diffeomeorphic to a compact spin manifold with boundary
$Z$. If $\mathcal{D}=(D(\theta))_{\theta\in T}$ is a family of
Dirac operators with invertible boundary family
$\mathcal{D}_0\equiv ((D(\theta)_0)_{\theta\in T}$ then there
exists a well defined index class $\Ind_{{\rm APS}}
(\mathcal{D})\in K^0 (T)=K_0 (C^0 (T))$. We extract higher indeces
out of this index class by taking the Chern character $\ch
(\Ind_{{\rm APS}} (\mathcal{D}))\in H^*_{dR} (M)$ and coupling it
with the homology class $[\Phi]\in H_* (M,\RR)$ defined by a
closed current $\Phi$ . The index formula can already be given at
the level of Chern character and reads: $$ \ch (\Ind_{{\rm APS}}
(\mathcal{D}))=\int_{{\rm fiber}} \widehat{A}(T_V(\widehat{M}),
\nabla^V) -\frac{1}{2} \widetilde{\eta}
(\mathcal{D}_0)\;\;\;\text{in}\;\;\; H^*_{dR}(T) $$ with
$\widetilde{\eta}(\mathcal{D}_0) \in \Omega^* (T):= C^\infty
(T,\Lambda^* T)$ the Bismut-Cheeger eta form defined by the
boundary family, $T_V (\widehat{M})$ the vertical tangent bundle
to the fibration and $\nabla^V$ a certain connection on it.
 Notice
that the right hand side of the index formula involves {\it
smooth} differential forms.

\smallskip
If $T=$ point, then our  $(T\rtimes \Gamma)$-proper manifold is a
Galois covering $\Gamma\to \widehat{M} \to M$ of a compact
manifold with boundary $M$ and the analytic object of interest to
us is a $\Gamma$-invariant Dirac operator on $\widehat{M}$; we
keep denoting this operator by $\mathcal{D}$. Under an
invertibility assumption on the boundary operator $\mathcal{D}_0$
there is a well defined index class $\Ind_{{\rm APS}}
(\mathcal{D})\in K_0 (C^*_r \Gamma)$. The reduced group
$C^*$-algebra plays the role here of the continuous functions on
the base of the fibration. As explained in \cite{Lott I}
\cite{Lott II} this analogy can be pushed much further: thus there
exist a "smooth" subalgebra $\Bi_\Gamma$, $$\CC\Gamma\subset
\Bi_\Gamma\subset C^*_r \Gamma\,,$$ playing the role of the
$C^\infty$ function on the base of the fibration, and an algebra
of "smooth" noncommutative differential forms $\widehat{\Omega}_*
(\Bi_\Gamma)$. Under our invertibility assumption there is a well
defined higher eta invariant $\widetilde{\eta}(\mathcal{D}_0)$
associated to the boundary family $\mathcal{D}_0$, an element in
$\widehat{\Omega}_*(\Bi_\Gamma)$ modulo the closure of the graded commutator
$[\widehat{\Omega}_*(\Bi_\Gamma),\widehat{\Omega}_*(\Bi_\Gamma)]$.
The higher eta
invariant was defined by Lott in \cite{Lott I} under the
assumption that $\Gamma$ is of {\it polynomial growth}. For the
general case see  \cite{Lott III}  and  the Appendix in
\cite{LPBSMF} . The "smooth" algebra $\Bi_\Gamma$ is in fact
dense and holomorphically closed in $C^*_r \Gamma$, so that $K_*
(\Bi_\Gamma)=K_* (C^*_r \Gamma)$.
 The graded differential algebra  $\widehat{\Omega}_* (\Bi_\Gamma)$ defines
noncommutative de Rham homology groups $\widehat{H}_*
(\Bi_\Gamma)$ and there is a well defined Chern character $ \ch
:K_* (\Bi_\Gamma)=K_* (C^*_r \Gamma)\rightarrow \widehat{H}_*
(\Bi_\Gamma)$, see \cite{Karoubi}. One defines higher indeces for
$\mathcal{D}$ by pairing the Chern character of the index class,
$\ch(\Ind_{{\rm APS}} (\mathcal{D}))$, with an element in the
cohomology of $ \Bi_\Gamma$, i.e. with  a {\it closed graded trace
$\Phi$ on $\widehat{\Omega}_* (\Bi_\Gamma)$}. The result proved in
\cite{LPMEMOIRS} (see also \cite{LPBSMF} Appendix) gives a formula
for these higher indeces; the formula already holds at the level
of Chern character, as an equality in $\widehat{H}_*
(\Bi_\Gamma)$, and reads $$\ch(\Ind_{{\rm APS}} (\mathcal{D}))
=\left[  \int_M \widehat{A}(M,\nabla^{{\rm LC}})\wedge \omega -
\frac{1}{2} \widetilde{\eta}(\mathcal{D}_0)\right]
\quad\text{in}\quad\widehat{H}_* (\Bi_\Gamma)$$ with $\omega\in
\Omega^* (M) \widehat{\otimes} \Omega_*(\CC\Gamma)$ an explicit
bi-form, see \cite{Lott I}. In particular, if $\Phi$ is a closed
graded trace on $\Omega_* (\Bi_\Gamma)$ then the higher index
$\Ind_\Phi (\D):=<\ch(\Ind_{{\rm APS}} (\mathcal{D})) ,[\Phi]>$ is
well defined, $\omega_\Phi:=<\omega,\Phi>$ is a smooth
differential form on $M$, $<\widetilde{\eta}(\mathcal{D}_0),\Phi>$
is also well defined and the following formula holds:
\begin{equation}
\Ind_\Phi (\D)=\int_M \widehat{A}(M,\nabla^{{\rm LC}})\wedge
\omega_\Phi - \frac{1}{2}<\widetilde{\eta}(\mathcal{D}_0),\Phi>.
\end{equation}
\medskip

We now back to the general case addressed in this paper:
$\widehat{M}\to T$ is a $\Gamma$-equivariant fibration, $\Gamma\to
\widehat{M}\to M$ is a Galois covering of a manifold with boundary
and $\D=(D_\theta)_{\theta\in T}$ is a $\Gamma$-invariant family
of Dirac operators. On the basis of the last two examples it is
clear that in order to develop a higher Atiyah-Patodi-Singer index
theorem we need first of all to fix a "smooth" subalgebra $\Ti$,$$
C^\infty (T)\rtimes \Gamma\subset \Ti\subset C^0 (T)\rtimes_r
\Gamma\,,$$ together with an algebra of smooth noncommutative
differential forms. These two steps are in fact already necessary
in the closed case, in order to develop a {\it local} higher index
theory: we can thus follow closely the work of Gorokhovsky and
Lott, although, for technical reasons having to do both with the
convergence of the higher eta invariant and the construction
of a suitable rapidly-decreasing-parametrix, we need to assume that
$\Gamma$ is of polynomial growth.
(For more on this, see the last remark in subsection \ref{rapid2}.)
The "smooth" subalgebra $\Ti$ is
then defined in terms of infinite sums $\sum_{g\in \Gamma} f_g
(\theta) g$ but with a growth condition on the coefficients $f_g
(\theta)\in C^\infty (T)$: $$\sup_{\theta\in T, g\in \Gamma}
(|f_g| (1+\|g\|)^N) < \infty \,\,\,\forall N \,.$$ The
noncomutative differential forms $\widehat{\Omega}(T,\Bi_\Gamma)$ are
defined in terms of infinite sums $$\sum \alpha_{g_0, g_1,
\cdots,g_\ell} g_0 dg_1 \cdots dg_\ell\,,\quad  \alpha_{g_0, g_1,
\cdots,g_\ell} \in \Omega^k (T)$$ with a similar growth condition
on the differential forms $ \alpha_{g_0, g_1, \cdots,g_\ell}
\in \Omega^k (T) $ and
with the product taking into account the fact that $\Gamma$ acts
on $T$. All this is explained in Section \ref{heta}. The higher eta
invariant associated to the boundary family $\D_0$, denoted
$\widetilde{\eta}_{<e>} (\D_0)$, is  an element in $\widehat{\Omega}_{*,<e>}
(T,\Bi_\Gamma)$ (modulo graded commutators), the subalgebra of
elements $$\sum_{g_0 g_1 \cdots g_\ell=e} \alpha_{g_0, g_1,
\cdots,g_\ell} g_0 dg_1 \cdots dg_\ell$$ concentrated on the
trivial conjugacy class. The definition of $\widetilde{\eta}_{<e>}
(\D_0)$ employs the Gorokhovsky-Lott superconnection, the heat
kernel associated to it and a suitable supertrace; the latter is
defined on the space of $\Gamma$-invariant families of smoothing
operators with $\widehat{\Omega}_* (T,\Bi_\Gamma)$-coefficients and has
values  in $\Omega_{*,<e>} (T,\Bi_\Gamma)$ (modulo graded
commutators). Let $\Phi$ a closed graded trace on
$\widehat{\Omega}(T,\Bi_\Gamma)$; we assume that $\Phi$ is concentrated on
the trivial conjugacy class, i.e. that it is non-zero only on the
subalgebra  $\widehat{\Omega}_{*,<e>} (T,\Bi_\Gamma)$. Then the number
$<\widetilde{\eta}_{<e>} (\D_0),\Phi>$ is well-defined. One can
also define a higher index $\Ind_\Phi (\D)$. If $\Gamma$ acts by
isometries on $T$,  then $\Ti$ is stable under holomorphic
functional calculus and the definition of  $\Ind_\Phi (\D)$
proceeds in a way similar to Galois coverings, i.e. through a
noncommutative Chern character \`a la Karoubi. We explain all this
in Section \ref{isomcase}  and in  Appendix C we prove the following  higher
Atiyah-Patodi-Singer index formula:

\smallskip
\n
 {\it there exists an explicit {\it current} $\omega_\Phi$ on
$M:=\widehat{M}/\Gamma$ such that}
\begin{equation}\label{theformula}
\Ind_\Phi (\D)= <\widehat{A}(T\mathcal{F},\nabla^{{\rm
LC},\mathcal{F}}),\omega_\Phi>-\frac{1}{2}<\widetilde{\eta}_{<e>}
(\D_0),\Phi>\,.
\end{equation}
In the first summand on the right hand side the pairing between
differential forms and currents appears.

\smallskip
 If $\Gamma$ does not act by isometries then the
definition of higher index $\Ind_\Phi (\D)$ is more involved and
the proof of formula $\ref{theformula}$ much more complicated. The
definition is based on ideas of Connes (\cite{Co}, page 229),
Nistor \cite{Nistor}
and
Gorokhovsky-Lott; the proof is based on  arguments used by
Gorokhovsky-Lott in order to prove the local version of Connes'
index theorem and on $b$-calculus techniques developed in Sections
\ref{rapid} and \ref{bsuper} of the present paper. The proof of formula
\ref{theformula} occupies all of Section \ref{mainth}.

\smallskip
The next section of the paper, Section \ref{genetale},
is devoted to a generalization of
the above results to general $G$-proper manifolds, with $G$ an
\'etale groupoid satisfying a suitable polynomial growth
condition. We employ such a result in order to establish a higher
Atiyah-Patodi-Singer index formula on certain foliations with
boundary; this is done in Section \ref{applfoliation}.

\smallskip
Finally, in Appendix A we define the rapidly decreasing $b$-calculus, in Appendix B
we define the space of $b$-smoothing operators with differential form coefficients
and in Appendix C we give a proof of our theorem in the isometric case.

\medskip
\n {\bf Acknowledgements.}  It is  a pleasure to thank Sasha
Gorokhovsky, Vincent Lafforgue, Hitoshi Moriyoshi and George
Skandalis for useful explanations and discussions.\\ Part of this
work was done while the first author was visiting the university
of Savoie, he would like to thank K.Kurdyka for the kind
hospitality. The paper was completed while the second author was
visiting {\it Institut de Math\'ematiques de Jussieu} and  he
would like to thank the members  of the {\it \'Equipe d'Alg\`ebres
d'op\'erateurs et representations} for the stimulating working
conditions he enjoyed during his visit.

\smallskip
\n
Both authors are members of
 the RTN ``Geometric Analysis''
of the European Community. This research project was partially
supported by a CNR-CNRS cooperation project, by {\it Ministero
Istruzione Universit\`a e Ricerca} (Italy), by {\it Institut de
Math\'ematiques de Jussieu} and by the RTN ``Geometric Analysis''.
 We thanks these institutions
for their support.

\section{{\bf $\Gamma$-equivariant fibrations}}\label{fibration}

\subsection{Geometric data}\label{data}$\;$

\medskip
Let $\Gamma$ be a finitely
generated discrete group. Let $T$ be
a smooth closed compact connected manifold
on which
$\Gamma$ acts on the right. Let $\widehat{M}$ be
a manifold with boundary
on which $\Gamma$ acts freely, properly and
cocompactly on the right: the quotient space $M=\widehat{M}
/\Gamma$ is thus a smooth compact manifold with boundary.
We assume that $\widehat{M}$ fibers over $T$ and that
the resulting fibration
$$
\pi: \,\widehat{M} \rightarrow T
$$
is a {\it $\Gamma$-equivariant fibration} with  fibers
$\pi^{-1}(\theta),\theta \in T,$ that are
transverse to $\partial \widehat{M}$ and of dimension $2k$.
Notice that each fiber is a smooth manifold with boundary;
we shall also denote the typical fiber of
$\pi: \,\widehat{M} \rightarrow T $ by $Z$.
We choose a $\Gamma$-invariant exact $b$-metric
\cite{Melrose} on the vertical $b$-tangent
bundle ${}^b TZ $. Finally, we  assume
the existence of a
$\Gamma-$equivariant spin structure on ${}^b TZ $
that is fixed once and for all.
We denote  by
$S^Z \rightarrow \widehat{M} $ the associated spinor bundle.

\smallskip
The compact manifold with boundary $M$ inherits a foliation
${\cal F}$, with leaves equal to the image of the fibres
of $\pi: \,\widehat{M} \rightarrow T
$ under the quotient map $\widehat{M}\rightarrow M=\widehat{M}/\Gamma$.
Notice that the foliation ${\cal F}$ is transverse to the
boundary of $M$.

\medskip
\noindent {\bf Example.} Let $X$ be a compact manifold with
boundary and let $\Gamma\rightarrow \widetilde{X}\rightarrow X$ be a
Galois cover of $X$. Let $T$ be a smooth compact manifold on which
$\Gamma$ acts by diffeomorphisms. We consider
$\widehat{M}=\widetilde{X}\times T$, $ \pi=\text{projection onto
the second factor}$, $M=\widetilde{X}\times_{\Gamma}
T:=(\widetilde{X}\times T)/\Gamma$ where we let $\Gamma$ act on
$\widetilde{X}\times T$ diagonally. The leaves of the foliation
${\cal F}$ are the images of the manifolds $\widetilde{X}\times
\{\theta\}$, $\theta\in T$. The foliated manifold $(M,{\cal F})$
is usually refered to as a {\it foliated $T$-bundle.}

\s
\n
As a particular example of this construction consider a closed
smooth closed riemann surface $\Sigma$ of genus $g>1$
and let $\Gamma=\pi_1 (\Sigma)$, a discrete subgroup of $PSL(2,\RR)$.
Let $\{p_1,\dots,p_k\}$ points in $\Sigma$ and let $D_j\subset \Sigma$
be  a small open disc around $p_j$. Let $D=\cup_{j=1}^{k} D_j$.
Then we can
consider $X:=\Sigma\setminus D$, $\Gamma\rightarrow\widetilde{X}
\rightarrow X$ the Galois cover induced by the universal cover
$\HH^2\rightarrow \Sigma$, $T=S^1$, with $\Gamma$ acting
on $S^1$ by fractional linear transformations.

\subsection{The groupoid $T\rtimes \Gamma$}\label{groupoid}
$\;$

\medskip  We consider the groupoid  $G = T \crpr \Gamma$
with set of morphisms $G^{(1)}=T\times \Gamma$ and  base $G^{(0)}=
T$. The range and source maps are respectively given by: $$
\forall (\theta,g) \in T \times \Gamma,\;
r(\theta,g)=\theta,\;\;\; s(\theta,g)= \theta\cdot g\,. $$ The
composition is defined as follows: $$(\theta,g)\cdot
(\theta^\prime,g^\prime)=(\theta,g
g^\prime)\;\;\text{if}\;\;\theta^\prime=\theta g\,.$$ The inverse
of $(\theta,g)$ is  $(\theta g,g^{-1})$. For more about groupoids
we refer the reader to \cite{Co}.

The algebraic cross-product $C^\infty_c(T) \crpr \Gamma$ is, by
definition, the set of functions $\sum_{g \in \Gamma} t_g(\theta)
g$ such that only a finite number of the $ t_g(\cdot) \in
C^\infty_c(T)$ do not vanish identically. We shall identify any
function  $f$ having compact support $$ f: T \crpr \Gamma
\rightarrow \CC $$ $$ (\theta,g) \rightarrow f(\theta,g) $$ with $
\sum_{g \in \Gamma} f(\theta,g) g$. Then one has: $$
\sum_{g^\prime  \in \Gamma} f^\prime (\theta, g^\prime) g^\prime
\cdot \sum_{g  \in \Gamma} f (\theta, g) g\,=\, \sum_{h \in
\Gamma} \left( \sum_{g\in\Gamma} f^\prime
(\theta,g^\prime)\,f(\theta\cdot g^\prime,(g^\prime)^{-1} h)
\right) h $$ where where we recall that $$g^\prime \cdot
(f(\theta,g) g) = f(\theta\cdot g^\prime,g) g^\prime g\,.$$ The
algebraic cross-product $C^\infty_c(T) \crpr \Gamma$ will also be
denoted by $C^\infty_c(T \crpr \Gamma)$. One can introduce the
reduced $C^*$-algebra $C^*_r (T\rtimes \Gamma)$ associated to the
groupoid $T\rtimes \Gamma$ as a suitable completion of the
algebraic cross product $C^\infty_c (T) \crpr \Gamma$. See
\cite{Co}. It is well known, and easy to check, that  there is a
natural isomorphism between the reduced $C^\ast-$algebra
$C_r^\ast(T\rtimes \Gamma)$ of the groupoid $T\rtimes \Gamma$ and
the cross-product algebra $C^0(T)\rtimes_r \Gamma$ (see, for example,
\cite{MoS})
; we
shall henceforth identify these two
$C^*$-algebras.

What we have described  is an example of a proper cocompact
$G$-manifold $P$ with $G$ an \'etale groupoid, see \cite{Co} (page
137) for the definition. In our case
$$G=T\rtimes\Gamma\,,\quad\quad G^{(0)}=T\,,\quad\quad
(\alpha:P\to G^{(0)})\equiv (\pi:\wm \to T).$$ We shall deal with
the general case in Section \ref{genetale}.

\subsection{$\Gamma$-equivariant families of operators}\label{g-families}$\;$

\medskip

We consider a $\Gamma$-equivariant complex hermitian vector bundle
$\widehat{V} \rightarrow \widehat{M} $ endowed with a
$\Gamma-$invariant $b-$hermitian connection $\widehat{\nabla}$ satisfying
$\widehat{\nabla}_{x \partial_x}=0$ on the boundary $\partial \widehat{M}$.
We then set $\widehat{E}= S^Z \otimes \widehat{V}= \widehat{E}^+ \oplus
\widehat{E}^-$ which defines a smooth $\Gamma-$invariant family of  $\ZZ_2-$graded hermitian
Clifford modules on the fibers $\pi^{-1}(\theta), \, \theta \in T$.
We then get a smooth family of $\Gamma-$invariant  $\ZZ_2-$graded Dirac
type operators
$$  D(\theta)= \begin{pmatrix} 0 & {D^-(\theta)}
\cr {D^+(\theta)}  & 0 \cr
\end{pmatrix},\; \theta\in T
$$
 acting fiberwise on
$C^\infty_c( \widehat{M},\, \widehat{E})$.
 Moreover in a collar neighborhood ($ \sim [0,1] \times \partial \pi^{-1}(\theta) =
\{(x, y)\}$) of  $\partial \pi^{-1}(\theta) $ we may write:
$$
 {D^+(\theta)}= \sigma (x \partial_x + D_0(\theta)\,)
$$ where $ D_0(\theta)$ is the induced boundary Dirac type operator acting
on
$$
 C^\infty(\partial \pi^{-1}(\theta) ,\, \widehat{E}^+_{|_{\partial \pi^{-1}(\theta)}}\,).
$$ Observe that our family can also be thought as a longitudinal
operator on $(M,{\cal F})$ acting on the sections of
$E:=\widehat{E}/\Gamma$.

\subsection{$C^*_r(T\rtimes \Gamma)$-Hilbert modules}\label{h-modules}$\;$

\medskip
We shall now describe how the $\Gamma$-equivariant family
$(D(\theta))_{\theta\in T}$ defines a $C(T)\rtimes_r
\Gamma-$linear operator on  suitable $C(T)\rtimes_r
\Gamma-$Hilbert modules.

\n
Recall that  $\widehat{E}$ is a $\Gamma-$equivariant hermitian vector
bundle over $\widehat{M}$ so that for each $(g,p) \in \Gamma \times \widehat{M}$
there is a unitary linear map
$U_{g,p}:\, \widehat{E}_{p\cdot g} \rightarrow \widehat{E}_{p}$. Then
we endow $C^\infty_c(\widehat{M}, \widehat{E} )$ with the structure of
a left $C^\infty_c(T) \crpr \Gamma-$module by setting
for any $s \in C^\infty_c(\widehat{M}, \widehat{E} )$ and
$\sum_{g \in \Gamma} f(. ,g) g \in C^\infty_c(T) \crpr \Gamma$
\begin{equation} \label{U}
\forall p \in \widehat{M},\; \;
\left( \sum_{g \in \Gamma} f(. ,g) g \,\cdot s \right)\, (p):=
 \sum_{g \in \Gamma} f(\pi(p), g) (R^\ast_g  s\,)(p)
\end{equation}
 where  $(R^\ast_g s)(p) = U_{p,g} (s(p\cdot g))$; observe that
$R^\ast_g \circ R^\ast_{g^\prime} = R^\ast_{gg^\prime}$.
\begin{lemma}  \label{D} The family of Dirac operators $(D(\theta))_{\theta \in T}$ acting fiberwise
on $C^\infty_c(\widehat{M}, \widehat{E} )$
 defines a left $C^\infty_c(T \crpr \Gamma)-$linear endomorphism $\D$ of  $C^\infty_c(\widehat{M}, \widehat{E} )$.
\end{lemma}

\begin{proof}
 Using the above notations we have:
$$ \D \bigl( \, \sum_{g \in \Gamma} f(. ,g) g \,\cdot s \,
\bigr)(p) =
 \sum_{g \in \Gamma} f(\pi(p), g) \, D(\pi(p)) (R^\ast_g  s\,) \, (p)
$$ where we have used the fact that $(D(\theta))_{\theta \in T}$
is a
 family of operators, i.e.
commutes with the natural action of $C^\infty(T)$. Since  the
family $(D(\theta))_{\theta \in T}$ is $\Gamma-$equivariant,
 the right hand  is by definition equal to
$$ \left( \sum_{g \in \Gamma} f(. ,g) g \,\cdot \D( s)\right) \,
(p) $$ which proves the lemma.
\end{proof}

Let $\dot{C}^\infty_c (\wm,\we)$ the space of sections of compact
support vanishing of infinite order at $\pa \wm$. We  define the
$C^0(T) \rtimes_r \Gamma-$hermitian product of two sections $s$
and $s^\prime$ of $\dot{C}^\infty_c(\widehat{M}, \widehat{E} )$ by
setting:
\begin{equation} \label{hermitianproduct}
\langle s;s^\prime\rangle \,=\,
\sum_{g\in \Gamma} \langle s;s^\prime\rangle (\theta,g) g
\end{equation}
where $\forall (\theta,g) \in T \times \Gamma $:
$$
 \langle s;s^\prime\rangle(\theta,g)\,=\,
\int_{\pi^{-1}(\theta.g)}\, <R^\ast_{g^{-1}} (s)(y) ; s^\prime(y) >_{\widehat{E}}
d {\rm Vol}^b_{\pi^{-1}(\theta\cdot g)}(y)
$$ where $ d{\rm Vol}^b_{\pi^{-1}(\theta\cdot g)}(y)$ denotes the ${\rm b}-$riemannian
density in the fiber.

We shall sometimes write, more shortly, $d {\rm Vol}^b (y)$ (the
domain of integration appears already in the integral). With the
definition (\ref{hermitianproduct}) one has:
\begin{lemma}
$$
\langle t \cdot s;s^\prime\rangle\,=\, t\cdot \langle s;s^\prime\rangle,
\quad\quad \forall t \in C^\infty_c(T \crpr \Gamma)\,.
$$
\end{lemma}
\begin{proof} We may assume that $t= t_\gamma(.) \gamma$
for $\gamma \in \Gamma$.
Then for any $(\theta, g) \in T \times \Gamma$ we have:
$$
\langle t_\gamma (.) \gamma \cdot s;s^\prime\rangle(\theta,g)=
\int_{\pi^{-1}(\theta\cdot g)}\, <R^\ast_{g^{-1}} (t_\gamma
(\pi(.) ) R^\ast_\gamma s)(y) ; s^\prime(y) >_{\widehat{E}}
d{\rm Vol}^b_{\pi^{-1}(\theta\cdot g)}(y)
$$
$$=t_\gamma (\pi(y\cdot g^{-1})) \;\int_{\pi^{-1}(\theta\cdot g)}\,
 < (R^\ast_{g^{-1}\gamma} s)(y) ; s^\prime (y) >_{\widehat{E}}
d{\rm Vol}^b_{\pi^{-1}(\theta g)}(y)
$$
but this last term is equal to
$
\langle s;s^\prime\rangle(\theta\cdot \gamma, \gamma^{-1}g)
t_\gamma (\theta )
$ so that
$$
\sum_{g \in \Gamma}
\langle t_\gamma (.) \gamma \cdot s;s^\prime\rangle(\theta,g) \,g =
\sum_{g \in \Gamma} \langle s;s^\prime\rangle(\theta. \gamma, \gamma^{-1}g)
t_\gamma (\theta ) \gamma \, (\gamma^{-1}g ).
$$
But we can write the right hand side  as
$$\sum_{h\in \Gamma} <s;s^\prime>(\theta\gamma,h)t_\gamma(\pi(\theta))
\gamma h$$
which is equal to
$$t_\gamma(\pi(\theta))\gamma
\sum_{h\in \Gamma} <s;s^\prime>(\theta,h)h$$
which is exactly
$t\cdot \langle s;s^\prime\rangle$ as required. The lemma is proved.
\end{proof}

\begin{definition}
The completion of $\dot{C}^\infty_c(\widehat{M}, \widehat{E} )$
 with
respect to the hermitian scalar product (\ref{hermitianproduct})
defines a $C^0(T) \rtimes_r \Gamma-$Hilbert module denoted by
$L^2_{b, C^0 (T) \rtimes_r \Gamma}(\wm,\widehat{E})$.
\end{definition}

Similarly one defines the b-Sobolev spaces $H^m_{b, C^0(T)
\rtimes_r \Gamma}(\wm,\widehat{E})$ for $m \in \ZZ$  where the
$L^2$-condition is defined with respect to the space ${\rm
Diff}^m_{b,\rtimes}(\wm,\we)$
  of
$\Gamma-$invariant fiberwise $b-$differential operators of order
$m$. If $m\geq 0$ then $$ H^m_{b, C^0(T) \rtimes_r
\Gamma}(\wm,\widehat{E})=\{s\in L^2_{b, C^0 (T) \rtimes_r
\Gamma}(\wm,\widehat{E}); Pu\in L^2_{b, C^0 (T) \rtimes_r
\Gamma}(\wm,\widehat{E}) \,\forall P\in {\rm
Diff}^m_{b,\rtimes}(\wm,\we)\}$$ If $m<0$ then $$ H^m_{b, C^0(T)
\rtimes_r \Gamma}(\wm,\widehat{E})= {\rm
Diff}^{-m}_{b,\rtimes}(\wm,\we)\,L^2_{b, C^0 (T) \rtimes_r
\Gamma}(\wm,\widehat{E})$$ Our $\Gamma$-equivariant family
$(D(\theta))_{\theta\in T}$ defines a morphism $$\D:
\dot{C}^\infty_c (\wm,\we)\longrightarrow \dot{C}^\infty_c
(\wm,\we)$$ of pre-hilbert
 $(C^0(T) \rtimes_r\Gamma)$-modules which extends for each $m\in \ZZ$ to a bounded
morphism $$ H^m_{b, C^0(T) \rtimes_r \Gamma}(\wm,\widehat{E})
\longrightarrow  H^{m-1}_{b,
C^0(T)\rtimes_r\Gamma}(\wm,\widehat{E})\,.$$

 We end this section by a useful

\medskip
\begin{lemma}\label{standard} For each $ m\in\ZZ$,
the Hilbert-modules $H^m_{b, C^0(T) \rtimes_r
\Gamma}(\wm,\widehat{E})$ are standard  i.e. isomorphic to
$l^2(C^0(T)\rtimes_r \Gamma) $.
\end{lemma}

\begin{proof} We shall only consider the case $m=0$. There exists
a positive integer $N$ and open connected subsets
$U_i \subset U_i^\prime$ of $T$ ($1 \leq i \leq N$) with the following
properties. Each $U_i$ is relatively compact in $U_i$,
$\cup_{1\leq i \leq N} U_i = T $, and for each $i\in \{1,\ldots,N\}$
the restriction of the fibration $\pi$ to
$\pi^{-1}(U_i^\prime)$ is trivial: $\pi^{-1}(U_i^\prime) \simeq
U_i^\prime \times Z$. Denote by $\mu$ the given $\Gamma-$invariant
 riemannian measure and $v_f$ the volume of a fundamental domain
 for the action of $\Gamma$ on $\wm$. Then using an induction
 argument on $i\in \{1, \ldots,N\}$, we may find an open connected
 subset $W \subset Z$ and  open
 subsets $V_i^\prime \subset \pi^{-1}(U_i^\prime)$  ($1 \leq i \leq N$)
 with the following three properties:

 \n 1)  $ V_i^\prime \simeq U_i^\prime \times W$ for each $i\in \{1, \ldots,N\}$
  and $ \mu( \cup_{1 \leq i \leq N} V_i^\prime ) \leq
 {1 \over 2} v_f$

 \n 2) $ \forall \gamma \in \Gamma \setminus \{e\}, \,
 \forall i \in \{1,\ldots , N \},\; V_i^\prime\cdot \gamma \cap V_i^\prime =
 \emptyset$.

\n 3) $ \forall \gamma \in \Gamma, \,
 \forall i,j \in \{1,\ldots , N \},\,
 {\rm { with}}\, i \not=j,\; V_i^\prime\cdot \gamma \cap V_j^\prime =
 \emptyset$.

 Then, for each $i\in \{1,\ldots , N \}$, we may find a sequence
 $(e^i_l)_{l\in \NN}$ of elements of
 $C^0(\overline{V_i^\prime},\widehat{E})$, vanishing in a neighborhood
 of the boundary,
  such that
 for any $\theta \in U_i^\prime$ and any $l, l^\prime \in \NN$:

$$\int_{\pi^{-1}(\theta)}
 < e^i_l (y) ; e^i_{l^\prime} (y) >_{\widehat{E}}
d{\rm Vol}^b_{\pi^{-1}(\theta )}(y) = \delta_{l,l^\prime}.
$$
 Next we consider for each $i\in \{1,\ldots , N \}$ a smooth
function $\phi_i \in C^\infty_c(U_i^\prime)$ such that $\phi_i =
1$ on $U_i$ and we set for each $l \in \NN$: $$ \xi_l \,=\, {
1\over \sqrt{\sum_{i=1}^N \phi_i^2}} \sum_{i=1}^{n} \phi_i e^i_l
\in L^2_{b, C^0(T) \rtimes_r \Gamma}(V_i^\prime,\widehat{E}) $$
Then, using the above properties, one checks that for any $l,
l^\prime \in \NN$ one has $$<\xi_l ; \xi_{l^\prime} > = \delta_{l,
l^\prime} \; {\rm in }\; C^0(T) \rtimes_r \Gamma. $$ Then the
$\xi_l,\, l\in \NN$ generate a closed Hilbert submodule $\cal{H}$
of $L^2_{b, C^0(T) \rtimes_r \Gamma}(\wm,\widehat{E})$ which is
isomorphic to $l^2(C^0(T) \rtimes_r \Gamma ) $ and such that $$
L^2_{b, C^0(T) \rtimes_r \Gamma}(\wm,\widehat{E}) = \cal{H} \oplus
\cal{H}^\perp. $$ The Lemma therefore follows from Kasparov's
stabilization theorem, see \cite{Bla}.
\end{proof}

\subsection{Foliations}$\;$\label{foliations}

\def \crpr {>\!\!\!\triangleleft\,}
Now, for the convenience of the reader, we connect the
noncommutative picture we have described above with the more
familiar noncommutative picture arising from the foliated manifold
$(M,\mathcal{F})$, as in the work of A. Connes.
We follow closely Morioshi-Natsume
\cite{Moriyoshi-Natsume (1996)}. In this subsection only  we
assume, for simplicity, that $\widehat{M}= \widetilde{M} \times
T$, the fibration $\pi$ is given by the projection
 $\pi: \widetilde{M} \times T \rightarrow T$ and the group
 $\Gamma$ acts on the  manifold with boundary $\widetilde{M}$ (and on $T$)
 in such a way that the quotient $\widetilde{M}/\Gamma$ is a smooth compact
 manifold with boundary. So the product foliation
 $\widetilde{M} \times \{\theta\}$, $\theta \in T$ descends
 to a foliation $\cal{F}$ on $M$. We shall assume that the $\Gamma-$action
 on $T$ satisfies the following:

\s
 \noindent {\bf Condition.} For $\gamma \in \Gamma$, if there exists an
 open subset $U$ of $T$ such that $\gamma(\theta) = \theta$ for
 any $\theta \in U$, then $\gamma$ is the identity element of $\Gamma$.

\s
\n
 This condition guarantees that the holonomy groupoid $G$
 of $(M, \cal{F})$ is Hausdorff and given by:
 $$
 G \simeq ( \widetilde{M} \times \widetilde{M} \times T)/\Gamma
 $$ where $\gamma \in \Gamma$ acts by $(y,z, \theta).\gamma =
 (y.\gamma, z.\gamma,  \theta.\gamma)$. The source map $s$ and the
 range map $r$ are given by:
 $$
 r([y,z, \theta])=[y, \theta],\quad s([y,z, \theta])=[z, \theta].
 $$

 Fix a $b-$metric on $\widetilde{M}/\Gamma$, the lifting to
 $\widetilde{M}$ of this $b-$metric induces a left Haar system
 $\{\nu^y\}$ on the groupoid $G$, see \cite{Ren}. We now recall the definition
 of the $C^\ast-$algebra of the foliation $\cal{F}$ with coefficient
 in the hermitian vector bundle $E\rightarrow M$ whose lift
 to $\widehat{M}$ is $\widehat{E}$. Denote by
 $\dot{C}^\infty_c(G,E)$ the space of all compactly supported smooth sections of the bundle
 $(s^\ast(E))^\ast \otimes r^\ast(E)$ vanishing of infinite order at the boundary.
 The space $\dot{C}^\infty_c(G,E)$ has
 a $\ast-$algebra structure:
 $$
 (f_1 \ast f_2)(\gamma) = \int_{G^{r(\gamma)}}f_1(\gamma^\prime)
 f_2({\gamma^{\prime}}^{-1}\gamma) \, d\nu^{r(\gamma)}(\gamma^\prime),
 $$
 $$
 f^\ast(\gamma)=(f(\gamma^{-1}))^\ast.
 $$
where we recall that $G^x:=\{\gamma\in G; r(\gamma)=x\}$.
Now let $\lambda_\theta, \theta \in T$ be the strictly positive lifted
 $b-$density on $\widetilde{M}\times \{\theta\}$ and set
 $$
 H_\theta=L^2_b(\widetilde{M}\times \{\theta\};
  \widehat{E}_{|\widetilde{M}\times \{\theta\}}, \,\lambda_\theta).
  $$ Then the collection $\mathcal{H}=(H_\theta)_{\theta \in T}$
  together with the space of compactly supported vanishing near the boundary
  continuous sections of the bundle
   $\widehat{E}$ over $\widetilde{M} \times T$, defines a continuous
   field of Hilbert spaces over $T$.  The $\Gamma-$action on
   $\widetilde{M} \times T $ and $\widehat{E} $ gives rise to an action
   on $\mathcal{H}$ denoted by $\xi \rightarrow \gamma \xi$
   for $\gamma \in \Gamma$ and $\xi $ a section of $\mathcal{H}$. The space
   ${\rm End}_\Gamma(H)$ of $\Gamma-$equivariant bounded measurable fields
   of operators $K=(K_\theta), K_\theta \in B(H_\theta),$ is
   a $C^\ast-$algebra, where the norm is given by
   $$
   ||K||= \sup\{||K_\theta||;\, \theta \in T\}.
   $$ There is a faithful representation $\rho:\dot{C}^\infty_c(G,E) \rightarrow
   {\rm End}_\Gamma(\mathcal{H})$ defined by
   $$
   \forall f \in \dot{C}^\infty_c(G,E),\quad
   [\rho(f)_\theta \xi](y)= \int_{} f(y,z, \theta) \xi(z) d\lambda_\theta(z),
   $$ where $\xi \in H_\theta$. The norm closure of
   $ \dot{C}^\infty_c(G,E)$ with respect to the norm
   $$
   ||f||=||\rho(f)||= \sup\{\,||\rho(f)_\theta||;\; \theta \in T \}
   $$ is by definition the $C^\ast-$algebra $C^\ast(M, \cal{F}, E)$
   of the foliation $(M, \cal{F})$ with coefficient  $E$.

   Now we define a right action of $C^\infty_c(T\rtimes \Gamma)$
   on the space $\dot{C}_c^\infty(\widetilde{M}\times T; \widehat{E})$ of
   sections vanishing of infinite order at the boundary as follows.
   $$
   (\xi f)(z,\theta) =\sum_{\gamma \in \Gamma}
   f(\theta \gamma^{-1}, \gamma) \xi(z\gamma^{-1},\theta \gamma^{-1} ), \quad
   \xi \in \dot{C}_c^\infty(\widetilde{M}\times T; \widehat{E}), \quad
   f
   \in C^\infty_c(T\rtimes \Gamma).
   $$ A $C^\infty_c(T\rtimes \Gamma)-$valued inner product $<. ; .>$
   on $ \dot{C}_c^\infty(\widetilde{M}\times T; \widehat{E})$ is defined by
   $$
   <\xi_1 ; \xi_2>(\theta, \gamma) =
   \int_{\widetilde{M} \times \{\theta\}}
   <\xi_1(z,\theta) ; \xi_2(z\gamma,\theta \gamma)>_{\widehat{E}}
   d\lambda_\theta(z),
   $$ where $< . ; .>_{\widehat{E}}$ is the hermitian scalar product
   of the vector bundle $\widehat{E}$.
   Recall now that $C^\infty_c(T\rtimes \Gamma)$ is a dense sub-algebra
   of the reduced $C^\ast-$algebra $C(T)\rtimes_r \Gamma$.
   Denote as before
by $L^2_{b,C^0 (T)\rtimes_r \Gamma}(\widetilde{M}\times T,\widehat{E})$ the
   $C(T)\rtimes_r \Gamma-$Hilbert module obtained by taking the
   $C(T)\rtimes_r \Gamma-$completion of $ \dot{C}_c^\infty(\widetilde{M}\times T; \widehat{E})$
   with respect to the norm
   $||\xi||= \sqrt{||<\xi ; \xi >||_{C(T)\rtimes_r \Gamma}}$.

   The fundamental point here is then the following:
    the  left action of
   $ \dot{C}_c^\infty(G; \widehat{E})$ on $L^2_{b,C^0 (T)\rtimes_r \Gamma}(\widetilde{M}\times T,\widehat{E})$ given by
    $$
   (f\ast\xi)(y,\theta)=\int_{\widetilde{M}\times \theta}
   f(y,z,\theta) \xi(z,\theta) d\lambda_\theta(z)
   $$ (where we have used the identification of $G$ with
   $(\widetilde{M}\times \widetilde{M}\times T)/\Gamma$), extends to a
faithful
   representation of $C^* (M, \cal{F}, E)$ and  the image of
   this representation is precisely the space
   $\KK_{C(T)\rtimes _r \Gamma}(L^2_{b,C^0 (T)\rtimes_r \Gamma}(\widetilde{M}\times T,\widehat{E}) )$ of
   $(C(T)\rtimes_r \Gamma)-$compact operators on $L^2_{b,C^0 (T)\rtimes_r \Gamma}(\widetilde{M}\times T,\widehat{E})$.
The proof of this fact is exactly as in \cite{Moriyoshi-Natsume
(1996)}, Proposition 2.4.

\subsection{Invertible families}\label{invertible}$\;$

\medskip
We consider a closed $\Gamma$-manifold $\widehat{N}$ fibering over
$T$. The fibration $\pi:\widehat{N}\longrightarrow T$ is assumed
to be $\Gamma$-equivariant and the fibers to be of dimension
$2k-1$. One example of this situation is given by the boundary
$\pa\widehat{M}\to T$ of a fibration with boundary as in the
previous subsections, but we want to proceed in full generality
now. We endow $\widehat{N}$ with a $\Gamma$-invariant riemannian
metric and we denote by $d(\,,\,)$ the induced distance. The group
$\Gamma$ is assumed to act freely properly and cocompactly on
$\wn$ so that the quotient  $N:=\widehat{N}/\Gamma$ is a compact
smooth manifold. Moreover  we assume that the fibers of $\pi$
carry a $\Gamma-$equivariant spin structure and denote by
$S\rightarrow \wn$ the associated spinor bundle. We also consider
a $\Gamma$-equivariant hermitian complex vector bundle
$\widehat{V} \rightarrow \wn$ endowed with a $\Gamma-$invariant
connection. We then set
 $\widehat{F}=S \otimes \widehat{V}$ on $\widehat{N}$.
We get in this way a $\Gamma$-equivariant family $(D_0
(\theta))_{\theta\in T}$ acting along the fibers of
$\pi:\widehat{N}\longrightarrow T$.

Some of the results in this paper are proved under the following
assumption:

\medskip
\n
{\bf Hypothesis A.} {\it There exists a real $\epsilon>0$ such that
for any $\theta\in T$, the $L^2$-spectrum
of $D_0 (\theta)$ acting on $L^2(\pi^{-1}(\theta),\wf_{|_{\pi^{-1}(\theta)}})$
does not meet $]-\epsilon,\epsilon[$.}

\smallskip
\n
{\bf Example.}
\\ Let $Y$ be a {\it spin} compact manifold without boundary
and let $\Gamma\rightarrow \widetilde{Y}\rightarrow Y$
be a Galois cover of $Y$. Let $T$ be a smooth compact manifold
on which $\Gamma$ acts by diffeomorphisms. We consider
$\widehat{N}=\widetilde{Y}\times T$, $ \pi=\text{projection onto
the second factor}$, $N=\widetilde{Y}\times_{\Gamma} T
:=(\widetilde{Y}\times T)/\Gamma$
where we let $\Gamma$ act on $\widetilde{Y}\times T$ diagonally.
We consider a hermitian vector bundle $V$ on $N$ with hermitian
connection $\nabla^V$ and its lift
$\widehat{V}$ on $\widehat{N}$.  We denote by
$\widehat{V}_\theta$ the restriction of $ \widehat{V}$
to $\widetilde{Y}\times \{\theta\}$.
We endow $\widehat{V}$ with
the pulled-back metric and with the pulled-back  connection.
We now assume that $Y$ admits a metric of
{\it positive scalar curvature} and we fix such a metric $g_Y$;
we let $g_{\widetilde{Y}}$ be the pulled-back metric.
Let $\widetilde{D}$ be the associated Dirac operator on $\widetilde{Y}$.
We consider the family $$D(\theta)=\widetilde{D}_{\widehat{V}_\theta}$$
where on the right hand side we are considering the twisted
Dirac operator.
This is a $\Gamma$-equivariant family; if the  scalar
curvature of $g_Y$ is big enough, then, by Lichnerowicz
formula, we easily obtain that $(D(\theta))_{\theta\in T}$
satisfies Hypothesis A.\\We now give a particular example
showing the existence of a {\it boundary} family satisfying
Hypothesis A.

\n
Let $\Gamma$ be a finitely presented group. There exists
 a $(2k+1)$-dimensional closed spin manifold $N$
such that $2k+1\geq 7$ and  $\Gamma=\pi_1(N)$.
We consider $X_1=B^9\times N$, with $B^9$ the
unit ball in the standard euclidean space $\RR^9$.
We perform surgery along a submanifold
of $X_1\setminus \pa X_1$ of codimension at least 3
and we obtain a manifold with boundary equal to $Y=S^8\times N$
and fundamental group $\Gamma$. Let $\widetilde{X}$ be the
universal cover of $X$, let  $\widetilde{Y}=\pa\widetilde{X}$
and let $\wm=\widetilde{X}\times T$,
$M=(\widetilde{X}\times T)/\Gamma$. Fix a
hermitian bundle $V$ on $M$ and consider its
pull-back $\widehat{V}$ to $\wm$, as in the example given above.
Consider now the family $$D(\theta)=\widetilde{D}_{\widehat{V}_\theta}$$
with $\widetilde{D}$ denoting now the Dirac operator on $\widetilde{X}$.
This is a family acting on the section of the bundle $\we=S\otimes \widehat{V}$
with $S$ equal to the spinor bundle on $\widetilde{X}$.
Let $\widetilde{D}_0$ be the Dirac operator on
$\widetilde{Y}=\pa\widetilde{X}$.
The boundary family associated to  $(D^+ (\theta))_{\theta\in T}$ is given
by the operator $\widetilde{D}_0$ twisted by
$(\widehat{V}_\theta)|_{\pa \widetilde{X}\times\{\theta\}}$.
Let us fix a metric of positive scalar curvature on $Y$
and let us extend this metric to $X$; if the scalar curvature
of $Y$ is big enough then the boundary family $(D_0 (\theta))_{\theta\in T}$
will
satisfy Hypothesis A.
A particular
 example of this situation is $\Gamma=\ZZ^k$
acting on $T=(S^1)^k$, the $k$-dimensional torus,
by
$$(n_1,\dots,n_k)\cdot (e^{i\theta_1},\dots,e^{i\theta_k})=
(e^{i\theta_1+2i\pi n_1\alpha_1},\dots,e^{i\theta_k+2i\pi n_k\alpha_k})$$
where $\alpha_1,\dots,\alpha_k$ are $k$  fixed irrational numbers.
\\
Another particular example is given by taking
$\Gamma=\pi_1(\Sigma_g)$, where
$\Sigma_g$ equal to a closed riemann surface of genus $g\geq 2$.
Then one can take $T=S^1$ with $\Gamma$ acting by fractional
linear tranformations.

\medskip
We introduce the space $\Psi^*_\pi(\widehat{N},\widehat{F})$
of smoothly varying families of pseudodifferential operators
along the fibers.
The Schwartz kernel of an element $P\in \Psi^*_\pi(\widehat{N},\widehat{F})$
is a distribution on the fibred product $\widehat{N}\times_\pi \widehat{N}$
with the usual conormal singularities on the fibre diagonal.
The space of $\Gamma$-equivariant elements in
$\Psi^*_\pi(\widehat{N},\widehat{F})$ will be denoted by
$\Psi^*_{\rtimes} (\widehat{N},\widehat{F})$.
An element $P\in \Psi^*_\pi(\widehat{N},\widehat{F})$ has compact
$\Gamma$-support if its Schwartz kernel has compact support
in $(\widehat{N}\times_\pi \widehat{N})/\Gamma$: we denote {\it the algebra}
of $\Gamma$-equivariant pseudodifferential operators of compact
$\Gamma$-support
by $\Psi^*_{\rtimes,c}(\widehat{N},\widehat{F})$.
For example
$$(D_0 (\theta))_{\theta\in T}\in
{\rm Diff}^1_{\rtimes}(\widehat{N},\widehat{F})\subset
\Psi^1_{\rtimes,c}(\widehat{N},\widehat{F}).$$

It is important to understand what Hypothesis A implies for the
$C^0(T)\rtimes_r \Gamma$-linear operator $\D_0$ induced by the
family $(D_0 (\theta))$, acting on $L^2_{C(T)\rtimes_r \Gamma}.$

\begin{proposition}\label{invertibility}
\item{1]} The operator $\D_0$ defines a regular unbounded operator
on the $C^0 (T)\rtimes_r \Gamma$-Hilbert module
$L^2_{C^0(T)\rtimes_r \Gamma}$.
\item{2]} If Hypothesis A holds then
this operator is $L^2_{C^0(T)\rtimes_r \Gamma}$-invertible, with
inverse induced by the $\Gamma$-equivariant family of operators
$\{D_0 (\theta))^{-1}\}_{\theta\in T}$.
\item{3]} One can write
$D_0 (\theta)^{-1}=A(\theta)+R(\theta)$ with $\{A
(\theta)\}_{\theta\in T}\in
\Psi^1_{\rtimes,c}(\widehat{N},\widehat{F})$ and $\{R
(\theta)\}_{\theta\in T}\in
\Psi^{-\infty}_{\rtimes}(\widehat{N},\widehat{F})$. Moreover $\{R
(\theta)\}_{\theta\in T}$ extends to an operator $$\mathcal{R}\in
\BB(H^m_{C^0 (T)\rtimes_r \Gamma} (\widehat{N},\widehat{F} ),
H^k_{C^0 (T)\rtimes_r \Gamma} (\widehat{N},\widehat{F} ))$$
$\forall k,m\in\ZZ$.

\end{proposition}
\begin{proof}

\n
1] The fact that $\D_0$ defines a regular unbounded operator is
standard.

\n
2] Let us prove that $\D_0$ is $L^2_{C(T)\rtimes_r
\Gamma}$-invertible. Since $\D_0$ is regular, we may consider
$$\frac {\D_0} { \sqrt{\Id + \D_0^2}} \in \BB( L^2_{C^0 (T)\rtimes_r
\Gamma}(\widehat{N},\widehat{F})) $$ instead of $\D_0$.
 For each $\theta \in T$ let $H_\theta$ denote the
Hilbert space $H_\theta=
L^2(\pi^{-1}\theta,\widehat{F}|_{\pi^{-1}(\theta)})$. The
collection $\mathcal{H}=(H_\theta )_{\theta \in T} $ together with
$C^0_c(\widehat{N} ; \widehat{F})$ defines a continuous field of
Hilbert spaces over $T$. The space ${\rm
{End}}_\Gamma(\mathcal{H})$ of $\Gamma-$equivariant bounded
measurable fields of operators $A=(A_\theta )_{\theta \in T}$
($A_\theta\in H_\theta$) is a $C^\star-$algebra where the norm is
given by: $$ || A|| = \sup \{ ||A_\theta||, \; \theta \in T\}. $$
Moreover, one has a natural injective morphism of
$C^\star-$algebras: $$ J: \BB (L^2_{C(T)\rtimes_r \Gamma}(
\widehat{N},\widehat{F}   )) \rightarrow {\rm {End}}_\Gamma(H). $$
So for any $\tilde{A} \in \BB( L^2_{C(T)\rtimes_r
\Gamma}(\widehat{N},\widehat{F})  )$ the spectrum of $\tilde{A}$
coincides with the one of $J(\tilde{A} )$. Since Hypothesis A
means that $J(\frac {\D_0} { \sqrt{\Id + \D_0^2}})$ is invertible
we get immediately that $\D_0$ is $L^2_{C(T)\rtimes_r
\Gamma}$-invertible.

\n
3] The proof is standard and left to the reader.

\end{proof}

\section{{\bf \'Etale groupoids and the
$b$-calculus}}\label{bcalc}

\subsection{The fibre-$b$-stretched product.}\label{stretched}$\;$

\medskip
Recall that we are given a $\Gamma$-equivariant fibration
$Z\rightarrow \wm \rightarrow T$, with $Z$ a fixed manifold with
boundary: we follow Section \ref{fibration} for the notation
adopted here. The boundary of the model fiber $Z$ will decompose
as a disjoint union of connected components: $$\partial
Z=\sqcup_{\alpha\in A} W^\alpha$$ Thus for each $\theta\in T$ the
fiber $\pi^{-1}(\theta):=\wm_\theta$ has a boundary $$\partial
\wm_\theta=\sqcup_{\alpha\in A} W^\alpha_\theta$$ with
$W^\alpha_\theta$  diffeomorphic to $W^\alpha$. We consider
$$\widehat{W}^\alpha := \cup_{\theta\in T} (W^\alpha_\theta)\,.$$
This is a connected component of $\partial \wm$ (recall that $T$
is assumed to be connected).


We shall now define the space that will carry the Schwartz kernels
of the operators we are interested in.
Let
$$\wm\times_\pi \wm=\{(p,p^\prime)\in \wm\times \wm\,|\, \pi(p)=\pi(p^\prime)\}$$
be the fiber product of $\wm$ with itself. This is a fibration over $T$, with
fiber over $\theta$ equal to $\wm_\theta \times \wm_\theta$.
Consider $$B=
\cup_\theta (\sqcup_{\alpha\in A} W^\alpha_\theta \times
W^\alpha_\theta)$$

\begin{definition}\label{stretched}
The  $b$-stretched product of the $\Gamma$-equivariant
fibration $\wm\to T$ is, by definition, the blow-up of
$\wm\times_\pi \wm$ along $B$:
\begin{equation}\label{stretched-f}
[\wm\times_\pi \wm;B]
\end{equation}
\end{definition}


 Notice that the  $b$-stretched
product of the $\Gamma$-equivariant fibration $\wm\to T$ is
exactly the fibre-$b$-stretched product considered in \cite{MP I},
but with a non-compact fibre.
On the   space $[\wm\times_\pi \wm;B]$ there is a well
defined $\Gamma$-action obtained by lifting  the diagonal action of
$\Gamma$ on the fibration $\wm\times_\pi \wm$.

\medskip
\n
{\it For more on the $b$-stretched
product and on what follows the reader is invited to consult the basic
reference \cite{Melrose}.}

\subsection{The $b$-calculus
$\Psi^{*,\delta}_{b,\rtimes}(\widehat{M},\widehat{E})$}\label{thebcalculus}
$\;$

\medskip
The small fiber-$b$-calculus $\Psi^*_{b,\pi}(\wm,\we)$ is defined,
exactly as in \cite{MP I}, in terms of the  fiber
$b$-stretched product $[\wm\times_\pi \wm;B]$ and the lifted fiber
diagonal $\Delta_{b,\pi}$.
We define the small {\it
$\Gamma$-invariant fibre-$b$-calculus} as
\begin{equation}\label{calculus1}
\Psi^*_{b,\rtimes}(\wm,\we):=\{P\in \Psi^*_{b,\pi}(\wm,\we)\, | \,
R_g^* \circ P=P\circ R^*_g  \,\;\forall \;g\in \Gamma \}
\end{equation}
 The Schwartz kernel of an operator in
$ \Psi^*_{b,\rtimes}(\wm,\we)$ is  a $\Gamma$-invariant distribution
on $[\wm\times_\pi \wm,B]$.  The subspace of operators
with  Schwartz kernel
 {\it compactly supported} in $[\wm\times_\pi \wm,B]/\Gamma$
will be denoted by $\Psi^*_{b,\rtimes,c}(\wm,\we)$;
as in the closed case
 we shall say that elements in the latter space have compact $\Gamma$-support.
The subspace of $b$-differential operators is denoted
$\Diff^*_{b,\rtimes}(\wm,\we)\,.$
We have already remarked that the family of Dirac operators $(D(\theta))_{\theta\in T}$
introduced in Section \ref{fibration}, defines an element in
$\Diff^1_{b,\rtimes}(\wm,\we)\subset \Psi^1_{b,\rtimes, c}(\wm,\we)$.

\medskip
Let us fix a $\Gamma$-invariant trivialization $\nu$
of the positive normal bundle $N_+\partial\wm$
to the boundary of $\pa\wm$. Let $x\in C^\infty(\wm)$ be
a $\Gamma$-invariant boundary defining function such that
$dx(\widehat{\nu})=1$.
Let $$\overline{N_+\partial\wm}\simeq \pa\wm\times [-1,1]$$
denote its fibre compactification.
The notion of indicial operator and  indicial family
are as in \cite{Melrose}. The indicial operator of
$P\in \Psi^{m}_{b,\rtimes}(\wm,\we)$, denoted $I(P)$,
is the $\RR^+$-invariant element in
$\Psi^{m}_{b,\rtimes}(\overline{N_+\partial\wm},\we_{|_{\pa\wm}})$
obtained by restricting the kernel of $P$ to the front face.
The indicial operator defines a short exact sequence
\begin{equation}\label{indicial}
0\longrightarrow \rho_{{\rm bf}}
\Psi^{m}_{b,\rtimes}(\wm,\we)
\longrightarrow
\Psi^{m}_{b,\rtimes}(\wm,\we)
\smash{\mathop{\longrightarrow}\limits^{I(\;)}}
\Psi^{m}_{b,\rtimes,I}(\overline{N_+\partial\wm},\we_{|_{\pa\wm}})
\longrightarrow 0
\end{equation}
with the last space denoting the $\RR^+$-invariant elements in $
\Psi^{m}_{b,\rtimes}(\overline{N_+\partial\wm},\we_{|_{\pa\wm}})$.
The indicial family of $P\in \Psi^{m}_{b,\rtimes}(\wm,\we)$ is
obtained by fiber-Mellin transform from $I(P)$, once the
trivialization $\nu$ has been fixed. It defines an  entire  map
with values $\Gamma$-invariant pseudodifferential operators along
the fibres of $\pa\wm$: $$\CC\ni z\rightarrow I_{\nu } (P,z)\in
\Psi^{m}_{\rtimes}(\pa\wm,\we_{|_{\pa\wm}})$$ where
$\Psi^{m}_{\rtimes}(\pa\wm,\we_{|_{\pa\wm}}) $
denotes the set of $\Gamma-$equivariant families of pseudo-differential
operators of order $m$ acting in the fibers of $\pi: \partial \wm \rightarrow T$.

\medskip
 Proceeding  as in \cite{MP I} (Appendix) and
\cite{LPMEMOIRS} (Section 7) we can also define the appropriate
full calculi with bounds
$$\Psi^{*,\delta}_{b,\pi}(\wm,\we):=
\Psi^{*}_{b,\pi}(\wm,\we)+\Psi^{-\infty,\delta}_{b,\pi}(\wm,\we)+
\Psi^{-\infty,\delta}_{\pi}(\wm,\we)
$$ for $\delta>0$. Here
$\Psi^{-\infty,\delta}_{b,\pi}(\wm,\we)$ is defined as in
\cite{LPMEMOIRS}, right after formula (7.3), but with
$[\wm\times_\pi \wm,B]$ replacing $\widetilde{M}^2_b$; similarly
the residual space $\Psi^{-\infty,\delta}_{\pi}(\wm,\we)$ is
defined as
\begin{equation}\label{rem}
\rho_{lb}^\delta \rho_{rb}^\delta H^\infty_{b,{\rm
loc}}(\wm\times_\pi \wm,\we\boxtimes \we^*).
\end{equation}

\smallskip
\noindent
{\bf Remark.}
We remark as in \cite{MP I} (Appendix) that $\Psi_{b,
\pi}^{-\infty, \delta}(\wm , \we)$ can also be obtained by
doubling $[\wm \times_{\pi} \wm, B]$ across the front face $bf$,
thus obtaining ${\cal{D}}_{bf} [\wm \times_{\pi} \wm, B]$ and then
declaring that $A \in \Psi_{b,\pi }^{-\infty, \delta}(\wm , \we)$
if there exists $\epsilon>0$ and
 $$\widetilde{A}\in\rho_{{\rm tot}}^{\delta+\epsilon}
H^\infty_{b,\, {\rm loc}} \bigl({\cal{D}}_{bf} [\wm \times_{\pi}
\wm, B] ; (E \boxtimes E^\ast)_{ \cal{D}} \bigr)$$ such that
$$\widetilde{A}|_{[\wm \times_{\pi} \wm, B]}=A\,.$$ In this
definition we have used the obvious inclusion $[\wm \times_{\pi}
\wm, B] \subset {\cal{D}}_{bf} [\wm \times_{\pi} \wm, B]$;
moreover we have denoted by $\rho_{{\rm tot}}$  a total boundary
defining function for ${\cal{D}}_{bf} [\wm \times_{\pi} \wm, B]$
and by $(E \boxtimes E^\ast)_{ \cal{D}}$ the bundle obtaned by
doubling $ E \boxtimes E^\ast$ across $bf$.

\smallskip

The
$\Gamma$-invariant elements 
in $\Psi^{*,\delta}_{b,\pi}(\wm,\we)$
will be denoted by $\Psi^{*,\delta}_{b,\rtimes}(\wm,\we)$;
thus
\begin{equation}\label{full-f}
\Psi^{*,\delta}_{b,\rtimes}(\wm,\we):=
\Psi^{*}_{b,\rtimes}(\wm,\we)+\Psi^{-\infty,\delta}_{b,\rtimes}(\wm,\we)+
\Psi^{-\infty,\delta}_{\rtimes}(\wm,\we)
\end{equation}

 The composition
rules for these calculi are as in the Appendix of \cite{MP I}
(Theorem 4); notice that since the fibres $\wm_\theta$ are
non-compact, for the composition of two operators
$P, Q\in \Psi^{*,\delta}_{b,\rtimes}(\wm,\we)$ to be
well defined it is necessary to assume that one of them is
compactly $\Gamma-$supported.
Because of our condition on the support, elements in
$\Psi^{*,\delta}_{b,\rtimes,c}(\wm,\we)$ can be composed exactly
as in Theorem 4 of (\cite{MP I},Appendix).

The mapping properties of elements in $\Psi^m_{b,\rtimes}(\wm,\we)$
are obtained by extending to the present context the arguments
in \cite{Melrose}. First of all, by proceeding
as in Lemma \ref{D} one can check that an
 element
in $\Psi^{m,\delta}_{b,\rtimes}(\wm,\we)$
defines a $C^\infty_c (T\rtimes \Gamma)$-linear operator
from  $C^\infty_c(\wm,\we)$ into $C^\infty(\wm,\we)$.
Moreover, the following Proposition (with
proof as in \cite{Melrose}) holds
\begin{proposition}\label{bounded}
Let $P\in \Psi^{m,\delta}_{b,\rtimes,c}(\wm,\we)$, $\delta>0$
Then for each $k\in\ZZ$, $P$ defines a bounded operator
\begin{equation}\label{bounded-f}
P:H^k_{b,C^0(T)\rtimes_r \Gamma}(\wm,\we)
\longrightarrow H^{k-m}_{b,C^0(T)\rtimes_r \Gamma}(\wm,\we)
\end{equation}
of
Hilbert $C^0(T)\rtimes_r \Gamma$-modules.
\end{proposition}


\subsection{Elliptic elements and $b$-parametrices}\label{parametrices}$\;$

\medskip
The goal of this subsection is to construct a parametrix
associated to a  $\Gamma$-equivariant family of odd Dirac
operators $D\in \Diff^1_{b,\rtimes}(\wm,\we)$ under the assumption
that the boundary family $D_0$ satisfies Hypothesis A of Section 2
with $\widehat{N}=\pa\wm$. The proof proceeds along the lines of
the parametrix construction given in \cite{Melrose}. See also
\cite{MP I} and \cite{LPMEMOIRS}.

\smallskip
\n
First notice that $D$ admits a symbolic  parametrix $Q_\sigma\in
\Psi^{-1}_{b,\rtimes,c}(\wm,\we)$, with $Q_\sigma$ odd; $Q_\sigma$
is obtained by proceeding as in \cite{CoLNM}.
If we write
$$Q=
\begin{pmatrix} 0 & Q^- \cr Q^+ &0  \cr
\end{pmatrix}
$$
we have
\begin{equation}\label{simbpar}
Q^+ D^+=\Id-S_{+,\sigma}\,,\quad\quad D^+ Q^+=\Id - S_{-,\sigma}
\end{equation}
with rests $S_{-,\sigma}$ and $S_{+,\sigma}$ in
$\Psi^{-\infty}_{b,\rtimes,c}(\wm,\we)$. Hypothesis A implies the
existence of the $L^2$-inverse of the indicial family of
$D^+(\theta)$,  $\forall \theta\in T$, i.e. the existence of the
$\Gamma$-equivariant family of operators $\{(D_0
(\theta)+i\lambda)^{-1}\}_{\theta\in T}$. We obtain in this way an
element $(D_0+i\lambda)^{-1}\,\in
\Psi^{-1}_{\rtimes}(\partial\wm,\we_{|_{\pa\wm}})$.
 By Proposition \ref{invertibility} we know that the family
$(D_0+i\lambda)^{-1}$ extends  for each $\lambda\in\RR$ to  a
bounded operator $(\D_0+i\lambda)^{-1}$ on $L^2_{C(T)\rtimes_r
\Gamma}(\partial\wm,\we_{|_{\pa\wm}})$; moreover,
$(D_0+i\lambda)^{-1} = A_\lambda + B_\lambda$ where $A_\lambda \in
\Psi^{-1}_{\rtimes,c}(\partial\wm,\we_{|_{\pa\wm}})$ uniformly in
$\lambda$ (i.e. the Schwartz kernel of $A_\lambda$ is included in
a fixed compact $K\subset  (\pa\wm\times_\pi \pa\wm)/\Gamma$) and
where $B_\lambda$ is given by a smooth Schwartz kernel and extends
to a bounded operator $\mathcal{B}_\lambda$ from $
H^k_{C(T)\rtimes_r \Gamma}(\partial\wm,\we_{|_{\pa\wm}})$ into $
H^l_{C(T)\rtimes_r \Gamma}(\partial\wm,\we_{|_{\pa\wm}})$ for any
$k,l \in \ZZ$. Let us consider the indicial family $$\RR\ni
\lambda\rightarrow (D_0+i\lambda)^{-1} \circ
I(S_{-,\sigma},\lambda)=A_\lambda \circ I(S_{-,\sigma},\lambda) +
B_\lambda \circ I(S_{-,\sigma},\lambda) \,,$$ which is well
defined since $S_{-,\sigma} \in
\Psi^{-\infty}_{b,\rtimes,c}(\wm,\we)$. We can take the inverse
Mellin transform of this indicial family and  obtain an element in
$\Psi^{-\infty,\delta}_{b,\rtimes,I}(\overline{N_+\partial\wm},\we_{|_{\pa\wm}})$;
we now consider the operator $Q^\prime\in
\Psi^{-\infty,\delta}_{b,\rtimes}$ corresponding to this operator
via the analogue of (\ref{indicial}) for the calculus with bounds.
The operator $Q^+ := Q^+_\sigma - Q^\prime$ belongs to
$\Psi^{-1,\delta}_{b,\rtimes}(\wm; \we^+,\we^-)$ and defines a
right inverse modulo $\rho_{{\rm
bf}}\Psi^{-\infty,\delta}_{b,\rtimes}$; this we call  a right
$b$-parametrix. Notice that
\begin{equation}\label{residual}
\rho_{{\rm bf}}\Psi^{-\infty,\delta}_{b,\rtimes}(\wm,\we)
\subset
\Psi^{-\infty,\delta}_{\rtimes}(\wm,\we)
\end{equation}
Thus $Q^-$ provides an inverse of
$D^+$ modulo $\Psi^{-\infty,\delta}_{\rtimes}(\wm,\we)$.
Standard arguments (\cite{Melrose}, p. 185)
show that this is also a left $b$-parametrix.
Thus: if Hypothesis A holds, then there exists $Q^+
\in \Psi^{-1,\delta}_{b,\rtimes}$ such that
\begin{equation}\label{pseudoinv}
D^+ Q^-=\Id-S_-\,,\quad Q^- D^+=\Id-S_+\,,\;\;\;\text{with}\;\;\;
S_-, S_+\in \Psi^{-\infty,\delta}_{\rtimes}
\end{equation} More information about $Q^\prime$ can be found in the next
subsection.

If $P\in\Diff^m_{b,\rtimes}(\wm,\we)$ is more generally an elliptic differential
$\Gamma$-invariant family, then we can find a
symbolic parametrix $Q_\sigma\in \Psi^{-m}_{b,\rtimes,c}(\wm,\we)$,
 see \cite{Co}.
If, in addition, the indicial family of $P$, $I_{\nu}(P,\lambda)$,
satisfies Hypothesis A
uniformly in $\lambda$, then we can
proceed as above and find an inverse modulo
$\rho_{{\rm bf}}\Psi^{-\infty,\delta}_{b,\rtimes}(\wm,\we)\subset
\Psi^{-\infty,\delta}_{\rtimes}(\wm,\we)$.

\s\n
{\it Summarizing}, we have shown that an elliptic element
in ${\rm Diff}^m_{b,\rtimes}(\wm,\we)$ with invertible
indicial family,
admits an
inverse $Q\in \Psi^{-m,\delta}_{b,\rtimes}$
modulo elements in $\Psi^{-\infty,\delta}_{\rtimes}$:
\begin{equation}\label{pseudoinv}
PQ=\Id-R_{0}\,,\quad\quad QP=\Id-R_{1}\,,\;\;\;\text{with}\;\;\;
R_{0}, R_{1}\in \Psi^{-\infty,\delta}_{\rtimes}
\end{equation}

\subsection{The $b$-index class}\label{bindexclass} $\;$

\m

The pseudoinverse $Q^+\in \Psi^{-m,\delta}_{b,\rtimes}$
associated to $D^+$ constructed in the previous subsection
will {\it not} be, in general,
of compact $\Gamma$-support. The problem comes from
 the correction term
$Q^\prime$, which  involves the {\it inverse} of the indicial family.
Indeed, the Schwartz kernel of $Q^\prime$ in the case
for example of Dirac operators, has the following
expression in a neighbourhood of the front face and in projective
coordinates $(s,y,y^\prime)$: (recall that $s=\frac {x} {x^\prime}$)
\begin{equation} \label{formula}
K(Q^\prime)(x,s,y,y^\prime)=
\phi(x) \int_\RR s^{ i \lambda} K(\,(D_0+i \lambda)^{-1} \circ
I(S_{\sigma ,-}, \lambda) \,) (y, y^\prime) d \lambda
\end{equation}
where $\phi$ is a given smooth function in $C^\infty([0,1] ; [0,1] )$
such that $\phi(x)=1$ for $0\leq x \leq \frac {1} {2}$ and
$\phi(x)=0$ for $x\geq \frac {3} {4}$ (recall that
$K(Q^\prime)$ vanishes identically outside a neighborhood
of the front face).
As a consequence, the rests $S_+, S_-$ (where
$S_-= S_{\sigma ,-}-D^+ Q^\prime$) are not of compact
$\Gamma$-support. In particular, we cannot conclude
that $S_+, S_-$ define $(C(T)\rtimes_r \Gamma)$-compact
operators: more work is  needed in
order to show that this is indeed the case and that
there is a well defined index class.

\n
We state first the result and we devote the next
subsection to its proof:

\begin{theorem}\label{c-index}
\item{1]} Let $D=(D(\theta)_{\theta\in T})$ a $\Gamma$-invariant
family of odd $\ZZ_2$-graded Dirac operators. If Hypothesis A
holds for the boundary family, then $D^+=(D^+(\theta))$ defines
$\forall m \in \NN^*$ a $C^0(T) \rtimes_r \Gamma$-linear bounded
operator $$ H^m_{b, C^0(T) \rtimes_r \Gamma}(\wm,\widehat{E}^+)
\longrightarrow  H^{m-1}_{b, C^0(T) \rtimes_r
\Gamma}(\wm,\widehat{E}^-)$$ which is invertible modulo $C^0(T)
\rtimes_r \Gamma $-compacts.
\item{2]} There is a well defined  index class $\Ind (D^+)$ in
$K_0(C^0(T) \rtimes_r \Gamma)$.
\end{theorem}

\subsection{Proof of Theorem \ref{c-index}}\label{proofthc-index} $\;$

\medskip
In order to simplify the notation we assume that $\we$ is
the product line bundle $\wm\times\CC\to \wm$.

\begin{lemma}\label{first}
Let $Q=Q_\sigma-Q^\prime$ be the parametrix constructed
in the previous section, with
$Q_\sigma$ of compact $\Gamma$-support. Then $\forall m\in\NN$,
$Q$
extends to a  bounded operator from $H^m_{b,C^0(T)\rtimes_r \Gamma}(\wm, \we)$
to $H^{m+1}_{b,C^0(T)\rtimes_r \Gamma}(\wm, \we)$.
\end{lemma}

\begin{proof} Since the Schwartz kernel of
$Q_\sigma$ is of compact $\Gamma$-support, we know that
$Q_\sigma$  extends to a  bounded operator from
$H^m_{b,C^0(T)\rtimes_r \Gamma}(\wm, \we) $
to $H^{m+1}_{b,C^0(T)\rtimes_r \Gamma}(\wm, \we)$.
Using Proposition \ref{invertibility} and the fact
that the schwartz kernel of $I(S_{\sigma, -}, \lambda )$ is
of compact $\Gamma$-support, one checks in a straightforward
way that for each
$\lambda \in \RR$:
$$
(D_0+i \lambda)^{-1} \circ I(S_{\sigma, -}, \lambda )
$$ is bounded from $H^m_{C^0(T)\rtimes_r \Gamma}(\partial \wm,
\we^-_{|\partial \wm} ) $
to $H^{m+1}_{C^0(T)\rtimes_r \Gamma}(\partial \wm,
\we^+_{|\partial \wm} )$.
Next we observe that near the boundary any element
of $L^2_{b,C^0(T)\rtimes_r \Gamma}(\wm, \we)$ may be
written as an element of
$L^2_b( [0,1], L^2_{C^0(T)\rtimes_r \Gamma}( \partial \wm, \we_{\partial \wm})
$ and that for $m\in \NN^*$ any element of
$H^m_{b,C^0(T)\rtimes_r \Gamma}(\wm, \we)$ maybe written as an element of
$$
L^2_b( [0,1], H^m_{C^0(T)\rtimes_r \Gamma}( \partial \wm, \we_{\partial \wm}))
+
H^1_b( [0,1], H^{m-1}_{C^0(T)\rtimes_r \Gamma}( \partial \wm, \we_{\partial
\wm})).
$$ Then using formula (\ref{formula}) one gets easily the Lemma.
\end{proof}

\begin{lemma}\label{second}
 The Schwartz kernel of $S_-$ satisfies
$K(S_-)(z, z^\prime)= \rho_{lb}^\epsilon \rho_{rb}^\epsilon
K(R^\prime)(z, z^\prime)$ for a suitable $\epsilon >0$
where  $R^\prime$ induces a bounded operator from
$H^m_{b,C^0(T)\rtimes_r \Gamma}(\wm, \we)$ to
$H^{m+1}_{b,C^0(T)\rtimes_r \Gamma}(\wm, \we)$ for
any $m\in \NN$.

\end{lemma}

\begin{proof} Recall that
$S_-= S_{\sigma, -} -D^+ Q^\prime$. Let
$\psi(z,z^\prime)$ be a smooth real valued function
on $\wm^2$ such that $\psi(z,z^\prime) = 1$ when both
$x=x(z)$ and $x^\prime=x(z^\prime)$ belong
to $[0, \frac {1} {100}]$ and $\psi(z,z^\prime) = 0$
when $x$ or $x^\prime$ is greater than $1$.
The required result is satisfied by the Schwartz kernel
$ (1-\psi(z,z^\prime)) K(S_-)(z, z^\prime) $.
So we are left to prove the Lemma for
$\psi(z,z^\prime) K(S_-)(z, z^\prime).$
With a small abuse of notation we identify
$K(S_-)$ with its lift to the stretched product.
 Using formula
(\ref{formula}) and employing
projective coordinates $(x,s,y, y^\prime)$
on the stretched product we have
$$
(\beta^*\psi) K(S_-)(x,s,y, y^\prime) =
(\beta^*\psi) K(S_{\sigma,-})(x,s,y, y^\prime)
-
$$
$$(\beta^*\psi) x\partial_x \phi(x)
\int_\RR s^{ i \lambda} K(\,(D_0+i \lambda)^{-1} \circ
I(S_{\sigma ,-}, \lambda) \,) (y, y^\prime) d \lambda
-
$$
$$(\beta^*\psi) \phi(x) \int_\RR s^{ i \lambda}  K(\,
I(S_{\sigma ,-}, \lambda) \,) (y, y^\prime) d \lambda.
$$
with $\beta$ equal to the blow-down map.

By construction, the indicial family of this operator
vanishes identically. Using Proposition \ref{invertibility} one
checks that for each real $\lambda$
the operators $ (D_0-i \lambda)^{-1} \circ
I(S_{\sigma ,-}, \lambda) $
belong to
$ B( \,H^m_{C^0(T)\rtimes_r \Gamma}(\partial \wm,
\we_{|\partial \wm}) ; H^{m+1}_{C^0(T)\rtimes_r \Gamma}(\partial \wm,
\we_{|\partial \wm}) \, ) $.
Following \cite{Melrose} (see Section 5.13 and Section 5.14 there)
one then checks that there exists $\epsilon >0$ such that one can write:
$$
\psi(z,z^\prime) K(S_-)(z, z^\prime)=\rho_{lb}^\epsilon \rho_{rb}^\epsilon
K(\beta_* R^\prime)(z, z^\prime)
$$ with  $(x,s) \rightarrow K(R^\prime)(x,s, y, y^\prime)$ induces
an element of
$$
H^\infty_b ( [0, 1]^2 ; B( \,H^m_{C^0(T)\rtimes_r \Gamma}(\partial \wm,
\we_{|\partial \wm}) ; H^{m+1}_{C^0(T)\rtimes_r \Gamma}(\partial \wm,
\we_{|\partial \wm}) \, ).
$$ One then gets easily the Lemma.
\end{proof}

\smallskip
Finally, one can show, proceeding as in
\cite{LPMEMOIRS} (Lemma 11.2) that the following
lemma holds:

\begin{lemma}\label{inclusion}
For each $\epsilon>0$, each $m\geq 0$ and each $k>0$ the inclusion
\begin{equation}\label{inclusion-f}
x^\epsilon H^{m+k}_{b, C^0(T) \rtimes_r \Gamma}(\wm,\widehat{E})
\hookrightarrow H^{m}_{b, C^0(T) \rtimes_r \Gamma}(\wm,\widehat{E})
\end{equation}
defines a $C^0(T) \rtimes_r \Gamma$-compact operator.
\end{lemma}

\smallskip
\n
{\bf End of the proof of Theorem \ref{c-index}.} It is an immediate consequence
of the previous  lemmas. More precisely, we can consider
\begin{equation}\label{bw}
{\rm IND}_m(D^+)=[p-p_0]\in K_0 (\KK (H^m_{b, C^0(T) \rtimes_r
\Gamma}))
\end{equation}
where $$p:= \begin{pmatrix} S_+^2 & {S_+(I+S_+)Q^+} \cr {S_- D^+} &
I-S_-^2 \cr
\end{pmatrix}\,,\quad p_0:= \begin{pmatrix} 0 & 0
\cr 0  & I \cr
\end{pmatrix} $$
and where $Q^+$ is a parametrix as above with rests $S_{\pm}$.
These index classes are all compatible through the natural
isomorphisms $$ K_0 (\KK (H^m_{b, C^0(T) \rtimes_r \Gamma}))\simeq
K_0 (\KK (H^\ell_{b, C^0(T) \rtimes_r \Gamma})). $$ Our index
class $\Ind (D^+)\in K_0 (C^0(T) \rtimes_r \Gamma)$ is defined as
the image of  any of these index classes under the natural
isomorphism $$ K_0 (\KK (H^m_{b, C^0(T) \rtimes_r \Gamma}))\simeq
K_0 (C^0(T) \rtimes_r \Gamma)\,.$$ We observe that Connes' index
class ${\rm IND}(D^+)\in K_0 (C^*_r (M,\cal F, E))$ is obtained by
using the K-Theory isomorphism $ K_0 (\KK (H^m_{b, C^0(T)
\rtimes_r \Gamma}))\simeq K_0 (C^*_r (M,\cal F, E))$ induced by
the isomorphism of algebras $\KK (H^m_{b, C^0(T) \rtimes_r
\Gamma})\simeq C^*_r (M,\cal F, E)$, see subsection
\ref{foliations}.

\medskip
\n
{\bf Remark.} It would be interesting to know whether the index class
we have just constructed can be recovered using  a
pseudodifferential calculus {\it on} a suitable groupoid. Recall
that for a single manifold with boundary, this is indeed possible.
See for example \cite{Monthubert}, \cite{NWX}.

\section{{\bf Rapid decay}}\label{rapid}

\subsection{Virtually nilpotent groups and the rapidly decreasing algebra}\label{rapid1}$\;$

\medskip
In the rest of the paper we shall make the following hypothesis:

\medskip
\noindent {\bf Hypothesis B.} {\it The group $\Gamma$ is virtually
nilpotent.}

\medskip

Following \cite{Go-Lo} we  shall now introduce  a {\it dense
subalgebra} ${\cal T}^\infty$ of $C^0(T) \rtimes_r \Gamma$  with
$$C^\infty(T)\rtimes \Gamma\subset {\cal T}^\infty\subset C^0(T)
\rtimes_r \Gamma\,.$$ Let us fix a word-metric on $\Gamma$,
$||\,\,||$, and let $$\Bi_\Gamma = \{f:\Gamma\rightarrow
\CC\,\,|\,\, \forall L\in\NN\,,\,
\sup_{\gamma\in\Gamma}(1+||\gamma||)^L |f(\gamma)|<\infty\}$$ be
the rapidly decreasing algebra in $C^*_r(\Gamma)$, a Fr\'echet
algebra.

\begin{definition} We define ${\cal T}^\infty$ as
$${\cal T}^\infty := C^\infty (T,\Bi_\Gamma)
\,.$$
\end{definition}

Elements in $\Ti$ can be written as $\sum t_\gamma \gamma$ where
the sum is now infinite but with the functions $t_\gamma$
satisfying the condition
\begin{equation}\label{decrease}
{\rm sup}_{\theta\in T,\gamma\in \Gamma}\left[ |t_\gamma |(\theta)
(1+||\gamma||)^N \right]<\infty,  \;\;\forall N\in\NN
\end{equation}
together with all their covariant derivatives.

\begin{definition} \label{decay-E}  Let $d$ be the distance on $\wm$ associated
with the lift on $\wm$ of an ordinary metric $g$ on $M$.
 Fix $z_0
\in \wm$. Then $C^\infty_{{\cal T }^\infty} (\wm,\widehat{E} )$ is
defined to be
 the subset of the elements $s \in C^\infty(\widehat{M}, \widehat{E})$
such that for any $N\in \NN$ and multi-index $\alpha$: $$
\sup_{z\in \wm} \bigl( ||\nabla^\alpha s(z)|| (1+d(z,z_0))^N
\bigr) < +\infty $$ where $\nabla$ denotes a $\Gamma-$invariant covariant derivative.
The space $C^\infty_{{\cal T }^\infty}
(\partial \wm,\widehat{E}_{|\partial \wm} )$ is defined similarly.
\end{definition}

\begin{lemma}\label{timodules}
Both $C^\infty_{{\cal T}^\infty}
  (\wm,\widehat{E} )$ and $ C^\infty_{{\cal T }^\infty}
(\pa\wm,\widehat{E}_{|\partial \wm } )$ are left $\Ti-$modules.
\end{lemma}

\begin{proof}
We assume that $\we$ is trivial. Let $s\in C^\infty_{\Ti}
(\wm,\we)$ and let $t=\sum t_\gamma \gamma\in \Ti$. We define the
module structure as follows: $$(t\cdot s)(z):=
\sum_{\gamma\in\Gamma} t_\gamma (\pi(z)) (\gamma\cdot s) (z)\,.$$
We need to check that
\begin{equation}\label{123}
\forall N\in \NN\,,\,\,\sup_{z\in \wm} |(t\cdot s)(z) | (1+
d(z,z_0))^N\,<\,\infty
\end{equation}
and similarly for the covariant derivatives. In order to check
(\ref{123})  we fix one $p\in \NN$ such that $\sum
(1+\|\gamma\|)^{-p}\,<\,\infty$. As $\Gamma$ is of polynomial
growth such $p$ always exists. Let $N \in \NN$, then fix $C_N>0$
such that: $$ {\rm for}\; ||\gamma || \leq \frac {d (z,z_0)}
{2},\; |s(z \cdot \gamma)| \leq C_N ( 1 +d(z,z_0))^{-N-p}, $$ $$
\sup_{z \in \wm}|s(z)| \leq C_N, \; \sup_{\theta \in T, \gamma \in
\Gamma} |t_\gamma(\theta)| \, \leq C_N$$ $$ \forall \gamma \in
\Gamma,\; \sup_{\theta \in T} |t_\gamma(\theta)| \leq C_N (1 + ||
\gamma||)^{-N-p}. $$ Then one has: $$ \sum_{\|\gamma\|\geq
\frac{d(z,z_0)}{2}} |t_\gamma (\pi (z))| |\gamma\cdot s (z)|
(1+d(z,z_0))^N\,\leq\, \sum_{\|\gamma\|\geq \frac{d(z,z_0)}{2}}
C_N^2 (1+\|\gamma\|)^{-p} \frac{(1+d(z,z_0))^N}{(1+\|\gamma\|)^N}
$$
$$\leq C^2_N 2^N \sum_{\gamma\in \Gamma} (1+\|\gamma\|)^{-p} <
\infty $$ Moreover
$$\sum_{\|\gamma\|\leq
\frac{d(z,z_0)}{2}} |t_\gamma (\pi (z))| |\gamma\cdot s (z)|
(1+d(z,z_0))^N\,\leq \,\sum_{\|\gamma\|\leq \frac{d(z,z_0)}{2}}
C_N^2 (1+d(z,z_0))^{-p}
\leq
$$
$$\sum_{\gamma\in \Gamma} C_N^2
(\frac{1}{1+2\|\gamma\|})^p \,<\,\infty
$$

This proves (\ref{123}); for the covariant derivatives the
argument is similar. The proof for $ C^\infty_{{\cal T }^\infty}
(\pa\wm,\widehat{E}_{|\partial \wm } )$ is the same.
The lemma is proved.
\end{proof}

\subsection{The refined $b$-index class}\label{rapid2}$\;$

\medskip
Once the rapidly decreasing algebra $\Ti$ is fixed, we can define
the rapidly decreasing $b$-calculus with bounds
$\Psi^{*,\delta}_{b,\Ti}$ where
$$\Psi^{*,\delta}_{b,\rtimes,c}\subset
\Psi^{*,\delta}_{b,\Ti}\subset \Psi^{*,\delta}_{b,\rtimes}.$$
Operators in $\Psi^{*,\delta}_{b,\Ti}$ are not of compact
$\Gamma$-support but have precise asymptotic properties with
respect to the $\Gamma$-action. The  definition is somewhat
technical and can be found in Appendix A (Section
\ref{rapidcalculus}).

Operators in the rapidly decreasing calculus have natural mapping
properties. First of all, an element in $\Psi^*_{b,\Ti}$ defines a
$\Ti$-linear endomorphism of $C^\infty_{\Ti}(\wm,\we)$. Moreover
the following proposition holds:

\begin{proposition}  If $\Gamma$ is virtually nilpotent and
$P\in\Psi^m_{b,\Ti}$, then $P:H^{k+m}_{b,C^0 (T)\rtimes_r
\Gamma}(\wm,\we) \rightarrow H^k_{b,C^0 (T)\rtimes_r
\Gamma}(\wm,\we)$ is bounded for each $k\in\ZZ$.
\end{proposition}

Next we consider $K\in \Psi^{-\infty, \epsilon}_{\Ti}(\wm,\we)
\cup \Psi^{-\infty, \epsilon}_{b, \Ti}(\wm,\we) $. Using the fact
that $K$ is a $\Gamma-$equivariant family and that the estimates
(\ref{es}) in Appendix A hold, one checks easily that $K$ defines
a $\Ti$-linear endomorphism of $H^\infty_{b, \Ti}(\wm,\we)$.
Moreover, if $K\in \Psi^{-\infty, \epsilon}_{\Ti}(\wm,\we)$
 then $K$ sends $H^\infty_{b, \Ti}(\wm,\we)$ into
$\rho^\epsilon H^\infty_{b, \Ti}(\wm,\we)$.

Finally, proceeding  as in \cite{Melrose} (see Section 5.16) and
as in \cite{MP I} (Theorem 4), one proves the following
composition rules.

\begin{proposition} For any $\epsilon > 0$, the spaces
$ \Psi^{-\infty, \epsilon}_{b, \Ti}(\wm,\we)$
 and $\Psi^{-\infty, \epsilon}_{\Ti}(\wm,\we) $ are two sided modules over
the small calculus $\Psi^{*}_{b, \Ti}(\wm,\we) $; moreover
$\Psi^{-\infty, \epsilon}_{\Ti}(\wm,\we) $ is a two sided-module
over $\Psi^{-\infty, \epsilon}_{b, \Ti}(\wm,\we) $. Both
$\Psi^{*}_{b, \Ti}(\wm,\we) $ and $ \Psi^{-\infty,
\epsilon}_{\Ti}(\wm,\we)$ are algebras and we have: $$
\Psi^{-\infty, \epsilon}_{b,\Ti}(\wm,\we) \circ \Psi^{-\infty,
\epsilon}_{b,\Ti}(\wm,\we) = \Psi^{-\infty,
\epsilon}_{b,\Ti}(\wm,\we) + \Psi^{-\infty,
\epsilon}_{\Ti}(\wm,\we). $$
\end{proposition}

We shall need the following

\begin{proposition}\label{pararapid}
Let $D^+=(D^+ (\theta))_{\theta\in T}$ a $\Gamma$-equivariant family
with boundary operator $D_0=(D_0(\theta))_{\theta\in T}$
satisfying Hypothesis A. Then there exists $Q^+\in\Psi^{-1}_{b,\Ti}$
such that
\begin{equation}\label{pararapid-f}
Q^+ D^+=\Id-S^+\,,\quad D^+ Q^+=\Id-S^-\,,\quad\text{with}\quad
S^\pm\in  \Psi^{-\infty, \epsilon}_{\Ti}
\end{equation}
\end{proposition}

\begin{proof}
We just need to observe that the correction term in the
construction of the true parametrix (see (\ref{formula})) is
rapidly decreasing. However this is an immediate consequence of
finite propagation speed estimates, exactly as in \cite{LPMEMOIRS}
Proposition 1.5 and Proposition 9.1.
\end{proof}

We are now in the position of
defining a {\it refined} index class.

\begin{definition}\label{refined}
Let
$D^+=(D^+  (\theta))_{\theta\in T}$ be a $\Gamma$-equivariant
family
with boundary operator $D_0=(D_0(\theta))_{\theta\in T}$
satisfying Hypothesis A. We define its index class
in $K_0 ( \Psi^{-\infty, \epsilon}_{\Ti}(\wm,\we))$ by
\begin{equation}\label{refined-f}
{\rm Ind}_{\Ti}(D^+)=[p-p_0]\in K_0 ( \Psi^{-\infty, \epsilon}_{\Ti}(\wm,\we))
\end{equation}
where $$p:= \begin{pmatrix} S_+^2 & {S_+(I+S_+)Q^+} \cr {S_- D^+} &
I-S_-^2 \cr
\end{pmatrix}\,,\quad p_0:= \begin{pmatrix} 0 & 0
\cr 0  & I \cr
\end{pmatrix} $$
\end{definition}

\m
\n
{\bf Remark.} There is an injection $$J:
\Psi^{-\infty, \epsilon}_{\Ti}(\wm,\we)\hookrightarrow
\KK (L^2_{b, C^0 (T)\rtimes_r \Gamma}(\wm,\we) )$$
and the index class defined in subsection \ref{bindexclass}
is simply the image of ${\rm Ind}_{\Ti}(D^+)$ under the morphism induced
by $J$ in K-theory.

\m
\n
{\bf Remark.}
Let $\Gamma$ be an arbitrary finitely generated group.
In the closed case Gorokhovsky and Lott  defines
a  subalgebra $C^\infty (\Gamma, \mathcal{B}^\omega)$
of $C^0 (T)\rtimes \Gamma$ by imposing an exponentially
rapid decay on the coefficients $t_\gamma$ appearing
in $\sum t_\gamma \gamma$. They can then construct
the algebra of rapidly degreasing smoothing operators
and prove the existence of a refined index class in
the K-theory group of such an algebra.

One might wonder why this program cannot be implemented
on a manifold with boundary.
The reason is once again in the construction of the true
parametrix and more precisely in formula (\ref{formula}),
where the resolvent family of the boundary operator appears.
One needs a notion of rapidly decreasing pseudodifferential
operators
containing the resolvent of the boundary operator,
uniformly in $\lambda$.
 At the moment it is only for the  groups of
polynomial growth that we know how to
define such a notion.

\bigskip
\section{{\bf Noncommutative differential forms  and
higher eta invariants}}\label{heta}

\subsection{Noncommutative differential forms}\label{heta1}$\;$

\medskip
We consider  Lott's space of noncommutative differential forms
$\Omega_\ast (\cal {B}^\infty_\Gamma )$  (see \cite{Lott II} or
\cite{LPMEMOIRS} page 22 for a definition).
We then make the following
\begin{definition} For any non negative integers $k$ and $l$ we set:
$$
\widehat {\Omega}_{k,l}(T, \Bi_\Gamma)=
\Omega^k(T) \widehat{\otimes} \Omega_l (\cal {B}^\infty_\Gamma )
$$ where $\Omega^k(T)$ denotes the Frechet space of smooth $k-$differential forms
over $T$ and $\widehat{\otimes} $ denotes a
complete projective graded tensor product. We also set:
$$\widehat {\Omega}_{\ast}(T, \Bi_\Gamma)=
\Pi_{k,l \in \NN} \widehat {\Omega}_{k,l}(T, \Bi_\Gamma).
$$
\end{definition}

Elements in  $\widehat {\Omega}_{k,l}(T, \Bi_\Gamma)$ are given by
sums  $$\sum \alpha_{g_0,g_1,\cdots,g_\ell}(\theta)g_0 dg_1 \cdots
dg_\ell$$ with $\alpha_{g_0,g_1,\cdots,g_\ell}(\theta) $ $k$-forms
on $T$ which are, in addition, rapidly decreasing, together with
their $\theta$-derivatives, with respect to $\| g_0 \|+\cdots
+\|g_k \|$.  Using Lemma \ref{decrease} and the fact that $\Gamma$
is virtually nilpotent one shows that $\widehat {\Omega}_{0}(T,
\Bi_\Gamma)$ may be identified with $\Ti$.

We let an  element $g$ of the group $\Gamma$ act on a form
$\omega\in \Omega^*(T)$ by considering the pullback via the
diffeomorphism defined by the action of  $g$. We observe that
$\widehat {\Omega}_{\ast}(T, \Bi_\Gamma)$ has a bi-module
structure over ${\cal T}^\infty$; this is described by the
following rules with $u $ denoting an element of $C^\infty(T)$ and
$u \gamma \in \Ti$:

$$(\,\alpha  \widehat{\otimes} g_0 dg_1\cdot \cdot \cdot dg_l\,)
\cdot u \gamma \,=\, (g_0\cdot \cdot \cdot g_l)^{\ast}(u) \alpha
\widehat{\otimes} g_0 dg_1\cdot \cdot \cdot dg_l \gamma $$
 $$u \gamma  \cdot (\,\alpha  \widehat{\otimes} g_0
dg_1\cdot \cdot \cdot dg_l\,) \,=\, u \gamma^\ast(\alpha)
\widehat{\otimes} \gamma g_0dg_1\cdot \cdot \cdot dg_l. $$

Moreover  $\widehat {\Omega}_{\ast}(T, \Bi_\Gamma)$ is
    a Fr\'echet  algebra whose product is defined by the following
rules:
\begin{equation}\label{forms}
(\,\alpha  \widehat{\otimes} g_0
dg_1\cdot \cdot \cdot dg_l\,)
\cdot (\,\beta  \widehat{\otimes} \gamma_0
d\gamma_1\cdot \cdot \cdot d\gamma_p\,) =
(-1)^{l.|\partial \beta|} \alpha \wedge
(g_0\cdot \cdot \cdot g_l)^{\ast}( \beta) \,  \widehat{\otimes} g_0
dg_1\cdot \cdot \cdot
dg_l      \gamma_0 d\gamma_1
\cdot \cdot \cdot d\gamma_p
\end{equation}
 where $\alpha$ and $\beta$ are homogeneous differentials forms on
$T$.  We then introduce, following \cite{Go-Lo}, the following
space  which will be the receptacle for the super-traces: $$
\overline{\widehat{\Omega}}_\ast(T, \Bi_\Gamma)=
\widehat{\Omega}_\ast(T, \Bi_\Gamma)  / \overline{
[\widehat{\Omega}_\ast(T, \Bi_\Gamma)\,,\,\widehat{\Omega}_\ast(T,
\Bi_\Gamma)]_t\;\;;} $$ this is a space of {\it noncommutative
differential forms modulo the closure of the space of graded
commutators.} We endow  $\overline{\widehat{\Omega}}_\ast(T,
\Bi_\Gamma)$ with the natural differential and we get a complex
whose homology is denoted $\widehat{H}_\ast(T,\Bi_\Gamma)$.

 This
homology does not coincide with
 the topological
noncommutative de Rham homology $\widehat{H}_\ast(\Ti)$ of $\Ti$
(see \cite{Karoubi}), simply because $\widehat{\Omega}_\ast(T,
\Bi_\Gamma)$ does not coincide with $\widehat{\Omega}_\ast(\Ti)$.
But, as we shall see later, there is a natural map from
$\widehat{\Omega}_\ast(\Ti)$ to $\widehat{\Omega}_\ast(T,
\Bi_\Gamma)$ .

\subsection{Operators with differential forms coefficients}\label{heta2}$\;$

\medskip
Following \cite{Go-Lo} (see also \cite{LPMEMOIRS} page 26), we
give the following

\begin{definition} \label{opdecay} We denote by
$$
\operatorname{Hom}^\infty_{{\cal T}^\infty}
(\,C^\infty_{{\cal T}^\infty} (\wn,\wf)\,;\,
\widehat{\Omega}_{k,l}(T,\Bi_\Gamma)\otimes_{{\cal T}^\infty}
C^\infty_{{\cal T}^\infty} (\wn,\wf)\,)
$$ the set  of left $\widehat{\Omega}_\ast(T, \Bi_\Gamma)-$linear operator $K$ acting on
$ \widehat{\Omega}_\ast(T, \Bi_\Gamma)\otimes_{\Ti} C^\infty_{{\cal T}^\infty} (\wn,\wf)$
and sending any
$F \in C^\infty_{{\cal T}^\infty} (\wn,\wf)$ to    $KF$ defined for any
$ z \in \wn$ by:
$$
(KF)(z)= \sum_{g_1,\dots , g_l\in \Gamma}
dg_1\cdot \cdot \cdot dg_l \int_{\pi^{-1}(z\cdot (g_1\cdot \cdot \cdot  g_l)^{-1})}
K_{g_1,\dots , g_l}(z,w) F(w)
dVol_{ \pi^{-1}( z\cdot (g_1\cdot \cdot \cdot  g_l)^{-1}) }(w)
$$
 where
$$K_{g_1,\dots , g_l}: (z,w) \in \wn\times \wn \rightarrow
\bigwedge^k(T^\ast_{\pi(z)}T) \otimes \Hom (\widehat{F}_w\,,\,
\widehat{F}_z)$$ vanishes for $\pi(z)\cdot (g_1\cdot\cdot\cdot
g_l)^{-1} \not= \pi(w)$ and defines a smooth section on $$\{ (z,w)
\in \widehat{ N} \times \widehat{ N} \;|\; \pi(z)\cdot (g_1\cdot
\cdot \cdot  g_l)^{-1} = \pi(w) \}.$$ Moreover, we require the
following decay estimates: for any fundamental domain $A$ of
$\widehat{ N} $, any  integer  $p>1$ and any differential operator $P \in {\rm
Op}\,(\wn\times_\pi \wn )$ $$ \sup_{(z,w) \in  \widehat{N}
\times_\pi  \widehat{ N} } \bigl[d(zg_1\cdot \cdot \cdot  g_l, A)
+ ||g_2||+ \cdot \cdot \cdot + ||g_{l}|| + d(w, Ag_1^{-1})
\bigr]^p |P_{z,w} K_{g_1,\cdot \cdot \cdot , g_l}(zg_1\cdot \cdot
\cdot g_l,w)| $$ is finite.
\end{definition}
 The  algebra of  differential operators
 ${\rm Op}\,(\wn\times_\pi \wn )$
 acting on $C^\infty(\wn\times_\pi \wn\,,\,
 \wf \boxtimes \wf^\star) $, which appears  above,  is
obtained by considering manifolds with empty boundary
in Definition
\ref{Op} in Appendix A.

\medskip
\n
{\bf Remarks.}

\n
 1) The $\Ti-$linearity is a strong assumption imposed on $K$
and on the $K_{g_1,\dots , g_l}$. It does not seem easy to
characterize the $\Ti-$linearity by simple formulas involving the
$K_{g_1,\dots , g_l}$ alone.  This fact of course  already appears
in the covering case ($T$ reduced to a point, see \cite{Lott I})
 when one deals with operators with noncommutative
differential forms coefficients. Let us simply observe that the
    condition on the support of the $K_{g_1,\dots , g_l}$
is an easy
consequence of the
rules
 (\ref{forms}) and of the $C^0(T)-$linearity of $K$.

 \n 2) Assuming $K$ to be $\Ti-$linear and
  proceeding
 as for Proposition
 2.3 of \cite{LPMEMOIRS}, one checks easily that Definition  \ref{decay-E} and
 the above decay estimate imply that $K$
indeed sends $C^\infty_{{\cal T}^\infty} (\wn,\wf)$ into
$\widehat{\Omega}_{k,l}(T,\Bi_\Gamma)\otimes_{{\cal T}^\infty}
C^\infty_{{\cal T}^\infty} (\wn,\wf)$

\n 3) For $k=l=0$, the operator $K$ induces a $C^0(T)\rtimes_r
\Gamma-$compact operator of $L^2_{C^0(T)\rtimes_r
\Gamma}(\wn,\wf)$.

\medskip
We fix  a function $\phi \in C^\infty_c(\widehat{ N})$ which
satisfies
$\sum_{\gamma \in \Gamma} \gamma\cdot\phi = 1$.
Let us denote the unit element of $\Gamma$ by $e$. For
any $\theta \in T$, one sets as in \cite{Go-Lo},   $ \STR_{<e>, Cl(1)} K (\theta) =$
\begin{equation} \label{strcl}
\sum_{g_0,\dots , g_l \in \Gamma,\, g_0\cdot \cdot \cdot g_l =e}
(dg_1\cdot \cdot \cdot dg_l ) g_0
\int_{\pi^{-1}(\theta)} \phi(w)  \Str_{Cl(1)}
K_{g_1,\cdot \cdot \cdot , g_l}(wg_0^{-1},w)
dVol_{\pi^{-1}(\theta)} (w)
\end{equation}
where the definition of the  supertrace $\Str_{Cl(1)}$
acting  on $(\wf_w)_\sigma=\wf_w \otimes {\rm Cl}(1)$ is recalled
in \cite{LPMEMOIRS} (page 21). Moreover,
 Definition \ref{opdecay} implies that
$w\rightarrow K_{g_1,\cdot \cdot \cdot , g_l}(wg_0^{-1},w) $ is
smooth on the fiber-diagonal in $\wn \times_\pi \wn$
and, under
the natural identification between
$\widehat{F}_w \simeq \widehat{F}_{wg_0^{-1}}$
(due to the $\Gamma$-invariance of $\wf$), that
$K_{g_1,\cdot \cdot \cdot , g_l}(wg_0^{-1},w) $ belongs to
$ \bigwedge^k(T^\ast_{\pi(wg_0^{-1})}T) \otimes \Hom (\widehat{F}_w\,,\,  \widehat{F}_w)$.

Notice that $ \STR_{<e>, Cl(1)} K \in \OM$ and that
 the definition of $ \STR_{<e>, Cl(1)} K $ extends obviously to
$$\Hom^\infty_{{\cal T}^\infty}(\,C^\infty_{{\cal T}^\infty} (\wn,\wf)\,;\,
\widehat{\Omega}_{\ast}(T,\Bi_\Gamma)
\otimes_{{\cal T}^\infty}C^\infty_{{\cal T}^\infty} (\wn,\wf)\,)
$$ which we define to be equal to
$$\Pi_{k,l \in \NN} \Hom^\infty_{{\cal T}^\infty}
(\,C^\infty_{{\cal T}^\infty} (\wn,\wf)\,;\,
\widehat{\Omega}_{k,l}(T,\Bi_\Gamma)\otimes_{{\cal T}^\infty}
C^\infty_{{\cal T}^\infty} (\wn,\wf)\,).
$$

\m
\n
{\bf Notation.} We shall henceforth use the following notation:
\begin{equation}\label{notation}
\Psi^{-\infty}_{\widehat{\Omega}_{\ast}(T,\Bi_\Gamma)}
(\wn,\wf)\,:=\,
\Hom^\infty_{{\cal T}^\infty}(\,C^\infty_{{\cal T}^\infty}
(\wn,\wf)\,;\, \widehat{\Omega}_{\ast}(T,\Bi_\Gamma)
\otimes_{{\cal T}^\infty}C^\infty_{{\cal T}^\infty}
(\wn,\wf)\,)
\end{equation}

\subsection{The higher eta invariant}\label{heta3} $\;$
\medskip

We fix a $\Gamma-$invariant horizontal distribution
$T^H\widehat{N}$ such that
$$
T\widehat{N} = T^H\widehat{N} \oplus T(\widehat{N}/T)
$$ and 
 consider the rescaled Bismut superconnection
in the odd-dimensional context (our notation
follows  \cite{MP I} or \cite{LPMEMOIRS}):
$$
\BB_s^{Bismut}=\sigma s \D_0+\nabla^{u}-\sigma {1\over 4 s}c(\tau),\;
s\in \RR^{+\ast}
$$ where $\D_0$ is the $\Ti$-linear Dirac operator introduced above,
$c(\tau) $ denotes the Clifford multiplication by the curvature
2-form $\tau$ of
$T^H\widehat{N} $ and  $\nabla^{u} $ is a certain unitary connection.

Then, following \cite{Go-Lo},  we
fix a function $h \in C^\infty_c(\widehat{M})$
such that $\sum_{\gamma \in \Gamma} \gamma\cdot h =1$
on $\widehat{N} $ and
 consider for each real $s>0$ the superconnection
$$\BB_s= \BB_s^{Bismut} +
\sum_{\gamma \in \Gamma} d\gamma \otimes h \gamma^{-1}
$$ which sends $C^\infty_{{\cal T}^\infty}(\widehat{N},\wf)$ into
$\widehat{\Omega}_{\ast}(T,\Bi_\Gamma)
\otimes_{{\cal T}^\infty} C^\infty_{{\cal T}^\infty}(\widehat{N},\wf)$. We
recall the following rule of computation valid for any
$\omega \otimes f \in  \widehat{\Omega}_{\ast}(T,\Bi_\Gamma)
\otimes_{{\cal T}^\infty} C^\infty_{{\cal T}^\infty}(\widehat{N},\wf)$:
$$
\BB_s(\omega \otimes f) = d\omega \otimes f + (-1)^{\partial \omega} \omega
\otimes \BB_s( f).
$$
In order to state a Lichnerowicz-type formula for $\BB_s^2$ we introduce
a few notations. Let $\{e_i\}_{i=1}^{2k-1}$ be a local
orthonormal basis for $T\pi^{-1}(\theta)$ and let
$\{c^i\}_{i=1}^{2k-1}$ be Clifford algebra generators, with
$(c^i )^2=-1$. Let $\{\tau^\alpha\}_{\alpha=1}^{{\rm dim}\, T}$ be a local
basis of $T^\ast T$ and let $E^\alpha$ denote exterior multiplication by
$\tau^\alpha$.

We then have the odd-dimensional analogue of formula (4.11) of \cite{Go-Lo} :
\begin{equation} \label{odd-lichne} \BB_s^2=
s^2D_0^2 + { 1\over 4} \sum_{i,j} \tilde{F}_{i,j}(\widehat{V})
[c^i,c^j] + \sum_{i,\alpha} \tilde{F}_{\alpha,i}(\widehat{V})\,
E^\alpha c^i + {1\over 4} \sum_{\alpha,\beta}
\tilde{F}_{\alpha,\beta}(\widehat{V}) [E^\alpha, E^\beta]
\end{equation} $$- s \sum_{\gamma \in \Gamma} d \gamma (\,
c(d^{vert} h) + {E}(d^{hor}h)\,) \gamma^{-1} -\sum_{\gamma,
\gamma^\prime} (\gamma\gamma^\prime \cdot h) (\gamma \cdot
h)d\gamma d\gamma^\prime (\gamma\gamma^\prime)^{-1} $$

Notice that  the $\Gamma-$invariance of the horizontal
distribution $T^H\wn$ implies that $d\BB_s/ ds$ is $\Ti$-linear
(i.e. ``vertical''). Now we define $\exp(-\BB^2_s)$  using Duhamel
formula around $\D_0^2$: $$ e^{-\BB^2_s} = e^{-s \D_0^2} +
\int_0^1 e^{-u_1 s^2 \D_0^2} (\D_0^2-\BB^2_s) e^{(-(1-u_1) s^2
\D_0^2} du_1 \,+ $$
$$ \int_0^1\int_0^{1-u_1} e^{-u_1 s^2
\D_0^2}(\D_0^2-\BB^2_s) e^{-u_2 s^2
\D_0^2}(\D_0^2-\BB^2_s) e^{-(1-u_1-u_2) s^2 \D_0^2} du_2 du_1+
\ldots $$
Now, since by
     Hypothesis A
$\D_0$ is invertible, we may apply to  $K(e^{-s^2D_0^2})(z,w)$
 the        finite propagation speed estimates
                           (see \cite{LPMEMOIRS} formula (2.14)). So, there
exists $\delta >0$ such that for any $a, N \in \NN$ one has:
$$
\forall (z,w) \in \wn \times_\pi \n\wn,\;
\forall s \geq 1,\; |( D_0^a e^{-s^2D_0^2})| ( z,w) \leq C(a,N) (1
+d(z,w) )^N \exp(-s^2\delta). $$ Of course,
$K(e^{-s^2D_0^2})(z,w)=0$ for $\pi(z)\not=\pi(w)$. Therefore,
for each $s>0$ $e^{-s^2 \D_0^2}$,  defines  an element in
$\Psi^{-\infty}_{{\cal T}^\infty} (\widehat{N}\,;\, \widehat{F})$

Using these estimates, one shows
as in  Section 3 of \cite{LPMEMOIRS}
 that
 the operators
  \begin{equation}\label{asinmem}
e^{-\BB^2_s}\;\;\text{and}\;\;
{d\BB_s \over ds}\,  e^{-\BB^2_s}\;\;\text{ both belongs to}\;\;
\Psi^{-\infty}_{\widehat{\Omega}_{\ast}(T,\Bi_\Gamma)} (\wn,\wf)
\end{equation}
uniformly with respect to $s\rightarrow +\infty$.

 Let $\R$ be the rescaling operator on $\OM$ which multiplies by
$(2 i \pi)^{-p}$ a form of degree $2p-1$ or of degree $2p$. We observe that
 $\R$ commutes with the differential $d$ of $\OM$. We
then introduce $$\widetilde{\eta}_{_{<e>}}(s)= \R(\STR_{<e>,Cl(1)}
(\,{d\BB_s \over ds}\,  e^{-\BB^2_s}\,))\,. $$

\medskip
\begin{proposition} \label{higheta} Under Hypothesis A and B,
the higher eta invariant
$$
\widetilde{\eta}_{_{<e>}} =
{2 \over \sqrt{\pi}} \int_0^{+\infty} \widetilde{\eta}_{_{<e>}}(s)\,ds\,\in\,\OM
$$ is well defined.
\end{proposition}

\begin{proof}  The large time behaviour is treated as
in \cite{LPMEMOIRS}, using (\ref{asinmem}).
Thus we just have to show the convergence of the
integral at $s=0$, but this is basically done as in the proof of
Proposition 26 page 219 of \cite{Lott II}, taking into account
theorem 2 in
\cite{Go-Lo}.
\end{proof}

\s\n {\bf Remark.}
In the covering case, the definition of higher eta invariant,
and the proof of the convergence of the integral defining
it,  are due to Lott
(see section 4 of  \cite{Lott II}).

\medskip Proceeding as in \cite{Lott II} page 220, one proves the following:
\begin{proposition} Consider a smooth path
of superconnections $\BB_s(r), 0 \leq r \leq 1$  associated to a
smooth path $D_0(r)$ ($0 \leq r \leq 1$) of Dirac type operators
satisfying Assumption A. We assume that $\BB(r)$ is obtained by
smoothly varying with respect to $r$, the metric on $N$, the
$\Gamma-$invariant metric in the fibers, the $\Gamma-$invariant
connection on  $\wf$, the horizontal distribution $T^H \wn$ and
the function $h \in C^\infty_c( N)$. Then one has the following
variational formula for the higher eta invariant
$\widetilde{\eta}_{_{<e>}}(r)$ associated to $\BB_s(r)$ as in
Proposition \ref{higheta}:
  $$
  \frac {d \widetilde{\eta}_{_{<e>}}(r)} {dr}\,=\,
 -{2 \over \sqrt{\pi}} \lim_{r\rightarrow 0^+}
\R\, \STR_{Cl(1)} (\,{d\BB_s(r) \over dr}\,  e^{-\BB^2_s(r)}\,)
  $$ in $\OM $ modulo exact forms.
\end{proposition}
\begin{remark} The local index theory for \'etale groupoids developed
in \cite{Go-Lo} (Section 4) insures the existence of the above limit.
\end{remark}

\section{{\bf The $b$-supertrace and the higher local index
theorem}}\label{bsuper}

\subsection{The $b$-supertrace of an element in
$\Psi^{-\infty, \delta}_{b,\widehat{\Omega}_{*}(T,\Bi_\Gamma)}
(\widehat{M}\,;\, \widehat{E})$}\label{bsuper1} $\,$

In the previous section we have introduced, in the closed case,
the space of rapidly decreasing smoothing operators with
differential form coefficients. We have denoted this space of
operators by $\Psi^{-\infty}_{\widehat{\Omega}_{*}(T,\Bi_\Gamma)}
(\widehat{N}\,;\, \widehat{F})$. Let us now pass to the case of
manifolds with boundary and to the notion of rapidly decreasing
$b$-smoothing operator with differential form coefficients. There
are as usual three spaces:
$$\Psi^{-\infty}_{b,\widehat{\Omega}_{\ast}(T,\Bi_\Gamma)}
(\widehat{M}\,;\, \widehat{E})\,,\quad
\Psi^{-\infty,\delta}_{b,\widehat{\Omega}_{\ast}(T,\Bi_\Gamma)}
(\widehat{M}\,;\, \widehat{E})\,\quad
\Psi^{-\infty,\delta}_{\widehat{\Omega}_{\ast}(T,\Bi_\Gamma)}
(\widehat{M}\,;\, \widehat{E})$$ depending on the boundary
behaviour of the coefficients $K_{g_1,\cdot \cdot \cdot , g_l}$ in
the Schwartz kernel of the operator $K:=\sum_{g_1,\cdot \cdot
\cdot , g_l\in \Gamma}dg_1\cdot \cdot \cdot dg_l \, K_{g_1,\cdot
\cdot \cdot , g_l}.$ We give the precise definition of these 3
spaces in Appendix B.

\medskip
We now proceed to the definition of the $b$-supertrace and its
properties. We fix a function $\phi \in C^\infty_c(\widehat{ M})$
which is constant in the normal direction near the boundary and
satisfies $\sum_{\gamma \in \Gamma} \gamma \cdot \phi = 1$.

As in Section \ref{bcalc}
 let us fix
a $\Gamma-$invariant trivialization $\nu \in
C^\infty(\partial\widehat{ M}, N_+\partial\widehat{ M})$ of the
normal bundle and $x \in C^\infty(\widehat{ M})$ a
$\Gamma-$invariant boundary defining function for
$\partial\widehat{ M}$ such that $dx\cdot\nu=1$ on
$\partial\widehat{ M}$. For any function $f \in C^\infty(\widehat{
M}\times_\pi \widehat{ M})$ we set for any $\theta \in T $, $$
\nuint_{\Delta_b(\theta)}\phi(w) f(w,w)\,
 dVol^b_{\pi^{-1}(\theta)}(w):=$$
$$ \lim_{\epsilon \rightarrow 0^+} [ \int_{x> \epsilon}
\phi(w) f(w,w)\, dVol^b_{\pi^{-1}(\theta)}(w) +\log \epsilon
\int_{\partial\pi^{-1}(\theta) } \phi(w) f(w,w)\,
dVol_{\partial \pi^{-1}(\theta)}(w) ]
$$
 where ${\Delta}_b (\theta) $ denotes the lifted diagonal of $\pi^{-1}(\theta)$ in
the associated $b-$streched product.

\medskip Now, for given $(k,l)\in \NN$ we consider a particular
element $K \in \Psi^{-\infty, \delta}_{b,\widehat{\Omega}_{\ast}(T,\Bi_\Gamma)}
(\widehat{M}\,;\, \widehat{E})$ thus
with Schwartz kernel  of the form
$$
K(z,w)= \sum_{g_1,\cdot \cdot \cdot , g_l\in \Gamma}dg_1\cdot \cdot \cdot dg_l \,
K_{g_1,\cdot \cdot \cdot , g_l}(z,w)
$$
 where for any $(z,w) \in  \widehat{ M} \times  \widehat{ M} $,  $K_{g_1,\cdot \cdot \cdot , g_l}(z,w)
\in
\bigwedge^k(T^\ast_{\pi(z)}T) \otimes Hom (\widehat{E}_w\,,\,  \widehat{E}_z)$ and
vanishes for $\pi(z)(g_1\cdot \cdot \cdot  g_l)^{-1} \not= \pi(w)$. Then we state the following
\begin{definition} \label{inte} For any $\theta \in T$ we set
$${}^b STR_{<e>}(K)(\theta):=$$ $$ \sum_{g_0,\cdot \cdot \cdot ,
g_l \in \Gamma,\, g_0\cdot \cdot \cdot g_l =e} (dg_1\cdot \cdot
\cdot dg_l ) g_0 \, \nuint_{ \Delta_b(\theta)} \phi(w)  Str
K_{g_1,\cdot \cdot \cdot , g_l}(wg_0^{-1},w) \,d{\rm Vol}^b_{
Z_\theta} (w) $$ Then ${}^b STR_{<e>}(K)$ defines an element of
$\OM$.
\end{definition}

Of course, this definition extends to any element of
$ \Psi^{-\infty, \delta}_{b,\widehat{\Omega}_{\ast}(T,\Bi_\Gamma)}
(\widehat{M}\,;\, \widehat{E})$.

\s
\n
{\bf Remark.} The previous definition is also valid for any
$K \in \Psi^{-\infty, \delta}_{\widehat{\Omega}_{\ast}(T,\Bi_\Gamma)}
(\widehat{M}\,;\, \widehat{E})$ but in this case the regularized integrals
in the above definition  are in fact  ordinary integrals because the indicial family of $K$ vanishes
identically.

Then, as in \cite{LPMEMOIRS} (see Proposition 13.5) one gets the following commutator
formula
\begin{proposition} \label{commutator} For any $K,K^\prime $ belonging to
$ \Psi^{-\infty, \delta}_{b,\widehat{\Omega}_{\ast}(T,\Bi_\Gamma)}
(\widehat{M}\,;\, \widehat{E})$ one has:
$${}^b STR_{<e>}[K,K^\prime]\,=\,
{\sqrt{-1} \over 2 \pi}
\int_\RR STR \bigl(\, {\partial \over \partial \lambda}I(K,\lambda)\circ I(K^\prime,\lambda)\,\bigr) d\lambda.
$$ If we replace $K$ by a differential operator $\in {\rm Diff}^1_{b, \Ti}(
\wm, \we)$ and $K^\prime$ by the composition of $K^\prime$ with
an element of the calculus with bounds $\Psi^{\ast,\delta}_{b, \Ti}(\wm, \we)$
then the same commutator formula is valid.

\end{proposition}

\bigskip Now we fix a $\Gamma-$invariant horizontal distribution $T^H\widehat{M}$ such that
$$ {{}}^bT\widehat{M} = T^H\widehat{M} \oplus
{{}}^bT(\widehat{M}/T) $$ and, as in \cite{Go-Lo} and \cite{MP I}
(section 9) we consider the Bismut superconnection $$
\AA_s^{Bismut}=s\D+\nabla^{u}-{1\over 4 s}c(\tau),\; s\in
\RR^{+\ast} $$ where $\D$ is the Dirac operator introduced in
Lemma \ref{D}, $c(\tau) $ denotes the Clifford multiplication by
the curvature 2-form $\tau$ of $T^H\widehat{M} $ and $\nabla^{u} $
is a certain unitary connection. Then, as in \cite{Go-Lo}  we fix
a function $h \in C^\infty_c(\widehat{M})$ constant in the normal
direction near the boundary such that $\sum_{\gamma \in \Gamma}
\gamma\cdot h =1$  on $\widehat{M} $ and
 consider for each real $s>0$ the superconnection
$$\AA_s= \AA_s^{Bismut} +\sum_{\gamma \in \Gamma} d\gamma \otimes h \gamma^{-1}
$$ which sends $C^\infty_{{\cal T}^\infty}(\widehat{M},\we)$ into
$\widehat{\Omega}_{\ast}(T,\Bi_\Gamma) \otimes_{{\cal T}^\infty}
 C^\infty_{{\cal T}^\infty}(\widehat{M},\we)$.
By developing in a straightforward way a $\Ti$-$b-$heat calculus as in
\cite{LPMEMOIRS} (section 10) one sees easily that
$e^{-s^2D^2} \in \Psi^{-\infty}_{b,{\cal T}^\infty}
(\widehat{M}\,;\, \widehat{E})$ for any $s>0$.
 Moreover, using a Duhamel expansion around
$e^{-s^2D^2}$ one checks that for any $s>0$
$$
e^{-\AA_s^2} \in \Psi^{-\infty}_{b,\widehat{\Omega}_{\ast}(T,\Bi_\Gamma)}
(\widehat{M}\,;\, \widehat{E})
$$

\bigskip Now as in \cite{MP I} (section 10) we consider for
any $s>0$  the induced boundary superconnection $$ \BB_s=s \sigma \D_0
+ \nabla^{u}_\partial - {\sigma  c(\partial \tau) \over 4 s}
+\sum_{\gamma \in \Gamma} d\gamma \otimes h \gamma^{-1} $$ with
$D_0$  the boundary Dirac operator of $D$. Then
according to Proposition \ref{higheta}
 the higher eta invariant
$$
\widetilde{\eta}_{<e>} =
{2 \over \sqrt{\pi}}\R\, \int_0^{+\infty}
STR_{<e>, Cl(1)} ( \, {d \BB_s \over ds } e^{-\BB^2_s}\,)\, ds \in \OM
$$ is well defined.
\medskip Next we recall from Section 2 of \cite{Go-Lo} the definition
of the  following
connection
$\nabla^{can}$.
 $$
 \forall f\in C^\infty_c(\wm),\;
 \nabla^{can} f = d^{\wm}\,f \oplus \sum_{\gamma \in \Gamma}
 d\gamma \otimes h(\gamma^{-1} \cdot f)
 .$$ Then $(\nabla^{can})^2$ acts on $C^\infty_c(\wm)$ as left multiplication
 by a $2-$form  $\Theta$ which commutes with $C^\infty_c(\wm)\rtimes \Gamma$.
 Explicitly,
 $$
 \Theta= \sum_{\gamma \in \Gamma} d^{\wm}(\gamma\cdot h) d\gamma \gamma^{-1}
 -
 \sum_{\gamma, \gamma^\prime \in \Gamma}
 (\gamma \gamma^\prime \cdot h ) (\gamma  \cdot h ) d \gamma
 d \gamma^\prime (\gamma \gamma^\prime)^{-1}.
 $$ Then put ${\rm ch}(\nabla^{can})= e^{-{\Theta \over 2 \pi i}}$,
 this Chern character  lies in
 \begin{equation} \label{space}
 \Pi_{k=0}^\infty \oplus_{l=0}^{+\infty} \Omega^{k-l}(\wm) \hat{\otimes}
 \Omega^{k+l}(\CC \Gamma).
 \end{equation}
Observe that for any function $\psi \in C^\infty_0(\wm,\CC)$, we
may write $$ \psi \,e^{-{\Theta \over 2 \pi i}} = \sum_{k\in \NN}
\omega_k $$ where each $\omega_k $ belongs to the algebraic tensor
product $ \Omega^\ast(\wm) \hat{\otimes} \Omega^\ast(\CC \Gamma)$ and is
of total degree $k$.

\subsection{The short time limit.}\label{shortsub}$\;$

\m

 Now we may state the higher local index theorem whose
proof proceeds along the line of the proofs of Theorem 2 of
\cite{Go-Lo}, (
see also  Theorem 13.6 of \cite{LPMEMOIRS} for the $b-$part ).
\begin{theorem} \label{localindex}
$$
\lim_{s\rightarrow 0^+} \R ({}^b STR_{<e>} \, e^{-\AA_s^2}) \,=
\int_Z \phi \,\widehat{A}\,(\nabla^{TZ})\,{\rm ch}(\nabla^{\widehat{V}}) \,
{\rm ch}(\nabla^{can}) \, \in
\OM
$$ where $Z$ denotes the typical
fiber of $\pi:\wm \rightarrow T$ and where we recall that $\phi \in C^\infty_c(\wm)$ denotes a function which is constant in the normal direction near
the boundary and such that $\sum_{\gamma \in \Gamma} \gamma \cdot \phi = 1$
\end{theorem}

Recall now that $M$ is endowed with a foliation
$\cal{F}$ whose leaves are the images by the covering map $p: \wm
\rightarrow M$ of the fibers of $\pi: \wm \rightarrow T$. Then,
let  $\Phi$ be a closed graded trace of degree $n$ on the graded
algebraic tensor product $\sum_{k+l=n} \Omega^k(T) {\otimes}
\Omega^l(\CC\Gamma)$.
 By proposition 2 of \cite{Go-Lo}
we may write $\Phi=\sum_{k+l=n} \tau_{k,l}$ where $\tau_{k,l}\in
C_k(T) \otimes C^l(\Gamma)$ ($ C_k(T)$ denoting the set of
currents of degree $k$); we shall assume that $\Phi$ extends as a
closed graded trace on $\widehat{\Omega}_\star ( T, \Bi_\Gamma)$.
  We shall now show that associated to $\nabla^{can}$
 and $\Phi$ there is a well defined current $\omega_\Phi$ on $M$.
  First, by pairing with respect
  to the $\Gamma-$variables the Chern character associated to
  $\nabla^{can}$ and $\Phi$,
   we obtain
  an element  $$< {\rm ch}\, \nabla^{can}\,;\, \Phi >_\Gamma\; \in\;\Omega^\ast(\wm) \otimes C_\ast(T)\,;$$
  then since
  $\pi:\wm \rightarrow T$ is a fibration we can associate to any
  $\omega \otimes \beta \in \Omega^\ast(\wm) \otimes C_\ast(T)$ an
  element $\omega \otimes \pi^\ast(\beta)$ of $\Omega^\ast(\wm) \otimes C_\ast(\wm)$.
Observe now that $ C_\ast(\wm)$ is a module over
$\Omega^\ast(\wm)$. Summarizing we can define $$
  A: \Omega^\ast(\wm) \otimes C_\ast(T) \rightarrow C_\ast(\wm).
  $$
by $A(\omega \otimes \beta) (\lambda)=
\pi^\ast(\beta)(\lambda\wedge\omega)$. Thus $A( <{\rm ch}\,
\nabla^{can}\,;\, \Phi >_\Gamma)$ is well defined as a current on
  $\wm$ and it is,  moreover, $\Gamma-$invariant.
We  denote by $\omega_\Phi$ the
  induced current on $M$.
Recall that $p$ denotes the covering map $p: \wm \rightarrow M $.
Thus, {\it  by the very definition of
$\omega_\Phi$},  one has for each $\alpha \in \Omega^*(M)$:
$$
< \alpha\,;\, \omega_\Phi > = < \phi \, p^*(\alpha)\,;\, A( <{\rm ch}\,
\nabla^{can}\,;\, \Phi >_\Gamma) >;
$$  we  refer to
  the proof of formula
(\ref{form}) below  for the independence
 of
the right hand side of the above formula
on the choice of $\phi$ satisfying
$\sum_{\gamma \in \Gamma} \gamma \cdot \phi =
1$.


Summarizing, we may state the following corollary:
  \begin{corollary} \label{corocurrent} Let $\Phi$ be a closed graded trace on
  $\widehat{\Omega}_*(T, \Bi_\Gamma)$ concentrated in the trivial conjugacy class. Then
  there is a current $\omega_\Phi$ on $M$ such that the following formula holds:
  $$
< \int_Z \phi \,\widehat{A}\,(\nabla^{TZ})\,{\rm ch}(\nabla^{\widehat{V}}) \,
{\rm ch}(\nabla^{can})\, ; \, \Phi >=
  <\widehat{A}(T\cal{F})\,{\rm ch}(\nabla^{V})\,,\, \omega_\Phi>
    $$
  where on the right hand side the pairing is the one between currents and
  differential forms on $M$.
  \end{corollary}
\begin{proof}
We write: $$ <\int_Z \phi \,\widehat{A}\,(\nabla^{TZ})\,{\rm
ch}(\nabla^{\widehat{V}}) \,{\rm ch}(\nabla^{can})\,;\, \Phi >=
< \phi \,\widehat{A}\,(\nabla^{TZ})\,{\rm
ch}(\nabla^{\widehat{V}})
 \,;\,A(< {\rm ch}\nabla^{can}\,;\,\Phi >_\Gamma)\,>
$$ where on the right hand side the pairing is the one between currents and
  differential forms on $\wm$.
Next we observe that $\widehat{A}\,(\nabla^{TZ})\,{\rm
ch}(\nabla^{\widehat{V}})$ is a $\Gamma-$invariant differential
form on $\wm$, it is then the pull back by the covering map $\wm
\rightarrow M$ of $\widehat{A}(T\cal{F})\,{\rm ch}(\nabla^{V}) $.
 We then get immediately by definition of $\omega_\Phi $:
 $$
 <\widehat{A}(T\cal{F})\,{\rm ch}(\nabla^{V})\,,\, \omega_\Phi>\,=\,
< \phi \,\widehat{A}\,(\nabla^{TZ})\,{\rm
ch}(\nabla^{\widehat{V}})
 \,;\,A(< {\rm ch}\nabla^{can}\,;\,\Phi >_\Gamma)\,>.
$$
The Corollary is thus proved.
\end{proof}

\s
\n {\bf Example 1 .} Let us describe $\omega_\Phi$ when $T=\{ pt \}$.
Then we are dealing with a single Galois covering $\Gamma\to
\widetilde{M}\to M$, with $M$ orientable. If $\Phi$ is a
graded closed trace on $\Omega_* (\Bi_\Gamma)$ and $\lambda$ is a
smooth form on $M$ then
$$ \forall \lambda \in \Omega^*(M),\;
\omega_\Phi(\lambda)=\int_M \lambda\wedge <\omega,\Phi>$$
with $<\omega ,\Phi>$
equal to a smooth differential form on $M$, see \cite{Lott II}.
 Therefore, in this case, $\omega_\Phi $ is simply
the integration against $ < \omega ,\Phi> \in \Omega^*(M)$.

\s
\n
{\bf Example 2 .} We now describe  $\omega_\Phi$ in another   simple
example. Assume that the $\Gamma-$equivariant fibration $\pi$ is
trivial, $\pi: \wm=\widetilde{M} \times T \mapsto T $ and that the
quotient $ \widetilde{M}/ \Gamma$ is a
 smooth compact manifold with
boundary. Consider a $\Gamma-$invariant measure $\mu$ on $T$
(since $\Gamma$ is of polynomial growth such a measure always
exists). We define $\Phi$ by:
\begin{equation} \label{exemple}
\forall \,\Xi\,=\,\sum_{\gamma \in \Gamma} \alpha_\gamma (\theta)
\, \gamma\, \in \Ti\,=\,\widehat{\Omega}_0(T, \Bi),\quad \Phi (\Xi) =
\mu(\alpha_e).
\end{equation}
Thus
\begin{equation}
\Phi=\mu\otimes \delta_e
\end{equation}
with $\delta_e$
the trivial
cyclic 0-cocycle on $\Bi_\Gamma$.
Since ${\rm ch}\, \nabla^{can}$ belongs to the space
(\ref{space}) and the component of ${\rm ch}\, \nabla^{can}$
living in $\Omega^0(\wm) \widehat{\otimes} \Omega^0(\CC \Gamma)$
is identically equal to $1$, one checks easily that $A( <{\rm
ch}\, \nabla^{can}\,;\, \Phi >_\Gamma)\,=\,  \mu\circ \pi_* $
where $ \pi_*$ denotes the integration in the fibers. Thus
$\omega_\Phi (\alpha)= \mu( \pi_* (\phi p^* (\alpha)))$. Recall
now  that the fibers of $\pi:\wm = \widetilde{M} \times T
\rightarrow T$ induce a foliation $\cal{F}$ on $M=\frac
{\widetilde{M} \times T} {\Gamma}$ and that $\mu$ induces an
invariant transverse measure, still denotes $\mu$, on
$(M,\cal{F})$. Let $\phi$ be any function belonging to $
C^\infty_c(\wm, [0,1])$ such that $\sum_{\gamma \in \Gamma} \gamma
\cdot \phi = 1$ on $\wm$. Then the Ruelle-Sullivan current $RS$ on
$(M,\cal{F})$ associated to the transverse measure $\mu$ is given
by: $$ \forall \alpha \in C^\infty(M, {\wedge}^* T^* M),\quad RS
(\alpha) = \mu ( \int_{\widetilde{M}} p^*(\alpha) \phi ) $$ where
$p:\wm \rightarrow M$ denotes the covering map. Summarizing, in
this particular case, \begin{equation}\label{omegaphi=rs}
\omega_\Phi=RS\,.
\end{equation}

\s
\n
{\bf Remark.}
 Notice that if
$\psi$ is any other function belonging to $ C^\infty_c(\wm,
[0,1])$ such that $\sum_{\gamma \in \Gamma} \gamma \cdot \psi = 1$
on $\wm$ then $  \mu(\int_{\widetilde{M}} p^*(\alpha) \phi) =\mu(
\int_{\widetilde{M}} p^*(\alpha) \psi).$ Indeed, using the fact
that $p^*(\alpha)$ is $\Gamma-$invariant, that $\Gamma$ acts by
isometries on $\widetilde{M}$ (by fixing a $\Gamma$-invariant
metric) and that $\mu$ is $\Gamma$-invariant,
 one simply
writes \begin{equation} \label{form} \mu(\int_{\widetilde{M}}
p^*(\alpha) \phi ) =\mu (\int_{\widetilde{M}} p^*(\alpha) \,\phi
  (\sum_{\gamma \in \Gamma} \gamma \cdot \psi) )=
  \mu (\sum_{\gamma \in \Gamma}\int_{\widetilde{M}}
  p^*(\alpha)\, \phi \,\gamma \cdot \psi)=
  \end{equation}
  $$
\mu(
  \sum_{\gamma \in \Gamma}\int_{\widetilde{M}}
  p^*(\alpha) \,\psi \, \gamma^{-1} \cdot \phi)   =
\mu(\int_{\widetilde{M}} ( p^*(\alpha) \, \psi
  \sum_{\gamma \in \Gamma}\gamma^{-1} \cdot \phi \,))\,=\,\mu(
  \int_{\widetilde{M}} p^*(\alpha) \psi).
   $$

\n
{\bf End of remark.}

\s
Now, we are going to interpret $< \widetilde{\eta}_{<e>} ; \Phi>$
as an $L^2-$type invariant. To this
  aim we consider first a closed odd-dimensional
Galois covering $\widetilde{N}\to N$
and the equivariant fibration $\pi: \widetilde{N}\times T\rightarrow T$.
We shall later apply our remarks to $\widetilde{N}\times T =\partial\widetilde{M}\times T$.
We still consider on  $T$  a $\Gamma$-invariant measure $\mu$.
First we construct a $L^2$-foliated trace  on
the space of  $\Gamma-$equivariant
  families of smoothing
  operators acting in the fibers of $\pi:\wn\rightarrow T$, with $\wn=\widetilde{N}\times T$.

  \begin{definition} \label{def:bok} Let $K=(K_\theta)_{\theta \in T}$ be a $\Gamma-$equivariant
  family of smoothing
  operators acting in the fibers of $\pi: \wn=\widetilde{N} \times
 T \rightarrow T$ on the sections of a equivariant
vector bundle $\wf$.
Then
  $tr\,K_\theta (z,z) \,d Vol_{\pi^{-1}(\theta)}$ defines a longitudinal
  measure, denoted $tr\, K_\theta (y,y) \,d Vol_{\pi^{-1}(\theta)}(y)$,  on the leaves of the foliation $\cal{F}$ of $N$
  induced by the image under the covering map $p:  \wn \rightarrow N$
  of the fibers $\widetilde{N}\times \{\theta\}$, $\theta \in T$.
  Denote still by $ \mu$ the invariant transverse measure induced by $\mu$ on
$(N, {\cal{F}} ) $.
The $L^2$-$\mu$-trace of $K=(K_\theta)_{\theta \in T}$, denoted $\Tr_{(2),\mu} K$,
is by definition
the integration of the longitudinal measure
$  tr\, K_\theta (y,y) \,d Vol_{\pi^{-1}(\theta)}(y)$
against the transverse measure $\mu$.
\end{definition}

\begin{lemma}\label{bok}
For the $L^2$-$\mu$-trace of $K=(K_\theta)_{\theta \in T}$ the following
formula holds:
  $$\Tr_{(2),\mu} K
=  \int_T ( \int \phi \,
  tr \,K_\theta(z,z) \,d Vol_{\pi^{-1}(\theta)}) d\mu(\theta).
  $$
\end{lemma}

\begin{proof} Let $F$ be a fundamental domain of $\widetilde{N}$
  for the $\Gamma-$action. We fix a point $m_0\in \widetilde{N}$ and
  consider the transversal $\cal{T}_0$ of $\cal{F}$ in $N$ induced
  by the immersion of $\{m_0\} \times T$  in $\wn$. For any
  $(z, \theta)  \in \wn \setminus
  ( \,(\cup_{\gamma \in \Gamma}
  \partial F\cdot \gamma ) \times T $ we consider the unique $\gamma \in \Gamma
  $ such that $z \cdot \gamma \in F \setminus \partial F$ and
  set
  $\chi( (z, \theta)\cdot \Gamma) = ( m_0 ,\theta \cdot \gamma) \in \cal{T}_0$. We define
the restriction of $\chi$ to
$$
\{(z, \theta)\cdot \Gamma, \; (z, \theta) \in  ( \,(\cup_{\gamma \in \Gamma}
  \partial F\cdot \gamma ) \times T \}
$$ to be any borel map such that
  $\chi( (z, \theta)\cdot \Gamma) = ( m_0 ,\widetilde{\theta})$ where there
exists $ \widetilde{z} \in \partial F$ such that $ (z,
\theta)\cdot \Gamma=(\widetilde{z},  \widetilde{ \theta})\cdot
\Gamma$. We have thus
  defined a borel map $\chi: N = \frac {\widetilde{N} \times
  T} {\Gamma} \rightarrow \cal{T}_0$ such
  that for any $x \in N $, $\chi(x)$ belongs to the leaf of $x$. Notice that the
following family
$$(\,tr \,K_\theta(z,z)
  \,d Vol_{\pi^{-1}( \theta)}(z)\,)_{\theta \in T}
$$ of $\Gamma-$invariant measures along each fiber of $\pi $ descends on $N$ and
thus  induces a longitudinal measure of $(N, \mathcal{F})$ that we  denote by
$$tr \,K_\theta(y,y)
  \,d Vol_{\pa\pi^{-1}( \theta)}(y).
$$
Then, by definition, (see \cite{MoS}):
  $$\Tr_{(2),\mu} K =
\int_T  \int 1_{\chi^{-1}(\{(m_0 ,  \theta)\})}(y) tr \,K_\theta(y,y)
  \,d Vol_{\pi^{-1}( \theta)}(y) d\mu(\theta).
  $$
   Using the fact that $\cup_{\gamma \in \Gamma}
  (\,\partial F\cdot \gamma \,)$ is of measure zero in $\widetilde{N}$,
  we observe that this last integral is equal to
  $$
  \int_T \int  1_F(z) tr \,K_\theta(z,z)\,
  d Vol_{\pi^{-1}( \theta)}(z) d\mu(\theta).
  $$ Since $\sum_{\gamma \in \Gamma} 1_F \cdot \gamma =1$ on
  $\widetilde{N}$, we can use the same reasoning
  as in the proof of formula (\ref{form}) in order to replace in the previous integral
  $1_F$ by $\phi$. This proves the Lemma.
 \end{proof}

We {\it define} the $L^2$-foliated eta invariant of the family $D_0$
as:
$$
\eta_{(2),\mu}\,:=\, \frac{1}{\sqrt{\pi}} \int_0 ^\infty \eta_{(2),\mu}(t) dt
\;\;\text{
with}\;\;\eta_{(2),\mu}(t) =\Tr _{(2),\mu} (D_0 \exp (-t^2 D_0 ^2)).$$
Now, recall that
  $\widetilde{\eta}_{<e>}$ and $\Phi$ are respectively
  defined in Proposition
 \ref{higheta} and formula (\ref{exemple}). Then, by applying Duhamel's formula
 one sees immediately that
 $$
 \sqrt{\pi}< \widetilde{\eta}_{<e>} ; \Phi > = \int_0^{+\infty}
\int_T (\, \int \phi \,
  tr \,K( D_0 e^{-t^2 D_0^2})(z,z) \,d Vol_{\pi^{-1}(\theta)} \,dt\,) d\mu(\theta)
 $$ where of course $K (D_0 e^{-t^2 D_0^2}) ( \cdot , \cdot )$ denotes
 the (fiberwise) Schwartz kernel of $D_0 e^{-t^2 D_0^2}$. But then
 Lemma \ref{bok} shows that
$$ < \widetilde{\eta}_{<e>} ; \Phi > = \eta_{(2),\mu}
$$
 Summarizing, we have proved the following:\\
{\it if $\Phi$ is given
by \eqref{exemple} then}
$$
< \widehat{A}(T\mathcal{F}) ; \omega_\Phi > - \frac {1} {2}
<  \widetilde{\eta}_{<e>} ; \Phi >\; =\;
< \widehat{A}(T\mathcal{F}) ; RS >
 - \frac {1} {2} \eta_{(2),\mu}
$$

\n {\bf Remark.} Recall that the measured eta-invariant $\eta_\mu$
of \cite{Ra} is given by:
 $$
 \sqrt{\pi}\eta_\mu=
\int(\, \int_0^{+\infty} tr\, K( D_\pa e^{-t^2 D_\pa^2})(y,y) \,d Vol_{\pi^{-1}(\theta)}(y) \,dt
\,) d\mu(\theta)
$$
where $D_\pa$ is the induced boundary longitudinal Dirac operator along
the leaves of $(\partial M, {\cal{F}})$, $D_0$ being of course
the lift of $D_\pa$. In general the two longitudinal measures
$$
tr\, K( D_\pa e^{-t^2 D_\pa^2})(y,y) \,d Vol_{\pi^{-1}(\theta)}(y) \not=
tr\, K( D_0 e^{-t^2 D_0^2})(y,y) \,d Vol_{\pi^{-1}(\theta)}(y)
$$ are different.  So in general $\eta_\mu \not= \eta_{(2),\mu}$. This is of course well known when
$T$=point (indeed, on a Galois covering the $L^2$-eta invariant is in general
not equal to
the eta invariant on the base).

\n {\bf End of remark.}

\n
{\bf End of example 2.}

\s
\n
{\bf Example 3.} Let $C$ be a $\Gamma$-invariant closed current on
$T$. Consider $\Phi=C\otimes \delta_e$.
Then, proceeding as in the previous example, we obtain
$$\omega_\Phi = \mathcal{C}$$
with $\mathcal{C}$ the foliated closed current on $M$ induced by $C$.

\medskip

Let us go back to the case of a general $\Gamma-$invariant
fibration $\pi: \wm \mapsto T$. In order to prove our main result,
we shall  need the following transgression formula:

\begin{proposition} Let $u>t>0$ and let $\widetilde{\eta}_{<e>}(s)$ be the
higher eta integrand introduced in Section \ref{heta}. Then
\begin{align}\label{tras}
\R ({}^b STR_{<e>} \, e^{-\AA_u^2})\,& =\, \R ({}^b STR_{<e>} \,
e^{-\AA_t^2}) -d\int_t^u\R {}^b STR_{<e>}(\,{d\AA_s \over ds}
e^{-\AA_s^2}\,)ds +\notag \\ & -{1\over 2} \int_t^u
\widetilde{\eta}_{<e>}(s) ds.
\end{align}
\end{proposition}
\begin{proof} One proceeds exactly as in the proof of
Proposition 11 of \cite{MP I}, using Proposition \ref{commutator}.
\end{proof}

Taking the limit as $t\rightarrow 0^+$ in the previous formula and
using Theorem \ref{localindex} one gets the following crucial
equality in $\OM$:

\begin{align} \label{transgression}
\R ({}^b STR_{<e>} \, e^{-\AA_u^2})\, &=\, \int_Z \phi
\,\widehat{A}\,(\nabla^{TZ})\, {\rm ch}(\nabla^{\widehat{V}})
\,{\rm ch}(\nabla^{can}) - {1\over 2} \int_0^u
\widetilde{\eta}_{<e>}(s) ds+ \notag \\ & - d\int_0^u\R ({}^b
STR_{<e>}(\,{d\AA_s \over ds} e^{-\AA_s^2}\,))ds.
\end{align}

\medskip

\section{{\bf A higher APS-index theorem: the isometric case}}\label{isomcase}

\subsection{Isometric actions.}\label{isomcase1}$\;$

\medskip
In this subsection we shall assume that $T$ is endowed with a
riemannian metric such that $\Gamma$ acts by isometries. This
corresponds to a type-II situation; we shall deal with
arbitrary $\Gamma$-actions in the next subsection.
The following Proposition gives an important consequence of
this assumption

\begin{proposition}\label{decrease-II}
Let us assume that the manifold $T$ admits a riemannian metric $g$
such that the group $\Gamma$ acts by isometries on $T$. Then $\Ti$
is  closed under the holomorphic functional calculus in $C^0(T)
\rtimes_r \Gamma$.
\end{proposition}

\n
\begin{proof}
$\Ti$ is equal to the set of the elements $\sum t_\gamma \gamma$
such that
\begin{equation}\label{decrease-f}
{\rm sup}_{\theta\in T,\gamma\in \Gamma}\left[ |\Delta^p t_\gamma
|(\theta) (1+||\gamma||)^N \right]<\infty \;, \;\;\forall p,
N\in\NN
\end{equation}
where $\Delta$ is the Laplace Beltrami operator of $T$.
Using  the hypothesis that $\Gamma$ acts by isometries and
standard arguments (see \cite{Lafforgue}) one checks that this set
is indeed a subalgebra and that  is closed under the holomorphic
functional calculus in $C^0(T) \rtimes_r \Gamma$. Alternatively,
one may check that $\Ti$ is
 a particular case
of a construction due to  Jiang (see section 3 of \cite{Jiang} ).
\end{proof}
In particular
\begin{equation}\label{isomor}
K_0 (\Ti)\simeq K_0 (C^0 (T)\rtimes \Gamma) .
\end{equation}

\smallskip
Using the fact that $\Ti$ is dense and stable under
holomorphic functional calculus and  proceeding as in the
proof of the decomposition theorem  of \cite{LPMEMOIRS} (see  Theorem 12.7
there) we
can   give a very explicit  decription of  the
index class in $K_0(\Ti)$.
Briefly, we can find $\epsilon > 0,$
 and $ \cal{L}_\infty$ [resp. $\cal{N}_\infty$] a sub-$\Ti-$module projective of
finite rank of $x^\epsilon H^\infty_{b, \Ti}(\wm, \we^+ )$
[resp. $x^\epsilon H^\infty_{b, \Ti}(\wm, \we^- )$]
such that
\begin{equation}\label{newindex}
\Ind (D^+)=[\cal{L}_\infty]-[\cal{N}_\infty]\;\;\text{in}\;\;K_0 (\Ti)
\end{equation}
Further properties of these sub-modules can be found in Appendix C.

\subsection{The APS index theorem for the groupoid
$T\rtimes \Gamma$ in the isometric case.}\label{isomcase2} $\;$

\medskip
Recall that ${\rm Ch}( \Ind D^+ )$ is defined (thanks to the work
of \cite{Karoubi}) by a representative in
$\widehat{\Omega}_*(\Ti)$ whereas the $b-$supertrace of the
superconnection heat kernel is defined by a representative in
$\widehat{\Omega}_*(T, \Bi_\Gamma)$. In order to establish a
connection between these two objects we state the following
proposition whose easy proof will be left to the reader.
\begin{proposition} \label{compa}

 \item{1).} By universality there is a natural morphism of algebras:
$$ j: \widehat{\Omega}_*(\Ti) \rightarrow \widehat{\Omega}_*(T,
\Bi_\Gamma). $$

\item{2).} Let $\cal{F}= \cal{L} \oplus \cal{N}$ be a
$\ZZ_2-$graded projective finitely generated $\Ti-$module. Let
$\nabla$ be a $\Ti-$connection of $\cal{F}$ with
value in $\widehat{\Omega}_1(\Ti)\otimes \cal{F}$ and
preserving both $\cal{L}$ and $\cal{N}$. Then $ \widetilde{\nabla}=
(j\otimes \Id_{\cal{F}}) \circ \nabla $ defines a connection of
$\cal{F}$ with values in $\widehat{\Omega}_1(T , \Bi_\Gamma)\otimes \cal{F}$.

\item{3).} The two supertraces (in the algebraic sense) $\R STR_{<e>}^{alg} e^{-\widetilde{\nabla}^2}$ and
$j_{<e>}(\R STR^{alg} e^{-{\nabla}^2})$ coincide  in
$\overline{\widehat{\Omega}}_*(T , \Bi_\Gamma)$ (ie modulo graded
commutators). Moreover for any $\Ti-$connection $\nabla_1$ of
$\cal{F}$ with values in $\widehat{\Omega}_1(\Ti)\otimes \cal{F}$
and preserving both $\cal{L}$ and $\cal{N}$, the two supertraces
$\R STR_{<e>}^{alg} e^{-\widetilde{\nabla}^2}$ and $\R
STR_{<e>}^{alg} e^{-{\nabla_1}^2}$ coincide modulo an exact form.
\end{proposition}

Now we may state the higher index theorem for the groupoid
$T\rtimes \Gamma$ and the proper$T\rtimes \Gamma-$manifold  $\wm$.

\medskip
\begin{theorem} \label{higherindex}
Let $\Gamma$ be of polynomial growth and let $D^+= (D^+
(\theta))_{\theta\in T}$ a $\Gamma$-equivariant family of Dirac
operators with inverible boundary family $D_0$. Then the following
formula holds in the homology $\widehat{H}_\ast(T, \Bi_\Gamma)$
(see Section \ref{heta1}): $$ j_{<e>}({\rm Ch}\, \Ind D^+)\,=\, \int_Z
\phi \,\widehat{A}\,(\nabla^{TZ})\,{\rm ch}(\nabla^{\widehat{V}})
\,{\rm ch}(\nabla^{can}) - {1\over 2} \widetilde{\eta}_{<e>} \in
\widehat{H}_\ast(T, \Bi_\Gamma) $$
 $j_{<e>}({\rm Ch}\, \Ind D^+) $  denotes the projection
of $j({\rm Ch}\, \Ind D^+) $ on the set of differential forms
$\sum_{\gamma_0 \gamma_1 \cdots\gamma_k =e}\omega_{\gamma_0, \cdots, \gamma_k} \otimes
 \gamma_0 d\gamma_1 \cdots d\gamma_k$
of $\widehat{\Omega}_\star ( T, \Bi_\Gamma)$ which are concentrated (with respect to the group
variables) in the trivial conjugacy class $<e>$ of $\Gamma$.
\end{theorem}

\begin{definition}\label{higherindexdefinition}
Let $\Phi$ be a closed graded trace on $\widehat{\Omega}_\star (
T, \Bi_\Gamma)$ as above. We define the higher $\Phi$-index of
$D^+$ as $$\Ind_\Phi (D^+)\,:= \, <j_{<e>}({\rm Ch}\, \Ind
D^+),\Phi>\,.$$
\end{definition}
\begin{corollary}\label{corolomega}
For the higher $\Phi$-index of a $\Gamma$-equivariant family as
above the following formula holds:
 $$  \Ind_\Phi (D^+) = <\widehat{A}(T\cal{F})\,{\rm ch}(\nabla^{V})
 \,,\,
\omega_\Phi >- <{1\over 2}  \widetilde{\eta}_{<e>}\,;\, \Phi > $$
\end{corollary}

The corollary follows from the Theorem, using the discussion in
Subsection \ref{shortsub}. The proof of the theorem can be found in Appendix C.

\section{{\bf A higher APS index theorem: the general  case}}\label{mainthsub}

\subsection{The Chern character of the index class}\label{chernsub}$\;$

\medskip
In the general case, when $\Gamma$ does not preserve any
Riemannian metric on $T$, it is not known whether the algebra
$\Ti$ is  stable (or not stable) under the holomorphic functional
calculus in $C^0(T) \rtimes_r \Gamma$. For this reason, the
definition of the Chern character of the index class and of  its
pairing with a closed graded trace on $\widehat{\Omega}_*
(T,\Bi_\Gamma)$ is more involved.
 We shall devote the beginning of this
subsection to a few preliminaries (following \cite{Co} and
\cite{Go-Lo}) leading to such a definition. We follow
the notations in \cite{Go-Lo}.

We consider the following algebras:
 $$ \cal
U=\Ti=C^\infty(T,\Bi_\Gamma)\,,\quad \Omega^* =\widehat{\Omega}_*
(T,\Bi_\Gamma)$$
$$
\widetilde{\cal U}=\Psi^{-\infty,\delta}_{\Ti}(\wm , \we^+) \,,\quad
\widetilde{\Omega}^*= \Psi^{-\infty,\delta}_{\widehat{\Omega}_*
(T,\Bi_\Gamma)}( \wm , \we^+).
$$ Using a
$\Gamma-$equivariant complex vector bundle isomorphism: $\we^+ \simeq \we^-$
we get natural identifications:
$$
\Psi^{-\infty,\delta}_{\Ti}(\wm , \we^+) \simeq
\Psi^{-\infty,\delta}_{\Ti}(\wm , \we^-),\quad
 \Psi^{-\infty,\delta}_{\widehat{\Omega}_*
(T,\Bi_\Gamma)}( \wm , \we^+) \simeq
 \Psi^{-\infty,\delta}_{\widehat{\Omega}_*
(T,\Bi_\Gamma)}( \wm , \we^-)
$$ that we shall use in this subection.
Then we consider the $\Ti-$module $\cal{E}= H^\infty_{b , \Ti}(\wm , \we)$ and
(with the notations defining
$\AA_s$) the connection $\nabla= \nabla^u + \sum_{\gamma \in \Gamma} d\gamma \otimes
 h \gamma^{-1}$.Observe that
$\nabla$ sends $\cal{E}$ into $\Omega_1(T, \Bi)
\otimes_{\Ti}\cal{E}$ and that one defines a degree one map  $
\widetilde{\nabla}: \widetilde{\Omega}^{*} \mapsto
\widetilde{\Omega}^{* + 1}$ by setting, for any $K \in
\widetilde{\Omega}^*$, $ \widetilde{\nabla} (K) = [ \nabla, K] \in
\widetilde{\Omega}^{* + 1}$, where $[ \nabla, K]$ denotes the
graded commutator. Now set $\Theta = \nabla^2 $, $\Theta$ sends
$\cal{E} $ into $\Omega_2(T, \Bi) \otimes_{\Ti} \cal{E} $ but is
not given by a smooth integral kernel in $\widetilde{\Omega}^*$,
nevertheless for any $K \in \widetilde{\Omega}^*$, $\Theta K  $
and $K
 \Theta$ both belong to $\widetilde{\Omega}^*$
One then checks that for any
$K \in \widetilde{\Omega}^*$, $\widetilde{\nabla}^2 (K)= \Theta K - K \Theta$.

Let $\Phi$ a closed graded trace on $\widehat{\Omega}_*
(T,\Bi_\Gamma)$, concentrated on the trivial conjugacy class. Then
$\Phi$ extends to a closed graded trace $\widetilde{\Phi}$ on
$\widetilde{\Omega}^*$: $$ \widetilde{\Phi}(K):= \Phi(\R \,{\rm
TR}_{<e>} (K)).$$ The forthcoming Lemma \ref{famous} states that
for any $K \in \widetilde{\Omega}^*$, $$\R \,{ \rm STR}_{<e>}\, [
\nabla, K] = d \,\R \,{ \rm STR}_{<e>} \, K\quad\text{ modulo
graded commutators}.$$ Thus, for any $K \in \widetilde{\Omega}^*$,
$\widetilde{\Phi}(\widetilde{\nabla} (K) ) = 0$. Next we recall
Connes' construction (\cite{Co} page 229) of a cycle $(
\widetilde{\widetilde{\Omega}}^*, d,
\widetilde{\widetilde{\Phi}})$ over $\widetilde{\cal{U}}=
\widetilde{\Omega}^0$. Let $X$ be a new odd formal variable of
degree $1$ and put: $$
\widetilde{\widetilde{\Omega}}^*=\widetilde{\Omega}^* \oplus
X\widetilde{\Omega}^* \oplus \widetilde{\Omega}X^* \oplus
X\widetilde{\Omega}^* X $$ with the new multiplication rules
$(\widetilde{\omega}_1 X) \widetilde{\omega}_2 =
\widetilde{\omega}_1 (X \widetilde{\omega}_2)=0$ and
$(\widetilde{\omega}_1 X )(X \widetilde{\omega}_2 ) =
 \widetilde{\omega}_1 \Theta \widetilde{\omega}_2.$ Then define
 a new graded trace on $\widetilde{\widetilde{\Omega}}^*$ by the formula:
 $$
 \widetilde{\widetilde{\Phi}}(\widetilde{\omega}_1 + X \widetilde{\omega}_2 +
 \widetilde{\omega}_3 X + X \widetilde{\omega}_4 X )=
 \widetilde{\Phi}(\widetilde{\omega}_1) + (-1)^{\partial\widetilde{\omega}_4}
 \widetilde{\Phi}( \widetilde{\omega}_4 \Theta).
 $$
 Define a differential $d$ on $\widetilde{\widetilde{\Omega}}^* $
 which is generated by the relations
 $$
 d \widetilde{\omega} = \widetilde{\nabla}(\widetilde{\omega} ) +
  X \widetilde{\omega} +
  (-1)^{\partial \widetilde{\omega}} \widetilde{\omega} X
 $$ and $d X=0$. One can check that $d^2=0$ and that
 $\widetilde{\widetilde{\Phi}}( d \widetilde{\widetilde{\omega}} ) =0$
 for any $\widetilde{\widetilde{\omega}}\in \widetilde{\widetilde{\Omega}}^* $.
 Therefore,
 $( \widetilde{\widetilde{\Omega}}^*, d, \widetilde{\widetilde{\Phi}})$ defines
 a cycle over
$\widetilde{\cal{U}}= \widetilde{\Omega}^0$ of dimension equal to
the degree of $\widetilde{\Phi}$. Let $ \widehat{\tau}_{\Phi}$ be
cyclic cocycle of $\widetilde{\cal{U}}$ defined by the character
(see \cite{Co} page 186) of this cycle.
 Let $\widetilde{\cal{U}}^+$
be the algebra obtained by adding a unit to $\widetilde{\cal{U}}$
with canonical homomorphism $\widehat{\pi}:
\widetilde{\cal{U}}^+\rightarrow \CC.$ Let us recall how $
\widehat{\tau}_{\Phi}$ induces a cyclic cocycle, still denoted $
\widehat{\tau}_{\Phi}$, on $\widetilde{\cal{U}}^+$. Set
 ${\widetilde{\widetilde{\Omega}}^0}_+ =
 \widetilde{\Omega^0} \oplus \CC= \widetilde{\cal{U}}^+$ and
 ${\widetilde{\widetilde{\Omega}}^k}_+={\widetilde{\widetilde{\Omega}}^k}$ for any $k\in \NN^\ast$.
 Then define $d:{\widetilde{\widetilde{\Omega}}^0}_+
 \rightarrow \widetilde{\widetilde{\Omega}}^1$
 by the formula $d(\omega_0 \oplus \lambda) = d\omega_0$ (ie, $d \lambda=0$ for
 any $\lambda \in \CC$),
 the differential from ${\widetilde{\widetilde{\Omega}}^k}_+$ to
 $\widetilde{\widetilde{\Omega}}^{k+ 1}_+$
 remaining unchanged for $k \geq 1$. If the degree of $ \widetilde{\Phi}$
 is zero
 then $\widehat{\tau}_{\Phi} $ is a trace on ${\widetilde{\Omega}^0}$ which
 extends as a trace on  ${\widetilde{\widetilde{\Omega}}^0}_+$ by the formula
 $ \widehat{\tau}_{\Phi}( \omega_0 \oplus \lambda) =\widehat{\tau}_{\Phi}(\omega_0)$.
If the degree of $ \widetilde{\Phi}$ is equal to $n \geq 1$ then
 one gets
 a cycle
 $( \widetilde{\widetilde{\Omega}}^*_+, d, \widetilde{\widetilde{\Phi}})$
 over $\widetilde{\cal{U}}_+$ which  allows to define
 a cyclic cocycle still denoted $\widehat{\tau}_{\Phi} $ on
$\widetilde{\cal{U}}^+$ by the formula: $$ \widehat{\tau}_{\Phi}(
\omega_0 \oplus \lambda_0, \cdots, \omega_n \oplus \lambda_n)=
\widehat{\tau}_{\Phi}( \omega_0, \cdots, \omega_n ). $$

 Consider
$p$ a projection in $M_2(\widetilde{\cal{U}}^+)$ such that
 $\widehat{\pi}( p)=\widehat{\pi}(p_0)$ where
$$p_0=\begin{pmatrix}  0 & 0\cr
0 & 1 \cr
\end{pmatrix},
$$
 then Connes has defined
the pairing $< [p-p_0], [\widehat{\tau}_{\Phi}]>$ (\cite{Co} page
224) between the $K-$theory class $[p-p_0] $ and the cyclic
cocycle $ \widehat{\tau}_{\Phi}$. We shall  use the following
lemma.
\begin{lemma}\label{connespairing}
 The operator $pe^{(p\nabla p)^2}p -p_0
e^{(p_0 \nabla p_0)^2} p_0 $ belongs to
$M_2(\widetilde{\Omega}^\ast)$ and there exists a constant $C(
{\rm deg} \,\Phi)$ depending only on the degree of $\Phi$ such
that: $$
 C( {\rm deg} \,\Phi)
 < [p-p_0], [\widehat{\tau}_{\Phi}]>
 =
<\, \R \,{\rm TR}_{<e>}\left( p e^{(p\nabla p)^2}p -p_0
e^{(p_0 \nabla p_0)^2} p_0 \right) , \Phi >
$$ where the definition of
${\rm TR}_{<e>}$ is given as in Definition \ref{inte} simply replacing there
$Str$ by $tr \otimes Tr_{M_2(\CC)}$.
\end{lemma}
\begin{proof} Let $(e_{i,j})_{1 \leq i, j \leq 2}$ be the basis
formed with elementary matrix of $M_2(\CC)$. Write $p=\sum_{i,j}
p_{i,j} e_{i,j}$. Then recall that: $$ \widehat{\tau}_{\Phi}
\sharp Tr (p, \cdots ,p) = \sum_{i_0,j_0,\cdots, i_n,j_n}
\widehat{\tau}_{\Phi}(p_{i_0,j_0},\cdots, p_{i_n,j_n}) \, Tr (\,
e_{i_0,j_0}\cdots e_{i_n,j_n}\,) . $$ Next using the fact (see
\cite{Go-Lo} just before (5.5)) that $p\, dp dp \,p $ equals the
curvature of $p\circ \nabla \circ p$ one proves easily the Lemma.

\end{proof}

\n

Now let $P$ be a parametrix of $D^+$ then there exists $\delta >0$
such that $S_+= \Id -P D^+$ and $S_-=\Id -D^+ P$ both belong
to $\Psi^{-\infty, \delta}_{\Ti}(\wm , \we^+) \simeq
 \Psi^{-\infty, \delta}_{\Ti}(\wm , \we^-)$ and,
in particular, $I(S_\pm ; \lambda) \equiv 0$.
As in \cite{Go-Lo} (section 5) and \cite{Co}
 we consider the two projections $p$ and
$p_0$ defined by:
$$
p=\begin{pmatrix} S^2_+ & S_+(\Id + S_+)P \cr
S_-D^+ & \Id-S^2_- \cr \end{pmatrix},\; \;
p_0=\begin{pmatrix} 0 & 0 \cr
0 & 1 \cr \end{pmatrix}.
$$ Recall (see \cite{Co} page 131) that the $K-$theory class
$ \Ind D^+\in K_0(\Psi^{-\infty, \delta}_{\Ti}( \wm , \we^+))$ is defined to be equal to
$[p-p_0] $.

Being motivated by section 5 of \cite{Go-Lo}  we give the
following:
\begin{definition}\label{Ch} We define  $\ch_{<e>} \Ind D^+$ to be equal to
$$\ch_{<e>} \Ind D^+ =\R \,{\rm TR}_{<e>}\left( pe^{(p\nabla p)^2}p -p_0
e^{(p_0 \nabla p_0)^2} p_0 \right).
$$
\end{definition}

\begin{definition}\label{hindex-d}
Let $\Phi$ a closed graded trace on $\widehat{\Omega}_*
(T,\Bi_\Gamma)$, concentrated on the trivial conjugacy class.
Let $\widetilde{\Phi}$ be its extension to
$\Psi^{-\infty,\delta}_{\widehat{\Omega}_*
(T,\Bi_\Gamma)}.$
We define the higher $\Phi$-index of $D^+$ as
$$
\Ind_\Phi (D^+)\,:=\,<{\rm Ch}_{<e>} (\Ind D^+)\,;\,\widetilde{\Phi} >$$
\end{definition}

\n
Lemma \ref{connespairing} shows that this definition is compatible with
the one of Connes.

\subsection{The Chern character of a superconnection.}\label{chernsupersub}$\;$

\medskip
We consider again
 $$\nabla= \nabla^u + \sum_{\gamma \in \Gamma} d\gamma \otimes h \gamma^{-1}$$ and
we define a new connection $\nabla^\prime = \nabla^{\prime,+} \oplus
\nabla^{\prime,-}$ acting on $ C^\infty_{\Ti}(\wm, \we)$ by:
$$
\nabla^{\prime,+}=S_+\nabla^{+} S_+ + P(\Id + S_-)\nabla^{-}D^+,
\; \nabla^{\prime,-}=\nabla^{-}.
$$
\n As in section 5 of \cite{Go-Lo} we consider
the superconnection $A(t)$ defined for $t\geq 0$ by
$$
A(t)=
\begin{pmatrix}
\nabla^{\prime,+} & t D^- \cr
t D^+ & \nabla^{\prime,-}
\end{pmatrix}
$$
On the one hand this connection is closely
related to the operators appearing in the definition of the index class;
on the other hand, precisely for this reason, the boundary behaviour
of $\nabla^\prime$ is rather complicated and it is a source of
technical complications in the proof of our main theorem.

\begin{definition}\label{chern-super}
For each $t\geq 0$
we define the Chern
character of $A(t)$, with values in
$\overline{\widehat{\Omega}}_* (T,\Bi_\Gamma)$,
 as  follows
\begin{equation}\label{chern-super-f}
\ch_{<e>}\,A(t)=\R\,\bigl(  TR_{<e>} S_+ e^{-A^2(t)}_{1,1} S_+ +
TR_{<e>} ( D e^{-A^2(t)}_{1,1} P(\Id + S_-) - e^{-A^2(t)}_{2,2} )
\bigr).
\end{equation}
\end{definition}
Notice that the first summand on the right hand side is an element
in $\Psi^{-\infty,\delta}_{\widehat{\Omega}(T,\Bi_\Gamma)}$,
thus the trace is well defined.
A simple computation shows that the operator $( D e^{-A^2(t)}_{1,1} P(\Id + S_-) - e^{-A^2(t)}_{2,2} )$, which a priori belongs to
$\Psi^{-\infty,\delta}_{b,\widehat{\Omega}(T,\Bi_\Gamma)}$, has vanishing
indicial family; thus the trace of the second summand on the right
hand side also exists.

We shall see below, in Proposition \ref{closed}, that $\ch_{<e>}\, A(t)$
is closed for $t\geq 0$.

\begin{lemma}\label{first}
The following formula holds:
\begin{equation}\label{first-f}
\ch_{<e>}\,A(0) = \ch_{<e>}\,\Ind D^+.
\end{equation}
\end{lemma}

\begin{proof}
This follows from  Definition \ref{Ch}
and formula  (5.29) in  \cite{Go-Lo}.
\end{proof}

\n
It is important to observe that
for $t>0$ the noncommutative differential form
$\R\,^{b} STR_{<e>}e^{-A^2(t)}$ is well defined.

\subsection{Main theorem and strategy of the proof.}\label{mainthsub}$\;$

\medskip
The main result of this paper
 is the following higher APS-index theorem

\n\begin{theorem} \label{higherindexIII}
\item {1)}
The following formula holds in the
homology $\widehat{H}_\ast(T, \Bi_\Gamma)$ (see Subsection \ref{heta1}): $$
{\rm Ch}_{<e>}\, (\Ind D^+)\,=\, \int_Z \phi
\,\widehat{A}\,(\nabla^{TZ})\,{\rm ch}(\nabla^{\widehat{V}})
\,{\rm ch}(\nabla^{can}) - {1\over 2}  \widetilde{\eta}_{<e>} \in
\widehat{H}_\ast(T, \Bi_\Gamma)
$$
\item{2)}
Let $\Phi$ a closed graded trace on $\widehat{\Omega}_*
(T,\Bi_\Gamma)$, concentrated on the trivial conjugacy class.
Then there exists a current $\omega_\Phi$ on $M$ such that
the following formula holds:
$$
 \Ind_\Phi (D^+)  =
<\widehat{A}(T\cal{F})\,{\rm ch}(\nabla^{V}) \,,\;\,
\omega_\Phi >- <{1\over 2}  \widetilde{\eta}_{<e>}\,;\,
\Phi >
$$
\end{theorem}

\smallskip
\noindent
{\bf Strategy of the proof of Theorem \ref{higherindexIII} .}\\
The proof of 1) can be divided into 4 steps.

\smallskip
\n {\bf Step 1.} For any real $t>0$ one has $$ \R ^{b} STR_{<e>}
e^{-\AA^2_t}=\int_Z \phi
\,\widehat{A}\,(\nabla^{TZ})\,{\rm ch}(\nabla^{\widehat{V}})
\,{\rm ch}(\nabla^{can}) - {1\over  \sqrt{\pi} }  \int_0^t
\widetilde{\eta}_{<e>}(s) ds.
$$ modulo an exact form.

\smallskip
\n {\bf Step 2.}
For any real $t>0$ one has
$$\R\, ^{b} STR_{<e>}
e^{-A^2(t)} =
\R\, ^{b} STR_{<e>}
e^{-\AA^2_t}+ B_1(t)$$ modulo an exact form, with $ B_1(t)$
a boundary term which goes  to zero as
  $t \rightarrow +\infty$ in $\overline{\widehat{\Omega}}_*(T, \Bi_\Gamma)$
endowed with the Hausdorff quotient topology induced by the
Frechet topology of ${\widehat{\Omega}}_*(T, \Bi_\Gamma)$.

\smallskip
\n {\bf Step 3.}
For $t>0$ one has $$\ch_{<e>}\,A(t)\,-\,\R\, ^{b} STR_{<e>}
e^{-A^2(t)}=B_2(t)$$ with $B_2(t)$ a boundary term going to zero
as $t \rightarrow +\infty$.

\smallskip
\n {\bf Step 4.} For any real $t>0$,
$$
\ch_{<e>}\,A(t)=\ch_{<e>}\,A(0)\,(=\ch_{<e>} \Ind D^+)
$$ modulo an exact form.

\smallskip
\n
{\bf Final step.}\\
First of all, one deduces 2) from 1) exactly
as explained in the proof of Corollary \ref{corocurrent}.\\
In order to prove 1) we proceed as follows:
by Step 4, Step 3, Step 2, Step 1 we have for each $t>0$, modulo exact forms,
\begin{align}
\ch_{<e>}\,\Ind D^+\:&=\: \ch_{<e>} A(t)\notag
\\ &=\: \R\,^{b} STR_{<e>}
e^{-A^2(t)}+B_2(t)\notag
\\ &=\:\R\, ^{b} STR_{<e>}
e^{-\AA^2_t}+ B_1(t)+B_2(t)\notag
\\ &=\:
\int_Z \phi
\,\widehat{A}\,(\nabla^{TZ})\,{\rm ch}(\nabla^{\widehat{V}})
\,{\rm ch}(\nabla^{can}) - {1\over  \sqrt{\pi} }  \int_0^t
\widetilde{\eta}_{<e>}(s) ds + B_1(t)+B_2(t)\notag
\end{align}
with $B_1(t)$ and $B_2 (t)$ additional boundary terms
that go to zero as $t\rightarrow +\infty$.
The index formula thus follows by taking $t\rightarrow +\infty$.
This finishes the sketch of the proof of the main theorem.

\medskip
We shall not use the following result but it is interesting (and
reassuring) in itself.
\begin{proposition}\label{closed}
For each real $t\geq 0$, $\ch_{<e>} A(t)$ is closed.
\end{proposition}

\begin{proof} (Sketch) Lemma \ref{famous} implies that
$$
d \R\,^{b} STR_{<e>}
e^{-\AA_t^2} = \R\,^{b} STR_{<e>}\,[\nabla, \,
e^{-\AA_t^2} ] = \R\,^{b} STR_{<e>}\,[\nabla -\AA_t, \,
e^{-\AA_t^2} ].
$$ Using Proposition \ref{commutator} and the fact that $D_0$ is invertible,
one checks (after an easy computation) that
$$\lim_{t\rightarrow +\infty}\, \R\,^{b} STR_{<e>}\,[\nabla -\AA_t, \,
e^{-\AA_t^2} ]\,=0.
$$ So $ \lim_{t\rightarrow +\infty}\,d \R\,^{b} STR_{<e>}
e^{-\AA_t^2} = 0$ and, Steps 2 and 3  imply that
$ \lim_{t\rightarrow +\infty}\,d \ch_{<e>} A(t) = 0$. The Lemma
is then a consequence of Step 4.
\end{proof}

\section{{\bf Proof of the main theorem}}\label{mainth}

We shall first prove
Step 3.
We shall then give a detailed proof of  Step 4; this is the most intricate step
and will require 4 lemmas. Then we shall only sketch the proof of Step 2,
 since it is similar to Step 4 but simpler.
Notice that Step 1 has already been proved.

\subsection{Proof of Step 3}\label{mainth1}$\;$

\medskip
We consider the noncommutative differential form
$$
\R\,^{b} TR_{<e>} \,[D^+,e^{-A^2(t)}_{1,1} P]\,+ \,
\R\,^{b} STR_{<e>}e^{-A^2(t)}.
$$
(Notice that it would be more rigorous to write
$$
\R\,^{b} TR_{<e>} D^+ e^{-A^2(t)}_{1,1} P -
\R\,^{b} TR_{<e>}e^{-A^2(t)}_{1,1} P D^+
$$ instead of
$\R\,^{b} TR_{<e>} \,[D^+,e^{-A^2(t)}_{1,1} P]$ as the operators
involved are not endomorphisms but rather homomorphisms from
the section of $\we^\pm$ to the sections of $\we^\mp$. )

\begin{lemma}\label{b}
For each real $t> 0$ one has:
$$
\ch_{<e>}\,A(t)=\R\,^{b} TR_{<e>} \,[D^+,e^{-A^2(t)}_{1,1} P]\,+ \,
\R\,^{b} STR_{<e>}e^{-A^2(t)}.
$$
\end{lemma}
\begin{proof} It is a simple computation using the definition
of $S_\pm$ and the fact that the indicial family of
$(D e^{-A^2(t)}_{1,1} P(\Id + S_-) - e^{-A^2(t)}_{2,2} ) $
is identically zero. The easy details are left to the reader.
\end{proof}

Using Proposition \ref{commutator},  an easy computation shows that:
 \begin{lemma}\label{extrab}
\begin{equation}\label{extrab-f}
\R\,^{b} TR_{<e>} \,[D^+,e^{-A^2(t)}_{1,1} P] =
 -{1 \over 2 \pi} \R TR \bigl(\,e^{-(\nabla_\partial^{-})^2-t^2 D_0^2}\,
 \int_\RR e^{-t^2\lambda^2} (D_0 + \sqrt{-1} \lambda \Id)^{-1} \,
 d\lambda \,\bigr).
 \end{equation}
 \end{lemma}

 \begin{lemma}\label{limitb2}
As $t\rightarrow +\infty$
\begin{equation}\label{limitb2-f}
\R TR \bigl(\,e^{-(\nabla_\partial^{-})^2-t^2 D_0^2}\,
 \int_\RR e^{-t^2\lambda^2} (D_0 + \sqrt{-1} \lambda \Id)^{-1} \,
 d\lambda \,\bigr)\longrightarrow 0
 \end{equation}
in the space of noncommutative differential forms $\overline{\widehat{\Omega}}_*
(T,\Bi_\Gamma)$.
 \end{lemma}

\begin{proof} Since $D_0$ is invertible,
  the following
finite propagation speed estimates (see \cite{Lott II} page 215)
are valid. There exists $\delta>0$ such that:
\begin{equation} \label{fp}
\forall s \geq 1,\; \forall (y,z) \in \partial \wm \times_\pi
\partial \wm, \,\; |e^{-s^2D_0^2}|(y,z) \leq C(N) (1 + d(y,z))^{-N}
e^{-s^2 \delta}.
\end{equation}
Now, using Duhamel's formula, we may write: $$
e^{-(\nabla_\partial^{-})^2-t^2 D_0^2}= e^{-t^2 D_0^2} + \int_0^1
e^{-u t^2 D_0^2} (-(\nabla_\partial^{-})^2) e^{-(1-u) t^2
D_0^2}\,du+ \cdots + $$ $$ \int_{\Delta_k} e^{-u_0 t^2 D_0^2}
(-(\nabla_\partial^{-})^2)
 e^{-u_1 t^2 D_0^2}(-(\nabla_\partial^{-})^2) \cdots
e^{-u_k t^2 D_0^2}  d\sigma(u_0,\cdots, u_k) + \dots $$ where
$\Delta_k$ denotes the $k-$simplex:
$$ \Delta_k=\{(u_0,\cdots,
u_k)\in [0, 1]^{k+1} /,\; \sum_{j=0}^k u_j=1 \}. $$ It is clear
that  for any $ (u_0,\cdots, u_k)\in \Delta_k$, at least one of
the $u_j$ satisfies $u_j \geq \frac {1} {k+1}$. Observe moreover
that standard  finite propagation speed estimates allow to show
that the Schwartz kernel of $(D_0 + \sqrt{-1} \lambda \Id )^{-1}$
is rapidly decreasing outside the diagonal. Using the estimates
\eqref{fp}, one checks easily that $$ K_k(t)= \int_{\Delta_k}
e^{-u_0 t^2 D_0^2} (-(\nabla_\partial^{-})^2)
 e^{-u_1 t^2 D_0^2}(-(\nabla_\partial^{-})^2) \cdots
e^{-u_k t^2 D_0^2}  d\sigma(u_0,\cdots, u_k) + \dots $$ belongs to
the space of Definition \ref{opdecay} and that the supremum
constants (for $ K_k(t) $) of the last line of  Definition
\ref{opdecay} are lower than $C e^{- \frac { \delta t^2} {2k+2}}$.
For a bit more details the reader may read the proof of Theorem
2.9 of \cite{LPMEMOIRS}. One then gets immediately the Lemma.

\end{proof}

\smallskip
\n
Step 3 follows from the above three lemmas.

\subsection{Proof of Step 4.}\label{mainth2}$\;$

\medskip
 The following lemma appears implicitly
in the literature, we feel it is appropriate to give a detailed proof.
We thank Sasha Gorokhovsky for useful explanations.

\begin{lemma} \label{famous}  Let
$K $ be an element of $\Psi^{-\infty}_{b,\widehat{\Omega}_* (T,\Bi_\Gamma)}(\wm ; \we)$.
 Then $$^{b}STR_{<e>}\,[\nabla, K] =
d  \,^{b}STR_{<e>}\, K $$ modulo graded commutators where we recall that
$\nabla = \nabla^u + \sum_{\gamma \in \Gamma} d\gamma
\otimes h \gamma^{-1}$.

\end{lemma}
\begin{proof} Recall that the differential $d$ of $\OM$ may be written
as $d=d_T + d_\Gamma$ where $d_T$ [resp. $d_\Gamma$] is the differential
corresponding to $\Omega^\ast(T)$ [resp. $\Omega^\ast(\Bi_\Gamma)$]. The fact
that $^{b}STR_{<e>}\,[\nabla^u , K] = d_T \,^{b}STR_{<e>}\, K $ is basically
proved in \cite{MP I}. So we have to prove that
$ ^{b}STR_{<e>}\,[A , K] = d_\Gamma \,^{b}STR_{<e>}\, K $ where we have set
$A=\sum_{\gamma \in \Gamma} d\gamma
\otimes h \gamma^{-1}$.

To shorten the notations we
 write $K(z,w)=\sum_{\omega} \omega K_\omega(z,w)$
 where it is understood (as in section
 \ref{bsuper1}) that $\omega=d g_1 \cdots d g_k$ where $g_1, \ldots , g_k \in \Gamma. $
 We set $\pi_1(\omega)= g_1 \cdots  g_k.$
Since $K$ is (graded-)$\widehat{\Omega}_* (T,\Bi_\Gamma)-$linear,
we observe that
$KA$ can be written as
$ (-1)^{|\omega|}\sum_{\gamma, \omega} d\gamma \gamma^{-1}
\omega K_{\omega}h^{\gamma}$ where $ K_{\omega}h^{\gamma}$ denotes the kernel
$K_{\omega}(z,w)h^{\gamma}(w)$ and $ h^{\gamma}(w)=h(w\cdot \gamma)$ .
  Now using the definition of the action of
  the connection $A$ on acting
  on $\widehat{\Omega}_* (T,\Bi_\Gamma) \otimes_{\Ti} C^\infty_{\Ti}(\wm,\we)$
  , one sees that $AK$ can be written as
$ \sum_{\gamma, \omega}(-1)^{|\omega|}\omega d\gamma \gamma^{-1}h^{\gamma}K_{\omega}$,
where $h^{\gamma}K_{\omega}$ denotes the kernel $h^{\gamma}(z)K_{\omega}(z,w)$.

Let $H$ be an operator (not necessarily $\Ti-$linear, such as
$A K $ or $K A$) given by a Schwartz kernel  $H(z,w)$  of the form
$H(z,w)= \sum_{\widetilde{\omega}}
\widetilde{\omega}H_{\widetilde{\omega}}(z,w)$ where it is understood that
$\widetilde{\omega}=g_0 dg_1\cdots dg_k$ and such that
the  $ H_{\widetilde{\omega}}(z,w)$ satisfy the decay estimates
of Definition \ref{12}.  Then  we set for any $\theta \in T$:
$$
TR_{<e>} H(\theta) =
\sum_{\widetilde{\omega}} \widetilde{\omega} \pi_1(\widetilde{\omega})^{-1}
\nuint H_{\widetilde{\omega}}
(w\pi_1(\widetilde{\omega}),w) \phi(w) d Vol^b_{\pi^{-1}(\theta)}(w)
$$ where $\pi_1(\widetilde{\omega})=g_0 g_1 \cdots g_k.$

Notice that in the formula expressing $K A$ we have
$d\gamma \gamma^{-1}=-\gamma^{-1}d\gamma$.
Hence for any $\theta \in T$ we get:

\begin{align}
TR_{<e>} [K,A](\theta)  \:&=\: \sum_{\gamma, \omega}\,
  (-1)^{|\omega|}d\gamma \gamma^{-1} \omega\pi_1(\omega)^{-1} \,
  \nuint  \phi(w) K_{\omega}(w\pi_1(\omega), w) h^{\gamma}(w)
  d Vol^b_{\pi^{-1}(\theta)}(w)
\notag
 \\ &-\: \sum_{\gamma, \omega}\omega d\gamma \gamma^{-1}\pi_1(\omega)^{-1} \nuint \phi(w)
  h^{\gamma}(w\pi_1(\omega)) K_{\omega}(w\pi_1(\omega), w)
  d Vol^b_{\pi^{-1}(\theta)}(w)
\notag
\end{align}
and modulo graded commutators this equals to
$$
  \sum_{\gamma, \omega} \omega\pi_1(\omega)^{-1} d\gamma \gamma^{-1}
  \nuint \phi(w)
  K_{\omega}(w\pi_1(\omega), w) h^{\gamma}(w)d Vol^b_{\pi^{-1}(\theta)}-
  $$
  $$
  \sum_{\gamma, \omega}\omega d\gamma \gamma^{-1}\pi_1(\omega)^{-1} \nuint \phi(w)
  h^{\gamma}(w\pi_1(\omega)) K_{\omega}(w\pi_1(\omega), w)
  d Vol^b_{\pi^{-1}(\theta)}
  $$

The second term, after the renaming $\gamma^\prime= \pi_1(\omega) \gamma$
(then replace  $\gamma^\prime$ by $\gamma$)  can be written as

  $$ \sum_{\gamma, \omega}\omega
  d(\pi_1(\omega)^{-1}\gamma) \gamma^{-1} \nuint \phi(w)h^{\gamma}(w) K_{\omega}
  (w\pi_1(\omega), w) d Vol^b_{\pi^{-1}(\theta)}=
  $$
  $$
  \sum_{\gamma, \omega}\omega
  d(\pi_1(\omega)^{-1})  \nuint \phi(w)h^{\gamma}(w)
  K_{\omega}(w\pi_1(\omega), w)d Vol^b_{\pi^{-1}(\theta)}+
  $$
  $$
  \sum_{\gamma, \omega}\omega \pi_1(\omega)^{-1}
  d\gamma \gamma^{-1} \nuint\phi(w) h^{\gamma}(w)
   K_{\omega}(w\pi_1(\omega), w)d Vol^b_{\pi^{-1}(\theta)}
  $$

Hence $TR_{<e>} [K,A](\theta)=-
\sum_{\gamma, \omega} \omega d(\pi_1(\omega)^{-1})  \nuint\phi(w)
h^{\gamma}(w)
K_{\omega}(w\pi_1(\omega), w)d Vol^b_{\pi^{-1}(\theta)}
=(-1)^{|\omega|+1}d TR_{<e>}(K)$, as $\sum_{\gamma \in
\Gamma} h^{\gamma}=1$ and $ d ( \omega \pi_1(\omega)^{-1}) =
(-1)^{|\omega|} \omega d\pi_1(\omega)^{-1} $.
\end{proof}

Now we recall from formulas (5.25) and (5.26) of \cite{Go-Lo} that
modulo operators with vanishing indicial families we have the
following two identities:

\begin{equation} \label{5.25}
A^2(t) \equiv \begin{pmatrix}   (\nabla^{\prime,+})^2 + t^2 D^-D^+ &
t[\nabla^{\prime,-}, D^-] + t(\nabla^{\prime,+} -\nabla^{\prime,-})D^- \cr
0 & D^+(\,(\nabla^{\prime,+})^2 + t^2 D^- D^+ \,) P \cr
\end{pmatrix}
\end{equation}

\begin{equation} \label{5.26}
e^{-A^2(t)} \equiv \begin{pmatrix} e^{-(\nabla^{\prime,+})^2 -t^2 D^-D^+} &
{\cal {Z}} \cr
0 & De^{-(\nabla^{\prime,+})^2 - t^2 D^-D^+}P \cr
\end{pmatrix}
\end{equation}
with
${\cal {Z}}$  given  by
$$-\int_0^1 e^{-u(\,(\nabla^{\prime,+})^2 + t^2 D^-D^+\,)}\left(t[\nabla^{\prime,-},D^-]+
t(\nabla^{\prime,+} - \nabla^{\prime,-})D^+\right)
e^{-(1-u)(\,(\nabla^{\prime,+})^2+ t^2 D^+D^-\,)}du$$

\bigskip

\begin{lemma} \label{der}
The following formula holds:
\begin{align}
&{ d \over dt} {\cal{R}}\, ^{b} STR_{<e>}\, e^{-A^2(t)}\:=\:
-{\cal{R}}\, ^{b}STR_{<e>}\,[A(t), A^\prime(t)e^{-A^2(t)}] \notag
\\&+\:
{\cal{R}} \int_\RR {t D_0\over \pi} e^{-(\,(\nabla^-_\partial)^2 + t^2 D_0^2 +
t^2 \lambda^2\,) }
d\lambda -
{ d \over dt} {\cal{R}} ^{b} TR_{<e>}[D, e^{-A^2(t)}_{1,1} P].
\notag
\end{align}

\end{lemma}
\begin{proof} Using Duhamel formula and Proposition \ref{commutator} one gets:
$${ d \over dt} {\cal{R}} ^{b} STR_{<e>}\, e^{-A^2(t)}=
-{\cal{R}} ^{b} STR_{<e>}\,\int_0^1 e^{-u A^2(t)} (\,{ d \over dt} A^2(t)\,)
e^{-(1-u) A^2(t)} =
$$
$$
-{\cal{R}} ^{b}STR_{<e>}\,[A(t), A^\prime(t)e^{-A^2(t)}]-
$$
\begin{equation} \label{eq}
{i \over 2 \pi}{\cal{R}}
STR_{<e>}\,\int_0^1 \int_\RR \bigl( {\partial \over \partial \lambda}
I( e^{-uA^2(t)}; \lambda) I( { d \over dt} A^2(t); \lambda)
 I( e^{-(1-u)A^2(t)}; \lambda)\bigr) du d\lambda.
  \end{equation}
  Due to the definition of $\nabla^\prime$ we have:
 $$
 I^2(\nabla^{\prime,+}, \lambda) = (i\lambda + D_0)^{-1}
 (\nabla_\partial^-)^2 (i\lambda + D_0).
 $$
  From (\ref{5.25}) and (\ref{5.26}) one gets the two following formulas:
 $$
 I( { d \over dt} A^2(t); \lambda) =
 \begin{pmatrix} 2 D_0^2 t + 2 t \lambda^2 & \ast \cr
 0 & 2 D_0^2 t + 2 t \lambda^2 \cr
 \end{pmatrix}
 $$
 \begin{equation} \label{I}
 I( e^{-uA^2(t)}; \lambda) = \begin{pmatrix} (i\lambda + D_0)^{-1}
 e^{-u ((\nabla_\partial^-)^2  +u t^2D_0^2+ut^2\lambda^2 )}\,(i\lambda + D_0) & \ast
 \cr
 0 & e^{-u (\nabla_\partial^-)^2  -u t^2D_0^2-ut^2\lambda^2} \cr
 \end{pmatrix}
 \end{equation}
 Then, a little computation (using the two previous equations) allows to see that:
 $$-{i \over 2 \pi}{\cal{R}}
STR_{<e>}\,\int_0^1 \int_\RR \bigl( {\partial \over \partial \lambda}
I( e^{-uA^2(t)}; \lambda) I( { d \over dt} A^2(t); \lambda)
 I( e^{-(1-u)A^2(t)}; \lambda)\bigr) du d\lambda=
 $$
 $$
 {1 \over 2 \pi} {\cal R} TR_{<e>}\,
 \int_\RR (2D_0^2 + 2 t \lambda^2) ( i \lambda +
 D_0)^{-1} e^{- (\nabla_\partial^-)^2  - t^2D_0^2-t^2\lambda^2} \, du d\lambda -
 $$
 $$
 {1 \over 2 \pi}{\cal R} TR_{<e>}\,\int_\RR ( i \lambda +
 D_0)^{-1} \int_0^1  e^{- u((\nabla_\partial^-)^2  + t^2D_0^2+t^2\lambda^2 )}
 (2D_0^2 + 2 t \lambda^2)
 e^{-(1- u)((\nabla_\partial^-)^2  + t^2D_0^2+t^2\lambda^2 ) } du d\lambda.
 $$ Computing a little more one gets:
 $${1 \over 2 \pi}{\cal R} TR_{<e>}\, \int_\RR (2D_0^2 + 2 t \lambda^2) ( i \lambda +
 D_0)^{-1} e^{- (\nabla_\partial^-)^2  - t^2D_0^2-t^2\lambda^2} \, du d\lambda=
 $$
 $$
 {\cal R} TR_{<e>}\, \int_\RR {1 \over  \pi} t D_0
  e^{- (\nabla_\partial^-)^2  - t^2D_0^2-t^2\lambda^2} \, du d\lambda.
 $$ Then, using lemma \ref{extrab} and Duhamel formula, one gets that:
 $${1 \over 2 \pi}{\cal R} TR_{<e>}\,\int_\RR ( i \lambda +
 D_0)^{-1} \int_0^1  e^{- u((\nabla_\partial^-)^2  + t^2D_0^2+t^2\lambda^2\,)}
 (2D_0^2 + 2 t \lambda^2)
 e^{-(1- u)((\nabla_\partial^-)^2  + t^2D_0^2+t^2\lambda^2 )} du d\lambda=
 $$
 $$
 { d \over dt}{\cal{R}} ^{b} TR_{<e>}[D, e^{-A^2(t)}_{1,1} P].
 $$ Then, from the two previous equations and from (\ref{eq}) one
 obtains the Lemma.
\end{proof}

\begin{lemma} \label{super} The following formula holds:
\begin{align}
{\cal{R}} ^{b}STR_{<e>}\,[A(t), A^\prime(t)e^{-A^2(t)}]\:&=
d {\cal{R}} ^{b}STR_{<e>}\, \begin{pmatrix}
0 & D^- \cr
D^+ & 0 \cr
\end{pmatrix} e^{-A^2(t)} \notag\\
&\:+\: { 1\over \pi}
\int_\RR {\cal{R}} TR_{<e>}\, t D_0 \,e^{-(\nabla_\partial^-)^2 - t^2 D_0^2 -t^2
\lambda^2}\, d\lambda.\notag
\end{align}
\end{lemma}
\begin{proof} Using the very definition of $A(t)$ one gets:
${\cal{R}} ^{b}STR_{<e>}\,[A(t), A^\prime(t)e^{-A^2(t)}]=$
$$
{\cal{R}} ^{b}STR_{<e>}\,[\begin{pmatrix} 0 & t D^- \cr
t D^+ & 0 \cr
\end{pmatrix}\,,\, \begin{pmatrix} 0 & D^- \cr
D^+ & 0 \cr
\end{pmatrix} e^{-A^2(t)}] \, +
$$
$$
{\cal{R}} ^{b}STR_{<e>}\,[\nabla , A^\prime(t)e^{-A^2(t)}]
+
{\cal{R}} ^{b}STR_{<e>}\,[\nabla^\prime-\nabla ,
\begin{pmatrix} 0 & D^- \cr
D^+ & 0 \cr
\end{pmatrix} e^{-A^2(t)}].
$$ Now using the definition of $\nabla^\prime$ and (\ref{I}) (with $u=1$) one gets:
$$
\begin{pmatrix} {\partial \over \partial \lambda} I(\nabla^{\prime,+} -\nabla^+;
\lambda) & 0 \cr
0 & 0 \cr \end{pmatrix}
I( A^\prime(t) e^{-A^2(t)} ; \lambda) \equiv 0
$$ therefore  Proposition \ref{commutator} shows that
$${\cal{R}} ^{b}STR_{<e>}\,[\nabla^\prime-\nabla ,
\begin{pmatrix} 0 & D^- \cr
D^+ & 0 \cr
\end{pmatrix} e^{-A^2(t)}] \equiv 0.
$$ Now using Proposition \ref{commutator} and (\ref{I}) (with $u=1$),
a simple computation shows that:
$${\cal{R}} ^{b}STR_{<e>}\,[\begin{pmatrix} 0 & t D^- \cr
t D^+ & 0 \cr
\end{pmatrix}\,,\, \begin{pmatrix} 0 & D^- \cr
D^+ & 0 \cr
\end{pmatrix} e^{-A^2(t)}] =
$$
$$
{ i \over 2 \pi}
\int_\RR {\cal R} STR_{<e>}\,
 \begin{pmatrix} -i t D_0 + t \lambda & 0 \cr
0 & i t D_0 + t \lambda \cr \end{pmatrix}
I(e^{-A^2(t)} ; \lambda)\,d\lambda\, =
$$
$${ 1\over \pi}
\int_\RR {\cal{R}} TR_{<e>}\, t D_0 \,e^{-(\nabla_\partial^-)^2 - t^2D_0^2 -t^2
\lambda^2}\, d\lambda.
$$ Lastly, Lemma \ref{famous} shows that:
$${\cal{R}} ^{b}STR_{<e>}\,[\nabla , A^\prime(t)e^{-A^2(t)}] = d {\cal R}
^{b}STR_{<e>}\, \begin{pmatrix} 0 & D^- \cr
D^+ & 0 \cr
\end{pmatrix}
e^{-A^2(t)}.
$$
The Lemma is proved.
\end{proof}

By combining the definition of $\ch_{<e>}\, A(t)$ and the Lemmas
\ref{der} and \ref{super} one gets
\begin{equation}\label{jls}
{d \over dt} \ch_{<e>} A(t) \,=
{d \over dt} (\, {\cal R} ^{b}TR_{<e>}\,[ D^+, e^{-A^2(t)}_{1,1} P] +
{\cal R} ^{b}STR_{<e>}\, e^{-A^2(t)}\,)\,=\,
\end{equation}
$$
 -
d {\cal R} ^{b}STR_{<e>}\,\begin{pmatrix} 0 & D^- \cr
D^+ & 0 \cr
\end{pmatrix} e^{-A^2(t)}.
$$
As observed in section 5.2 of \cite{Go-Lo},
${\cal R} ^{b}STR_{<e>}\,\begin{pmatrix} 0 & D^- \cr
D^+ & 0 \cr
\end{pmatrix} e^{-A^2(t)}$ is not integrable at $t=0$ but we are going to show
that we can add an exact form to it so that the sum does become
 integrable at $t=0$.
To this aim we shall need the following Lemma:
\begin{lemma} \label{heat} For any reals $s>0$ and $v>0$ one has:
$ [ \nabla^- ; e^{-s(v (\nabla^-)^2 + D^+ D^-)}\,]=$
$$
 -\int_0^s e^{-(s-u)(v (\nabla^-)^2 + D^+ D^-)} \,
[\nabla^- , v (\nabla^-)^2 + D^+ D^-] e^{-u(v (\nabla^-)^2 + D^+ D^-)} \,du.
$$
\end{lemma}
\begin{proof} Let us fix $v>0$. Then one checks easily that $X(s)=$
$$
[ \nabla^- ; e^{-s(v (\nabla^-)^2 + D^+ D^-)}] + \int_0^s e^{-(s-u)(v (\nabla^-)^2 + D^+ D^-)} \,
[\nabla^- , v (\nabla^-)^2 + D^+ D^-] e^{-u(v (\nabla^-)^2 + D^+ D^-)} du
$$ satisfies
$( {\partial \over \partial s} +  v (\nabla^-)^2 + D^+ D^- ) X(s)\equiv 0 $
and that $\lim_{s \rightarrow 0^+} X(s) = 0$. By a uniqueness argument
one then gets the Lemma.

\end{proof}

Now, following section 5.2 (and especially formulas
(5.31) and (5.32)) of \cite{Go-Lo} one sees that modulo
an integrable function at $t=0$, ${\cal R} ^{b}STR_{<e>}\,\begin{pmatrix} 0 & D^- \cr
D^+ & 0 \cr
\end{pmatrix} e^{-A^2(t)}$ is equal to
$$
{1 \over t} {\cal R} ^{b} TR_{<e>} \int_0^1 e^{-u( (\nabla^-)^2 + t^2 D^+ D^-)}
[\nabla^- ,  (\nabla^-)^2 + t^2 D^+ D^- ]
e^{-(1-u)( (\nabla^-)^2 + t^2 D^+ D^-)} du.
$$ But applying Lemma \ref{heat} with $s=t^2$ and $v={1\over t^2}$ and using
a linear change of variables to replace $\int_0^{t^2}$ by
$\int_0^{1}$ one sees that the previous term is equal to
$$
- {1 \over t} {\cal R} ^{b} TR_{<e>} [\nabla^-,
e^{-( (\nabla^-)^2 + t^2 D^+ D^-)}\,].
$$ Hence, using Lemma \ref{famous} one deduces that modulo graded commutators
 $$
{\cal R} ^{b}STR_{<e>}\,\begin{pmatrix} 0 & D^- \cr
D^+ & 0 \cr
\end{pmatrix} e^{-A^2(t)} + { 1 \over t}
d \, {\cal R} ^{b} TR_{<e>}e^{-( (\nabla^-)^2 + t^2 D^+ D^-)}
$$ is integrable at $t=0$. Therefore, having in mind (\ref{jls}), one sees
that we may write (up to graded commutators) that for any $t \in ]0,1],$
${d \over dt} \ch_{<e>} A(t) = d F(t) $ where
$t \rightarrow F(t)$ belongs $C^1( [0,1] ; \OM)$. Since
$t \rightarrow\ch_{<e>} A(t)$ is continuous at $t=0$ one obtains Step 4.

\subsection{Proof of Step 2}\label{mainth3}$\;$

The proof of the Step 2 consists in the
analysis of a
transgression formula for the homotopy $u\AA_t + (1-u) A(t)$
($0\leq u \leq 1$) as $t\rightarrow +\infty$. One takes care of the details
by using similar techniques as the ones used in the proof of  Step
 4 given above and also in section 14 of \cite{LPMEMOIRS}.

\section{{\bf General \'etale groupoids}}\label{genetale}

\medskip  We follow  Section 6 of \cite{Go-Lo}. Let
$G$ be a smooth Hausdorff etale groupoid with units $G^{(0)}$, recall
that the range maps $r:G\rightarrow G^{(0)}$ and
the source maps $s:G\rightarrow G^{(0)}$ are (by assumption) local
diffeomorphisms. We shall assume that $G^{(0)}$ is a
  smooth compact manifold possibly with boundary and that $G^{(0)}$
  is
  endowed with a  riemannian metric $g$ . Notice
  that the pull back by the map $ s:G\rightarrow G^{(0)}$
  of $g$ defines a   riemannian metric on $G$, we shall
  denote by $\Delta$ the corresponding  Laplace Beltrami
  operator on $G$.

  We consider a fiber bundle $\sigma: L\rightarrow G^{(0)}$ such that
  each fiber $L_\theta=\sigma^{-1}(\theta)$ is a complete length
  space with metric
  $d_\theta$. We assume that $G$ acts isometrically, properly and cocompactly
  on $L$. Let $i: G\rightarrow L$  be a $G-$equivariant map ($i$ is
  not necessarily
  continuous) sending $G_\theta$ into $L_\theta$ for any $\theta\in  G^{(0)}$.
   We assume
  that the preimage by $i$ of any compact set of $L$ has compact
  closure in $G$ and for any compact subset $K$ of   $G^{(0)}$ ($\subset G$),
  $i(K)$ has compact closure. One defines a length function
  on $G$ by
  $$
  \forall \gamma \in G,\; l(\gamma)= d_{s(\gamma)}(i(s(\gamma)), i(\gamma))
  $$ where we recall that $\gamma $ and $s(\gamma)$ belong to
  $G_{s(\gamma)}$. Moreover we assume  that
  \begin{equation} \label{poly-growth}
  \exists C, N \in \RR^+\; \forall \theta \in  G^{(0)}, \;
  \#\{\gamma \in G_\theta/\, l(\gamma) \leq R\} \leq C ( 1 + R)^N.
  \end{equation}
  This assumption is the generalization
  of the hypothesis "$\Gamma$ virtually nilpotent" in the case
  $G=T\rtimes \Gamma$, we recall that  Section 6.1 of \cite{Go-Lo}
   explains to which extent these assumptions generalize to $G$ the
   ones made for  $T\rtimes \Gamma$. Now we are going to define a
   smooth sub-algebra $\Ti(G) $ of $C^*_r(G)$
   which is the analogue for $G$  of
   $\Ti$ for $T\rtimes \Gamma$.
   We fix a Haar system in $G$.

\begin{definition}
   We denote by $\Ti(G) $ the sub-algebra of $C^*_r(G)$ whose elements are the
    functions $f\in C^\infty(G)$ such that
   for any $p, q \in \NN:$
   $$
   \sup_{\gamma \in G}
   \bigl( ( 1 + l(\gamma) + l(\gamma^{-1}))^q\,|\Delta^p f (\gamma)|  \bigr) <
   +\infty.
   $$
   \end{definition}

   \noindent {\bf Remark.} The fact  that $\Ti(G) $ is an algebra  is a consequence of the
   estimates: $ l(\gamma \cdot \gamma^\prime) \leq l(\gamma) + l(\gamma^\prime)$.

   Now for each $n\in \NN$, we denote (following \cite{Go-Lo}) by
   $G^{(n)}$ the set of $n-$chains of composable elements of $G$:
   $$
   G^{(n)}=\{(\gamma_1,\cdots,\gamma_n) \in G^n:\,
   s(\gamma_1)=r(\gamma_2),\cdots, s(\gamma_{n-1})=r(\gamma_n)\,\}.
   $$

   Then, for $(m,n) \in \NN^2$, we denote by $\widetilde{\Omega}^{m,n}(G)$
   the subset of sooth differential forms $\omega \in \Omega^m(G^{(n+1)})$
   such that for any $p,q \in \NN$:
   $$
   \forall p,q \in \NN,\; \sup_{(\gamma_0,\cdots,\gamma_n) \in G^{(n+1)}}
   \bigl( ( 1 + l(\gamma_0)+\cdots +l(\gamma_n) )^q \,
   |\Delta^p \omega (\gamma_0,\cdots,\gamma_n) |\bigr).
   $$ Then we denote by $\Omega^\ast(G)$ the quotient of
   $\Pi_{m,n\in \NN}\, \Omega^m(G^{(n+1)})$ by the set of differential forms
   which are supported on $\{(\gamma_0,\cdots,\gamma_n) \in G^{(n+1)},\;
   \gamma_j \;{\rm is\, a \, unit\, for \, some}\, j\, \}$.

   Now we consider a smooth $G-$manifold with boundary $P$ (see section
   II.10.$\alpha$ of \cite{Co}). That is, first of all,
   there is a fiber bundle $\pi:P\rightarrow G^{(0)}$ whose fibers are
   manifolds with boundary which are transverse to $\partial P$
   and of dimension $2k$. Given $\theta\in G^{(0)}$, we write
   $Z_\theta=\pi^{-1}(\theta)$. Putting
   $$
   P\times_{G^{(0)}}G = \{(p,\gamma) \in P \times G:\; p \in Z_{r(\gamma)}\},
    $$ we have a map $ P\times_{G^{(0)}}G \rightarrow P$, denoted
    $(p,\gamma) \rightarrow p\gamma$, such that
    $p\gamma \in Z_{s(\gamma)}$ and $(p\gamma_1)\gamma_2= p(\gamma_1\gamma_2)$
    for any $(\gamma_1 , \gamma_2) \in G^{(2)}$. We assume
    that $P$ is a proper $G-$manifold, ie that the map $P\times_{G^{(0)}}G
    \rightarrow P\times P$ given by
    $(p, \gamma) \rightarrow (p, p\gamma)$ is proper. Then the groupoid
    \def \G {\cal{G}}
    $\G= P\rtimes G$ with underlying space $P\times_{G^{(0)}}G $, with
    space of units $\G^{(0)}=P$ and maps $r(p,\gamma)=p$ and
    $s(p,\gamma)=p\gamma$, is a proper groupoid. We also assume
    that $G$ acts cocompactly and freely on $P$, equivalently
    $\G$ is a free and cocompact groupoid.
    Let us fix $h\in C^\infty_c( \G^{(0)}) $ (constant in the normal direction
    near the boundary) such that for all
    $\theta\in \G^{(0)}$
    $$
    \sum_{\gamma \in \G^\theta}\, h(s(\gamma)) =1.
    $$
    Then we recall that
    Gorokhovsky and Lott defined a connection
    $$
    \nabla^{can}:\;
    C^\infty_c(\G^{(0)}) \rightarrow \Omega^1_c(\G) \otimes_{C^\infty_c(\G)}
    C^\infty_c(\G^{(0)})
    $$ of the form $\nabla^{can}=\nabla^{1,0} \oplus\nabla^{0,1}$ where
    for any $F\in C^\infty_c(\G^{(0)})$,
    $\nabla^{1,0}(F) \in \Omega^1 (\G^{(0)} )$ is the de Rham differential of
    $F$ and $\nabla^{0,1}(F) \in \Omega^{0,1}_c (\G^{(0)} )
    \otimes_{C^\infty_c(\G)} C^\infty_c(\G^{(0)})$ is given by
    $$
    \forall \gamma_0 \notin \G^{(0)},\; \,(\nabla^{0,1}(F)\,)\,(\gamma_0)\,=\,
    F(r((\gamma_0) )\, h( s(\gamma_0)).
    $$

    Now we assume that the $2k-$dimensional fibers
     $\pi^{-1}(\theta),\, \theta \in G^{(0)}$
     carry a $G-$invariant spin structure and that there
    exists a
     $G-$invariant smooth family of $b-$metrics on the fibers
     $\pi^{-1}(\theta),\, \theta \in G^{(0)}$. We denote the typical fiber of
     $\pi$ by $Z$ and by $S^Z=(S^Z)^+ \oplus (S^Z)^-\rightarrow P$ the
     associated ($\ZZ_2-$graded) spinor bundle.
   We consider also
a $G-$equivariant complex hermitian vector bundle
$\widehat{V} \rightarrow P $ endowed with a
$G-$invariant $b-$hermitian connection $\widehat{\nabla}$ satisfying
$\widehat{\nabla}_{x \partial_x}=0$ on the boundary $\partial P$.
We then set $\widehat{E}= S^Z \otimes \widehat{V}= \widehat{E}^+ \oplus
\widehat{E}^-$ which defines a smooth family of  $\ZZ_2-$graded hermitian
Clifford modules in the fibers $\pi^{-1}(\theta), \, \theta \in G^{(0)}$.
Then we get a smooth family of $G-$invariant  $\ZZ_2-$graded Dirac
type operators
$$  D(\theta)= \begin{pmatrix} 0 & {D^-(\theta)}
\cr {D^+(\theta)}  & 0 \cr
\end{pmatrix},\; \theta\in G^{(0)}
$$
 acting fiberwise on
$C^\infty_c( P,\, \widehat{E})$.
 Moreover in a collar neighborhood ($ \sim [0,1] \times \partial \pi^{-1}(\theta) =
\{(x, y)\}$) of  $\partial \pi^{-1}(\theta) $ we may write:
$$
 {D^+(\theta)}= \sigma (x \partial_x + D_0(\theta)\,)
$$ where $ D_0(\theta)$ is the induced boundary Dirac type operator acting
on
$$
 C^\infty(\partial \pi^{-1}(\theta) ,\, \widehat{E}^+_{|_{\partial \pi^{-1}(\theta)}}\,).
$$

\medskip \n
In the rest of this paper we shall make the following

\noindent {\bf Hypothesis A} There exists a real $\epsilon >0$ such that for any
$\theta \in G^{(0)} $, the $L^2-$spectrum of $ D_0(\theta)$ acting on
$L^2(\partial \pi^{-1}(\theta) ,\, \widehat{E}^+_{|_{\partial \pi^{-1}(\theta)}}\,)
$ does not  meet $ ]-\epsilon, \epsilon[$.

\medskip
We set
$$
C^\infty_{\Ti(G)}(P,\we)= \Ti(G)\otimes_{C^\infty_c(G)}C^\infty_c(P,\we),\;
C^\infty_{\Ti(G)}(\partial P,\we) =
\Ti(G)\otimes_{C^\infty_c(G)}C^\infty_c(\partial P,\we_{|\partial P}).
$$ As in Section \ref{fibration}, the $(D(\theta))_{\theta \in G^{(0)}}$ [resp.
$(D_0(\theta))_{\theta \in G^{(0)}}$ ] induce a
a left $\Ti(G)-$linear endomorphism of $ C^\infty_{\Ti(G)}(P,\we)$
[resp. $C^\infty_{\Ti(G)}(\partial P,\we)$].

Now we fix a $G-$invariant horizontal distribution $T^H  P$ such that
$$
{{}}^bTP = T^H P \oplus {{}}^bT( P/T)
$$ and, as in \cite{Go-Lo} and \cite{MP I} (section 9) we
consider the Bismut superconnection
\begin{equation} \label{bismut}
\AA_s^{Bismut}=sD+\nabla^{u}-{1\over 4 s}c(\tau),\; s\in \RR^{+\ast}
\end{equation} where $D$ is the $\Ti(G)-$linear Dirac operator introduced above,
$c(\tau) $ denotes the Clifford multiplication by the curvature
2-form $\tau$ of
$T^H P $ and $\nabla^{u} $ is a certain unitary connection. Then,
as in \cite{Go-Lo}  we
 consider for each real $s>0$ the superconnection
$$\AA_s= \AA_s^{Bismut} + \nabla^{0,1}
$$ which sends $C^\infty_{\Ti(G)}(P,\we)$ into
$\Omega^\ast (G) \otimes_{\Ti(G)} C^\infty_{\Ti(G)}(P,\we)$.
By developing in a straightforward way a $b-$heat calculus as in
\cite{LPMEMOIRS} (section 10) one sees easily that
$e^{-s^2D^2} \in \Psi^{-\infty}_{b,{\cal T}^\infty(G)}
(P\,;\, \widehat{E})$ for any $s>0$ (the definition of this  space
 is completely analogous to the one given in Section \ref{bcalc}).
 Moreover, using a Duhamel expansion around
$e^{-s^2D^2}$ one checks that for any $s>0$
$$
e^{-\AA_s^2} \in \Psi^{-\infty}_{b,\Omega^\ast(G)}
(P\,;\, \widehat{E}).
$$

\bigskip Now as in \cite{MP I} (section 10) we consider for
any $s>0$  the induced boundary connection
$$
\BB_s=s \sigma D_0 + \nabla^{0,1}
 $$ where $D_0$ is the boundary Dirac operator of $D$
 introduced above, $\nabla^{u}_\partial$
is a certain unitary connection and $c(\partial \tau) $ is the boundary connection of
$c(\tau) $. We fix a function $\phi \in C^\infty_c(P)$ which is constant
in the normal direction near the boundary such that
$$
\forall p\in P,\,\;\; \sum_{\gamma \in G^{\pi(p)}}\, \phi(p \gamma)=1.
$$ Such a function $\phi$ was used in Sections \ref{rapid} and \ref{heta} in the definition
of the higher super traces.
Then proceeding as in the proof of
 Proposition \ref{higheta} we can show that
 the higher eta invariant
$$
\widetilde{\eta}_{<e>} =
{2 \over \sqrt{\pi}} \R\,\int_0^{+\infty}
STR_{Cl(1)} ( \, {d \BB_s \over ds } e^{-\BB^2_s}\,)\, ds \in \Omega_{ab}^\ast(G)
$$ is well defined where $\Omega_{ab}^\ast(G)$ denotes the quotient of
$\Omega^\ast(G)$ by
the closure of the space generated by the graded commutators. The differential
of $\Omega^\ast(G)$ allows to endow $\Omega_{ab}^\ast(G)$ with the structure of
a complex whose associated homology will denoted by $\widehat{H}_\ast( \Omega_{ab}^\ast(G)) $.
\medskip

Now we may state the higher local index theorem whose proof
proceeds along the line of the proofs of Theorem 2 of \cite{Go-Lo},  Theorem
13.6 of \cite{LPMEMOIRS} and of Theorem \ref{localindex}
\begin{theorem}
$$
\lim_{s\rightarrow 0^+} \R {}^b STR_{<e>} \, e^{-\AA_s^2} \,=
\int_Z \phi \,\widehat{A}\,(\nabla^{TZ})\,{\rm ch}(\nabla^{\widehat{V}}) \,{\rm ch}(\nabla^{can}) \, \in
\Omega_{ab}^\ast(G)
$$ where $Z$ denotes the typical fiber of $\pi:P \rightarrow T$ and
the connection $\nabla^{can}$ is defined in the section 6 of \cite{Go-Lo}.

\end{theorem}

\begin{corollary}   Let $\Phi$ be a closed graded trace on
  $\Omega_{ab}^\ast(G)$ concentrated in the trivial conjugacy class. Then
  there is a current $\omega_\Phi$ on $M$ such that the following formula holds:
  $$
  < \int_Z \phi \,\widehat{A}\,(\nabla^{TZ})\,{\rm ch}(\nabla^{\widehat{V}}) \,
{\rm ch}(\nabla^{can})\, ; \, \Phi >=
  <\widehat{A}(T\cal{F})\,{\rm ch}(\nabla^{V})\, \omega_\Phi>
  $$
  \end{corollary}

Recall we assume that the quotient $M=P/G$ is a
smooth compact manifold with boundary, it inherits a
foliation $\cal{F}$ whose leaves are the image by the map
$P\rightarrow P/G$ of the fibers of the fibration
$\pi: P \rightarrow G^{(0)}$. If we assume that $G^{(0)}$ admits
a $G-$invariant riemannian metric $g$ then we may state for the $\Ti(G)-$linear operator $D^+$
a result completely analogous to Theorem \ref{decomposition} so that
the index class $\Ind D^+$ is well defined in $K_0(\Ti(G))$.
The higher index $\Ind_\Phi (D^+)$ associated to a closed
graded trace on $\Omega^* (G)$ is defined as in
Definition \ref{higherindexdefinition}.
In this general case,
 proceeding as in Section \ref{isomcase}, we obtain the following analogue
of Theorem \ref{higherindex}
\medskip
\begin{theorem} \label{generalhigherindex} Assume that $G^{(0)}$ admits
a $G-$invariant riemannian metric $g$, then:
\begin{itemize} \item By universality there is a natural morphism of algebras
$$
j: \Omega^\ast(\Ti(G)) \rightarrow \Omega_{ab}^\ast(G).
$$
\item The following formula holds in the
 homology groups $\widehat{H}_\ast( \Omega_{ab}^\ast(G)) :$
$$
j_{<e>}({\rm Ch}\, \Ind D^+)\,=\,
\int_Z \phi \,\widehat{A}\,(\nabla^{TZ}) \,{\rm ch}(\nabla^{\widehat{V}})\,{\rm ch}(\nabla^{can}) -
{1\over 2}  \widetilde{\eta}_{<e>} \in \widehat{H}_\ast( \Omega_{ab}^\ast(G)).
$$
\item Let $\Phi$ denote a closed graded trace on $\Omega^\ast(G)$.
 Then
$$
\Ind_\Phi (D^+) =
<\widehat{A}\,(T\cal{F}) \,{\rm ch}(\nabla^{V})\,,\,
\omega_\Phi> -<{1\over 2}  \widetilde{\eta}_{<e>} \,;\, \Phi>.
$$
\end{itemize}
\end{theorem}

In the general case (i.e. when $G^{(0)}$ does not admit any
a $G-$invariant riemannian metric) we define
${\rm Ch}_{<e>}\, \Ind D^+)\in  \Omega_{ab}^\ast(G)$
 exactly as in
 Definition \ref{chern-super}.

Then, we
have the following Theorem whose proof is completely analogous to the one of Theorem \ref{higherindexIII}.

\n\begin{theorem}
\item {1)}
The following formula holds in the
homology $\widehat{H}_\ast(\Omega^\ast_{ab}(G))$ (see Section \ref{heta1}): $$
{\rm Ch}_{<e>}\, (\Ind D^+)\,=\, \int_Z \phi
\,\widehat{A}\,(\nabla^{TZ})\,{\rm ch}(\nabla^{\widehat{V}})
\,{\rm ch}(\nabla^{can}) - {1\over 2}  \widetilde{\eta}_{<e>} \in
\widehat{H}_\ast(\Omega^\ast_{ab}(G))
$$
\item{2).}
Let $\Phi$ a closed graded trace on $\Omega^\ast_{ab}(G)$,
 concentrated on the trivial conjugacy class.
Let $\widetilde{\Phi}$ be its extension to $\Psi^{-\infty,\delta}_{\Omega^\ast_{ab}(G)}
.$
Then one has:
$$\Ind_\Phi (D^+):=
<{\rm Ch} (\Ind D^+)\,;\,\widetilde{\Phi} > =
<\widehat{A}(T\cal{F})\,{\rm ch}(\nabla^{V}) \,,\,
\omega_\Phi>- <{1\over 2}  \widetilde{\eta}_{<e>}\,;\,
\Phi >
$$
\end{theorem}

\section{{\bf Applications to foliations}}\label{applfoliation}

\m
We give an application of the previous theorem to
a class of foliations.
Let $(Y, \mathcal{F})$ be a foliated  manifold with boundary
where the leaves of $ \mathcal{F} $ are transversal to the boundary $ \partial Y $
 and of
dimension $2k$. We assume that the holonomy groupoid
of $(Y, \mathcal{F})$ is Hausdorff, that
$ ^{b}T\cal{F} $ carries a riemannian $b-$metric
 and
that $ \mathcal{F} $ admits a leafwise spin structure, we shall denote
by $S^{\mathcal{F}}$ the spinor bundle. Moreover, we consider $V\rightarrow Y$
a complex hermitian vector bundle endowed with a hermitian connection.
These data fix a $b-$Dirac type leafwise operator $D$ acting on the
sections of $V\otimes S^{\mathcal{F}}$. Let $H$ denote the holonomy groupoid
of $(Y, \mathcal{F})$ and let
$s:H\rightarrow Y$ denote the source map. By lifting $D$ to $H$ we
 get a $H-$invariant
Dirac type operator on $H$ still denoted $D$. We shall denote by
$D_0$ the boundary operator on $\partial H$ of $D$ and assume that there exists
$\epsilon > 0$ such that for any $y \in \partial Y$, the $L^2-$spectrum
of $D_0$ acting on $L^2(\partial s^{-1}(y)\,,\, (V\otimes S^{\mathcal{F}})_{|
\partial s^{-1}(y)}\,)$ does not meet $]-\epsilon, \epsilon[$.
Then one may associate to $D$ an index class  $\Ind D^+ \in
K_0(C^\ast_r(H))$, following an argument of \cite{Co}, we are
going to use a complete transversal $T$ and  the associate etale
groupoid $H^T_T =G$ so as to interpret (via a Morita equivalence
invariance argument) $\Ind D^+$ as an element of
$K_0(C^\ast_r(G))$ and to apply to it Theorem
\ref{generalhigherindex}. In fact we need to make the following
assumptions:

{\it 1)} There exists a smooth compact submanifold
(possibly with boundary) $T$ of
$Y$ which defines a complete transversal of $(Y, \mathcal{F})$ and
we denote by $G=H^T_T$
the etale groupoid  formed of the elements of $H$ having
their extremities in $T$.

{\it 2)} Set $P=H_T=s^{-1}(T)$. We assume that the restriction
denoted  $\pi=s_{|P}: P \rightarrow T$ is a fibration over $T$.

{\it 3)} Denote by $i$ the inclusion $G \rightarrow H_T=s^{-1}(T)$ and
set
$L =P=s^{-1}(T)$ and $\sigma=s_{|P}$. Recall that for each $\gamma \in G$,
$s^{-1}(s(\gamma)) = G_{s(\gamma)}$ is a complete length metric space.
 Then  we define a length function
$l$ on $G$ by setting for any $\gamma \in G$,
$l(\gamma) = d_{s(\gamma)}(\,i(s(\gamma)) , i(\gamma)\,)$. We then assume that
 (\ref{poly-growth}) holds.

 With these assumptions we consider the algebra
 $\Ti(G)$ and the
 space $\Omega^\ast(G)$  as defined previously.

Let $\int_{\cal{F}}$ denote the Haefliger integration map from
$\Omega^\ast(Y)$ to $\Omega^\ast(T)/V$ where $V$ denote the
sub-vector space of $\Omega^\ast(T)$ generated by the differential
forms $\omega -h(\omega),\, h \in G$. Then Theorem
\ref{generalhigherindex} implies the following
\begin{theorem} Let $\tau$ be a holonomy-invariant  closed
transverse current of $(Y, \mathcal{F})$; it induces a closed graded trace
on $\Omega^\ast(G)$
still denoted $\tau$. Denote by $\widetilde{\tau}$
its extension to $\Psi^{-\infty,\delta}_{\Omega^\ast_{ab}(G)}
.$
 We then have: $$ <{\rm Ch}\, \Ind D^+)\,;\, \widetilde{\tau}> = <\int_{\mathcal{F}}
\,\widehat{A}\,(T\mathcal{F}) \,{\rm ch}(\nabla^{V})\, ;\, \tau > -
<{1\over 2}  \widetilde{\eta}_{<e>} \,;\, \tau>. $$
\end{theorem}


\section{{\bf Appendix A: the rapidly decreasing
$b$-calculus}}\label{rapidcalculus}

\medskip
In order to define the rapidly decreasing $b$-calculus we
introduce as before
 an auxiliary
metric $\hat{g}$ on $\wm$ for which $\wm$ and $\Gamma$ ($\Gamma$
viewed as a metric space with respect to the word-metric) become
{\it quasi-isometric}. The metric $\hat{g}$ is simply the lift to
$\wm$ of an ordinary metric on $M$. In the sequel we denote by
$d(\cdot,\cdot)$ the distance function associated to $\hat{g}$.

\noindent Let $\beta_\pi: [\wm\times_\pi \wm,B]
\rightarrow\wm\times_\pi \wm$ be the blow-down map associated to
the fiber-$b$-stretched product.

\medskip We  consider the fibration
$\pi: \wm\times_\pi \wm \rightarrow T$. Let $(U_j)_{1\leq j \leq
l}$ be a finite open cover of $T$ such that for each $j\in
\{1,\ldots,l\}$, the fibration $\pi: \wm\times_\pi \wm \rightarrow
T$ is trivial over $U_j$. For each $ j \in \{1, \ldots,  l\}$ let
us choose a trivialization
 $\chi_j:\, \pi^{-1}(U_j) \sim Z \times Z
\times U_j$ where $Z$ denotes the typical fiber of $\pi: \wm
\rightarrow T$. Let $X_j$ be any
 smooth vector field on $U_j \subset T$ with compact support on $U_j$, we still denote
 by $X_j$ the induced vector field on  $Z \times Z
\times U_j$,
 then $\chi_j^\ast(X_j)$ defines a differential operator acting on
 $C^\infty(\wm\times_\pi \wm\,,\,
\we \boxtimes \we^\star $. Moreover, for any $p,q\in \NN$ the
family of fiberwise operators $D^p(\theta)$ (resp. $D^q(\theta)$)
defines a differential operator $D_l^p( \cdot)$ (resp.
$D_r^q(\cdot)$) acting on any Schwartz kernel $K(z,w) \in
C^\infty(\wm\times_\pi \wm\,,\, \we \boxtimes \we^\star ) $ by the
formula: $D_l^p( \cdot)\cdot K(z,w) = D^p( \pi(z))\circ K(z,w) $
(resp. $ D_r^p( \cdot)\cdot K(z,w) =  K(z,w)\circ D^q( \pi(w))$).

\begin{definition} \label{Op} \item{1.} We shall denote
by $^{b}{\rm Op}\,(\wm\times_\pi \wm )$ the algebra of
differential operators acting on $C^\infty(\wm\times_\pi \wm\,,\,
\we \boxtimes \we^\star ) $ which is generated by all the above
operators $\chi_j^\ast(X_j)$, $D^p(\cdot)$, $D^q(\cdot)$.
\item{2.} We shall denote by $^{b}{\rm Op}\,( [ \wm \times_\pi \wm , B] )$
the span over $C^\infty([ \wm \times_\pi \wm , B] ) $ of
$\beta_\pi^\ast ( ^{b}{\rm Op}\,(\wm\times_\pi \wm )$ where
 $\beta_\pi^\ast$ denotes the lift of the blow-down map
 $\beta_\pi$.
\end{definition}


\begin{definition}\label{rapid}
 Let $\Gamma$ be virtually nilpotent and let
$P\in\Psi^m_{b,\rtimes}(\wm,\we)$. We shall say that $P$ is {\it
rapidly decreasing} outside an $\epsilon$-neighborhood of the
lifted fiber-diagonal in $\wm\times_\pi\wm$ if
 for each $ Q \in ^{b}{\rm Op}\,( [ \wm \times_\pi \wm , B] )$
 and any $q\in\NN$ we can find a constant $C_{Q,q}>0$
such that $\forall p \in [ \wm \times_\pi \wm , B]$
 such that
 $d(z,z^\prime)>\epsilon$ where $(z,z^\prime)= \beta_\pi (p)$
 $$
|Q(K_{P})(p)|(1+d(z,z'))^q < C_{Q,q}.$$

\end{definition}

\begin{definition}
The elements in  $\Psi^*_{b,\rtimes} (\wm,\we)$ which are rapidly
decreasing outside the lifted diagonal form a subalgebra which is,
by definition, the rapidly decreasing algebra
$\Psi^*_{b,\Ti}(\wm,\we)$.
\end{definition}

\begin{definition}

\n
 We denote by $H^\infty_{b, \Ti}(\wm, \we)$ the subset of the
elements $s \in H^\infty_{b, \,{\rm {loc}}}(\wm, \we)$ such that
for any $\phi \in C^\infty_c(\wm)$ and any integers $p,q\in \NN$:
$$ \sup_{\gamma \in \Gamma} \,\bigl(\, (1+ ||\gamma||)^p || \phi
\,(\gamma\cdot s)||_{H^q_b(\wm,\we)} \, \bigr) < +\infty $$ where
$H^q_{b,\, \Gamma}(\wm,\we)$ is defined using $\Gamma-$equivariant
$b-$differential operators of order $\leq q$ acting on
$C^\infty(\wm , \we)$ (nevertheless, notice that the sections
$\phi \,(\gamma\cdot s)$ have compact support in supp$\phi$.
\end{definition}

We shall now enlarge the algebra $\Psi^*_{b,\Ti}(\wm,\we)$ and
define the rapidly decreasing $b-$calculus with bounds $$ \Psi^{*,
\epsilon}_{b,\Ti}(\wm,\we):= \Psi^{*}_{b, \Ti}(\wm,\we) +
\Psi^{-\infty ,\epsilon}_{b,\Ti}(\wm,\we)+ \Psi^{-\infty,
\epsilon}_{\Ti}(\wm,\we)\,.$$

\begin{definition}
\n \item { 1.}
  We
denote by $\Psi^{-\infty, \epsilon}_{\Ti}(\wm,\we)$ the set of
left $\Ti-$linear operators acting on $ H^\infty_{b, \Ti}(\wm,
\we)$ defined by: $$ s\rightarrow \int_{\pi^{-1}(\theta)} K(z,y)
s(y) dVol^b_{\pi^{-1}(\theta)} \,=\,K( s)(z) $$ where the
 Schwartz kernel $K(z,y) \in \rho^\epsilon_{lb}
\rho^\epsilon_{rb} H^\infty_{b, \,{\rm{ loc}}}(\wm\times_\pi \wm,
\we\boxtimes\we^*)$ is  $\Gamma-$invariant and such that  for any
operator $P\in \;$ $ ^{b}{\rm Op}\,(\wm\times_\pi \wm )$
 and any $p\in \NN$: \begin{equation} \label{es}
   \sup_{R > 1} \, R^p || 1_{d(z,y) > R} \,
(\,P(\rho^{-\epsilon}_{lb} \rho^{-\epsilon}_{rb}K(\cdot,\cdot)\,)
||_{L^2_{b}(\wm\times_\pi \wm\,,\,\we\boxtimes\we^*)} \, < +\infty
\end{equation}
\n \item { 2.} Proceeding as in Remark  \ref{rem}, one can define
the space of operators
 $\Psi^{-\infty, \epsilon}_{b, \Ti}(\wm,\we) $ by
considering the doubled space $\cal{D}_{bf} ( [\wm \times \wm,
B])$ and the $\Gamma$-invariant elements of $$ H^\infty_{b, loc} (
\cal{D}_{bf} ( [\wm \times \wm, B]) ; (\we \boxtimes
\we^\star)_{\cal{D}} ) $$ satisfying estimates analogous to
(\ref{es}).
\end{definition}

\section{{\bf Appendix B: $b$-smoothing operators with
differential form coefficients.}}\label{bsmoothsect}

\medskip
Let $\epsilon \in (0,1)$ and let $$\cal O({\rm bf})=\{\,p\in
[\wm\times_\pi \wm,B]\; | \; d(\beta_\pi (p),B)<\epsilon \,\}.$$
 In $\cal O({\rm bf})$ the variables $r=x+x^\prime$
and $\tau= (x-x^\prime)/(x+x^\prime)$ together with the boundary
variables $(y,y')$ (see \cite{Melrose} Ch. 4) can be used. We set
 $\cal{C} ({\rm bf}):=[\wm\times_\pi\wm,B]\setminus \cal O({\rm bf})$
 and we identify it with its image in $\wm\times_\pi\wm$ under the blow-down map
$\beta_\pi$.

We shall now give a precise definion of the 3 spaces:
$$\Psi^{-\infty}_{b,\widehat{\Omega}_{\ast}(T,\Bi_\Gamma)}
(\widehat{M}\,;\, \widehat{E})\,,\quad
\Psi^{-\infty,\delta}_{b,\widehat{\Omega}_{\ast}(T,\Bi_\Gamma)}
(\widehat{M}\,;\, \widehat{E})\,\quad
\Psi^{-\infty,\delta}_{\widehat{\Omega}_{\ast}(T,\Bi_\Gamma)}
(\widehat{M}\,;\, \widehat{E})$$

\begin{definition} \label{12} For any $(k,l) \in \NN \times \NN$,
$\Psi^{-\infty, \delta}_{b,\widehat{\Omega}_{\ast}(T,\Bi_\Gamma)}
(\widehat{M}\,;\, \widehat{E}) (k,l)$ denotes
 the set of left $ \widehat{\Omega}_{\ast}(T,\Bi_\Gamma)$-linear operators:
$$ \widehat{\Omega}_{\ast}(T,\Bi_\Gamma) \otimes_{{\cal
T}^\infty}H^\infty_{b, \Ti}(\wm,\we)\longrightarrow
\widehat{\Omega}_{\ast}(T,\Bi_\Gamma) \otimes_{{\cal T}^\infty}
H^\infty_{b, {\cal T}^\infty} (\wm,\widehat{E}) $$ whose Schwartz
kernel $K$  can be written in the form $$ K(z,w)= \sum_{g_1,\cdot
\cdot \cdot , g_l\in \Gamma}dg_1\cdot \cdot \cdot dg_l \,
K_{g_1,\cdot \cdot \cdot , g_l}(z,w) $$
 where for any $(z,w) \in  \widehat{ M} \times  \widehat{ M} $,
$K_{g_1,\cdot \cdot \cdot , g_l}(z,w)
\in
\bigwedge^k(T^\ast_{\pi(z)}T) \otimes Hom (\widehat{E}_w\,,\,
\widehat{E}_z)$, vanishes for $\pi(z)(g_1\cdot \cdot \cdot
g_l)^{-1} \not= \pi(w)$ and is such that $$R^*_{((g_1 \cdot \cdot
\cdot  g_l),e)} K_{g_1,\cdot \cdot \cdot ,
g_l}\,\,:=\,\,K_{g_1,\cdot \cdot \cdot , g_l}(\cdot
g_1\cdot\cdot\cdot g_l,\cdot)$$ defines an element in
$\Psi^{-\infty, \delta}_{eb,\pi}(\wm,\we)$, the extended
fiber-$b$-calculus with bounds. We require these kernels to
satisfy the following estimates:

\noindent 1] If $ \cal{C} ({\rm bf}):= [\wm\times_\pi
\wm;B]\setminus \cal O({\rm bf})$ then for any fundamental domain
$A$ of $\wm$, for any $P\in ^{b}{\rm Op}\,(\wm\times_\pi \wm)$ and
any integer $p>1$
 $$ \sup_{\cal{C} ({\rm bf}) } \bigl[d(zg_1\cdot \cdot \cdot  g_l, A) + ||g_2||+ \cdot \cdot
\cdot + ||g_{l}|| + d(w, Ag_1^{-1}) \bigr]^p |P_{z,w}
R^*_{((g_1,\cdot \cdot \cdot , g_l),e)} K_{g_1,\cdot \cdot \cdot ,
g_l}| $$ is finite.

\noindent 2] For any $P \in {\rm Op}\,( \partial \wm\times_\pi
\partial \wm )$ and any integers $p,q_1,q_2\in \NN$, we have $$
\sup_{ \cal O ({\rm bf})} \bigl[d(yg_1\cdot \cdot \cdot  g_l, A) +
||g_2||+ \cdot \cdot \cdot + ||g_{l}|| + d(y^\prime, Ag_1^{-1})
\bigr]^p |P \partial^{q_1}_\tau \partial^{q_2}_r R^*_{((g_1,\cdot
\cdot \cdot , g_l),e)} K_{g_1,\cdot \cdot \cdot , g_l}| $$ is
finite.

\noindent 3] We set: $$\Psi^{-\infty,
\delta}_{b,\widehat{\Omega}_{\ast}(T,\Bi_\Gamma)}
(\widehat{N}\,;\, \widehat{F})= \oplus_{k,l \in \NN}
\Psi^{-\infty, \delta}_{b,\widehat{\Omega}_{\ast}(T,\Bi_\Gamma)}
(\widehat{N}\,;\, \widehat{F})(k,l). $$

\noindent 4] Proceeding as in Section \ref{rapid2}, one defines in a
similar way $\Psi^{-\infty,
\delta}_{\widehat{\Omega}_{\ast}(T,\Bi_\Gamma)} (\widehat{N}\,;\,
\widehat{F})$.
\end{definition}

\noindent 5] The space
$\Psi^{-\infty}_{b,\widehat{\Omega}_{\ast}(T,\Bi_\Gamma)}
(\widehat{N}\,;\, \widehat{F})$ is defined as in 1], 2] above but
with $C^\infty_{\Ti}(\wm,\we)$ appearing instead of $
H^\infty_{b,\Ti}(\wm,\we)$ and with the kernels $R^*_{((g_1 \cdot
\cdot \cdot g_l),e)}K_{g_1,\dots,g_l}$ in
$\Psi^{-\infty}_{eb,\pi}(\wm\we)$.

\section{{\bf Appendix C: a proof of theorem \ref{higherindex}}}\label{appendixc}

In this section we shall give a proof of Theorem
\ref{higherindex}. The proof is quite parallel to the proof of
Theorem 14.1 of \cite{LPMEMOIRS}
 so we shall only sketch the main lines of the proof. Notice that
 since \cite{LPMEMOIRS} deals with {\it right-}modules there is in
 \cite{LPMEMOIRS}
 a grading $\Upsilon$ in front of the Lott's connection, whereas here,
 since we deal with left $\Ti-$modules there is no such grading $\Upsilon$.

First of all we need  the following decomposition theorem:

\begin{theorem} \label{decomposition} We can find $\epsilon > 0,$
$ \cal{L}_\infty$ [resp. $\cal{N}_\infty$] a sub-$\Ti-$module projective of
finite rank of $x^\epsilon H^\infty_{b, \Ti}(\wm, \we^+ )$
[resp. $x^\epsilon H^\infty_{b, \Ti}(\wm, \we^- )$] with the following
properties:

\noindent 1) $ \cal{L}_\infty$ is free and $D^+(\cal{L}_\infty )
\subset \cal{N}_\infty$.

\noindent 2) As Frechet spaces: $$\cal{L}_\infty  \oplus^\perp
\cal{L}_\infty^\perp = x^\epsilon H^\infty_{b, \Ti}(\wm, \we^+
),\;\;\;\;  \cal{N}_\infty \oplus D^+(\cal{L}_\infty^\perp) =
x^\epsilon H^\infty_{b, \Ti}(\wm, \we^- ). $$

\noindent 3) The orthogonal projection $\Pi_{\cal{L}_\infty}$ of $
H^\infty_{b, \Ti}(\wm, \we^+ )$ onto $\cal{L}_\infty $ and the
projection $\Pi_{\cal{N}_\infty}$ of $ H^\infty_{b, \Ti}(\wm,
\we^- )$ onto $\cal{N}_\infty$ along $D^+(\cal{L}_\infty^\perp)$
are operators in $\Psi^{-\infty, \epsilon}_{\Ti}(\wm,\we)$.

\noindent 4) As Banach spaces $$(\,C^0(T)\rtimes_r \Gamma
\otimes_{\Ti} \cal{L}_\infty\, ) \oplus
(\,\overline{C^0(T)\rtimes_r \Gamma \otimes_{\Ti}
{\cal{L}_\infty}^\perp}\,)=L^2_{b,C^0(T)\rtimes_r \Gamma}(\wm,\we)
$$
 $$
(\, C^0(T)\rtimes_r \Gamma \otimes_{\Ti} \cal{N}_\infty \,) \oplus
(\, \overline{C^0(T)\rtimes_r \Gamma \otimes_{\Ti}
D^+(\,{\cal{L}_\infty}^\perp\,) }\,) =H^{-1}_{b,C^0(T)\rtimes_r
\Gamma}(\wm,\we). $$

\noindent 5) The operator $$ D^+: {\cal{L}_\infty}^\perp
\rightarrow D^+({\cal{L}_\infty}^\perp) $$ is invertible for the
Frechet topologies; the induced operator $$ D^+: \overline{
C^0(T)\rtimes_r \Gamma \otimes_{\Ti} {\cal{L}_\infty}^\perp}
\rightarrow \overline{C^0(T)\rtimes_r \Gamma \otimes_{\Ti}
D^+(\,{\cal{L}_\infty}^\perp\,) } $$ is invertible

\noindent 6) The operator $(D^{+})^{-1} \circ (
\Id-P_{\cal{N}_\infty}) $ belongs to $\Psi^{-1,
\epsilon}_{b,\Ti}$.
\end{theorem}

As a consequence of the Theorem we immediately obtain: $$ \Ind D^+
= [{\cal{L}_\infty}]- [{\cal{N}_\infty}] \in K_0(\Ti) \simeq
K_0(C^0(T)\rtimes_r \Gamma ), $$
as already stated in subsection \ref{isomcase2}.

 With the notations of Theorem
\ref{decomposition}, we set: $$ \cal{H}^+=
H^\infty_{b,\Ti}(\wm,\we^+) \oplus \cal{N}_\infty=
\cal{L}_\infty^\perp \oplus \cal{L}_\infty \oplus \cal{N}_\infty
$$ $$\cal{H}^-= H^\infty_{b,\Ti}(\wm,\we^-) \oplus \cal{L}_\infty
=D^+(\cal{L}_\infty^\perp ) \oplus \cal{N}_\infty \oplus
\cal{L}_\infty $$ $$\cal{H}=\cal{H}^+ \oplus\cal{H}^-. $$ Then we
define, for each real $\alpha >0$, the operator
$\R^+_\alpha:\cal{H}^+ \rightarrow \cal{H}^-$ by: $$ \R^+_\alpha
(f \oplus n) = (D^+f + \alpha n ) \oplus \alpha
\Pi_{\cal{L}_\infty}f $$ where $f \in \cal{L}_\infty^\perp \oplus
\cal{L}_\infty$ and $n \in \cal{N}_\infty$.

Next we define the operator $\R^-_\alpha:\cal{H}^- \rightarrow
\cal{H}^+$ by: $$ \R^-_\alpha(g \oplus l) = (D^-g \oplus \alpha l)
\oplus \alpha \Pi_{\cal{N}_\infty} g. $$ where $g\in
H^\infty_{b,\Ti}(\wm,\we^-)$ and $l \in \cal{L}_\infty$. Finally
we define: $$ \R_\alpha=\begin{pmatrix} 0 & \R_\alpha^- \cr
\R_\alpha^+  & 0 \cr
\end{pmatrix}
$$

Now we set $\cal{F}_\infty=\cal{L}_\infty \oplus \cal{N}_\infty$
and we define a $\Ti-$connection $\nabla_{\cal{F}_\infty}$ on
$\cal{F}_\infty$ by compressing $\AA_1-D+{1\over 4} c(\tau)$
 by $\Pi_{\cal{L}_\infty} \oplus \Pi_{\cal{N}_\infty} =\Pi_{\cal{F}_\infty}$. Then we
observe that for each real $u>0$ $$ (\AA_u-u D\oplus
\nabla_{\cal{F}_\infty}) + u\R_\alpha $$ defines a superconnection
operator on $\cal{H}$ and we set: $$ {}^b {\rm
ch}_{u,\alpha}(\we_{\Ti})= \R ({}^b STR_{<e>}\,
e^{-((\AA_u-D\oplus \nabla_{\cal{F}_\infty}) + u\R_\alpha)^2}). $$
Clearly we have:
\begin{equation} \label{chern}
{}^b {\rm ch}_{u,0}(\we_{\Ti})=\R ({}^b STR_{<e>} \,e^{-\AA_u^2})
- \R ({}^b STR_{<e>}\Pi_{\cal{F}_\infty}
e^{-\nabla_{\cal{F}_\infty}^2}\Pi_{\cal{F}_\infty})
\end{equation}
Now using Theorem \ref{decomposition} 3] one sees that $$ {}^b STR
\,\Pi_{\cal{F}_\infty}e^{-\nabla_{\cal{F}_\infty}^2}\Pi_{\cal{F}_\infty}=
 STR_{<e>}\, \Pi_{\cal{F}_\infty}e^{-\nabla_{\cal{F}_\infty}^2}\Pi_{\cal{F}_\infty}.
$$ From Theorem \ref{decomposition} 3] we get $$j_{<e>}({\rm Ch}
\Ind D^+)\,=\, j_{<e>}({\rm Ch}([\cal{L}_\infty]) - {\rm
Ch}([\cal{N}_\infty])). $$
 By using  Proposition \ref{compa} 3) and carrying out a tedious computation we prove
 easily the following lemma
 \begin{lemma} $$
 STR_{<e>}^{alg}  e^{-\nabla_{\cal{F}_\infty}^2}=STR_{<e>}
 \Pi_{\cal{F}_\infty}e^{-\nabla_{\cal{F}_\infty}^2}\Pi_{\cal{F}_\infty}.
 $$\end{lemma}
  The previous lemma and Proposition \ref{compa} thus imply that
$$
 j_{<e>}({\rm Ch}([\cal{L}_\infty]) -
{\rm Ch}([\cal{N}_\infty]))\,=\,\R \,STR_{<e>}\Pi_{\cal{F}_\infty}
e^{-\nabla_{\cal{F}_\infty}^2} \Pi_{\cal{F}_\infty}\in
\overline{H}_\ast(T,\Bi_\Gamma). $$ From the last three equations
we deduce that $$ j_{<e>}({\rm Ch} \Ind D^+)\,=\,\R
\,STR_{<e>}\Pi_{\cal{F}_\infty} e^{-\nabla_{\cal{F}_\infty}^2}
\Pi_{\cal{F}_\infty}\in \overline{H}_\ast(T,\Bi_\Gamma). $$ Hence
from equations (\ref{transgression}),(\ref{chern}) and the
previous one we get $j_{<e>}({\rm Ch} \Ind D^+)= $
\begin{equation} \label{last}
\int_Z \phi \,\widehat{A}\,(\nabla^{TZ}) \,{\rm
ch}(\nabla^{\widehat{V}}) \,{\rm ch}(\nabla^{can}) - {1\over 2}
\int_0^u \widetilde{\eta}_{<e>}(s) ds- d\int_0^u\R ({}^b
STR_{<e>}(\,{d\AA_s \over ds} e^{-\AA_s^2}\,))ds \,-\, {}^b {\rm
ch}_{u,0}(\we_{\Ti})
\end{equation}
 We shall prove Theorem \ref{higherindex} by taking the limit
$u\rightarrow +\infty$ of the right handside of the previous
equation and showing that ${}^b {\rm ch}_{u,0}(\we_{\Ti})
\rightarrow 0$ modulo
 $ d \OM$ as  $u\rightarrow +\infty$. Now, as in \cite{LPMEMOIRS} pages 100
 and 101 one can show that
 \begin{equation} \label{step}
 {}^b {\rm ch}_{u,\alpha}(\we_{\Ti})- {}^b {\rm ch}_{u,0}(\we_{\Ti})
 \in  d \OM
 \end{equation}

 If $\alpha$ is large enough then
 $\R^+_\alpha: \cal{H}^+\rightarrow \cal{H}^-$ is
 invertible and its inverse belongs to the $b-\Ti-$calculus with bounds. Thus
 for  $\alpha>0$ large and $u>1$ we consider the following connection
 on $\cal{H}^+$:
 $$
 \nabla_{u,\cal{H}^+} =
 (\R^+_\alpha)^{-1} \circ(\AA_u -u D
 \oplus \nabla_{\cal{F}_\infty})_{|\cal{H}^+}\circ \R^+_\alpha.
$$ Then we define a new connection $ \nabla_{u,\cal{H}}$ on $
\cal{H}$ by setting $$ \nabla_{u,\cal{H}}=
\nabla_{u,\cal{H}^+}\oplus -(\AA_u -u D \oplus
\nabla_{\cal{F}_\infty})_{|\cal{H}^-} $$ The two superconnections
$\nabla_{u,\cal{H}}$ and $\AA_u -u D
 \oplus \nabla_{\cal{F}_\infty}$ are of course homotopic through
 the path of connections parametrized by $\vartheta \in [0,1]$:
 $$
 \vartheta \nabla_{u,\cal{H}} +
 (1-\vartheta) (\AA_u -u D \oplus \nabla_{\cal{F}_\infty}).
 $$ Now, proceeding as in \cite{LPMEMOIRS} pages 102-103,
 one first proves that, modulo $d \OM$:
 $$
  {}^b {\rm ch}_{u,\alpha}(\we_{\Ti})=
  {}^b \R STR \,e^{-( \nabla_{u,\cal{H}} + u\R_\alpha)^2} + B_2(u, \alpha)
$$ Then, still as in  \cite{LPMEMOIRS} pages 102-103, one proves
that $$ \lim_{u\rightarrow +\infty}\, {}^b \R STR \,e^{-(
\nabla_{u,\cal{H}} + u\R_\alpha)^2}=0,\; \lim_{u\rightarrow
+\infty}\,B_2(u, \alpha)=0. $$  Theorem \ref{higherindex} is then
a consequence of the equations (\ref{step}), (\ref{last}) and of
the three previous ones.

\end{document}